\documentclass[12pt,authoryear]{elsarticle}
\journal{}

\makeatletter
\def\ps@pprintTitle{%
 \let\@oddhead\@empty
 \let\@evenhead\@empty
 \def\@oddfoot{\today\hfill\thepage}%
 \let\@evenfoot\@oddfoot}
\makeatother

\usepackage{amsmath,amssymb,amsfonts,amsthm}
\usepackage[bookmarksnumbered,colorlinks=true]{hyperref}
\usepackage[all]{hypcap}

\providecommand{\doi}[1]{\href{https://doi.org/#1}{DOI:#1}}
\usepackage{xurl} 
\renewcommand{\doi}[1]{%
 \href{https://doi.org/#1}{\nolinkurl{DOI:#1}}%
}

\usepackage{mathtools} 
\usepackage{booktabs} 
\usepackage{float} 
\usepackage{geometry} 
\usepackage{placeins}
\usepackage{tikz} 
\usetikzlibrary{arrows.meta}

\numberwithin{equation}{section}
\numberwithin{table}{section}
\numberwithin{figure}{section}

\theoremstyle{plain}
\newtheorem{theorem}{Theorem}[section]
\newtheorem{proposition}[theorem]{Proposition}
\newtheorem{lemma}[theorem]{Lemma}
\newtheorem{corollary}[theorem]{Corollary}

\theoremstyle{definition}
\newtheorem{remark}{Remark}[section]

\newcommand{\N}{\mathbb{N}}
\newcommand{\R}{\mathbb{R}}
\newcommand{\C}{\mathbb{C}}
\newcommand{\EE}{\mathsf{E}}
\newcommand{\bb}[1]{\boldsymbol{#1}}
\newcommand{\rd}{\mathrm{d}}
\newcommand{\ii}{\mathrm{i}}
\newcommand{\tr}{\mathrm{tr}}
\newcommand{\leqdef}{\vcentcolon=}

\begin{document}

\begin{frontmatter}

\title{\textbf{Spectral Analysis of Gaussian Integral Operators \\Arising in BHEP Tests}}

\author[a1]{Bruno Ebner}\ead{bruno.ebner@kit.edu}
\author[a2]{Dominic Edelmann}\ead{dominic.edelmann@dkfz-heidelberg.de}
\author[a1]{Norbert Henze}\ead{norbert.henze@kit.edu}
\author[a3]{Fr\'{e}d\'{e}ric Ouimet}\ead{frederic.ouimet2@uqtr.ca}
\author[a4]{Donald Richards}\ead{richards@stat.psu.edu}

\address[a1]{Institute of Stochastics, Karlsruhe Institute of Technology, Karlsruhe, Germany}
\address[a2]{Division of Biostatistics, German Cancer Research Center, Heidelberg, Germany}
\address[a3]{D\'{e}partement de math\'{e}matiques et d'informatique, Universit\'{e} du Qu\'{e}bec \`{a} Trois-Rivi\`{e}res, Trois-Rivi\`{e}res, Canada}
\address[a4]{Department of Statistics, Pennsylvania State University, University Park, USA\vspace{-5mm}}

\begin{abstract}
The Baringhaus--Henze--Epps--Pulley (BHEP) tests for multivariate normality are affine-invariant goodness-of-fit tests based on a Gaussian-weighted $L^2$ distance between empirical and Gaussian characteristic functions. In 1990, Henze and Zirkler expressed the limiting null distribution through the eigenvalues of an integral operator on the standard Gaussian space. In 1997, Henze and Wagner obtained a simpler covariance kernel and raised the problem of calculating the eigenvalues of the resulting operator on a Gaussian-weighted space. Although subsequent work treated the univariate case and numerical approximations in a few low dimensions, the complete all-dimensional spectral problem remained open. This paper determines both complete spectra for every dimension $d \in \N$ and every smoothing parameter $\beta > 0$. The two operators are shown to have the forms $\smash{\mathcal{X}_{\beta,d}^*\mathcal{X}_{\beta,d}}$ and $\smash{\mathcal{X}_{\beta,d}\mathcal{X}_{\beta,d}^*}$ for the same Hilbert--Schmidt operator $\mathcal{X}_{\beta,d}$. Consequently, their nonzero eigenvalues agree, including multiplicities, while the null space of the Henze--Zirkler operator is identified exactly. The Gaussian integral operator in the Henze--Wagner decomposition is diagonalized by Mehler's formula, and rotational symmetry confines the finite-rank correction to the sectors associated with spherical harmonics of degrees $0$, $1$, and $2$. The degree-$1$ and degree-$2$ eigenvalues are characterized by scalar transcendental equations, and the radial eigenvalues by an explicit pole-safe Fredholm determinant. The paper establishes nonnegativity, multiplicities, eigenfunction reconstruction, completeness, the trace identity, and a complete characterization of all exceptional pole cases.
\end{abstract}

\begin{keyword} 
BHEP test, covariance operator, eigenvalue problem, Fredholm determinant, Gaussian kernel, goodness-of-fit testing, Henze--Zirkler test, multivariate normality
\MSC[2020]{Primary: 45C05, 47G10; Secondary: 33C45, 33C55, 47A10, 47A55, 62E20, 62H15}
\end{keyword}

\end{frontmatter}

%
\section{Introduction}\label{sec:intro}
%

Testing multivariate normality is a classical composite goodness-of-fit problem, with procedures based on moments, projections, interpoint distances, empirical distribution functions, and integral transforms. Broad accounts of the subject are given by \citet{Henze2002Review} and \citet{EbnerHenze2020}. Among transform methods, the empirical characteristic function is especially attractive because a characteristic function always exists, uniquely determines the corresponding probability law, and can be compared with a model characteristic function without density estimation or a partition of the sample space. The functional asymptotic foundation of this approach includes the uniform consistency and Gaussian-process convergence theory established by \citet{FeuervergerMureika1977}.

In testing for normality, \citet{EppsPulley1983} introduced an omnibus univariate statistic obtained by integrating the squared modulus of the difference between the empirical and fitted normal characteristic functions against a Gaussian weight. \citet{BaringhausHenze1988} extended the construction to multivariate normality by applying it to standardized residuals, \citet{Csorgo1989} coined the acronym BHEP (for Baringhaus--Henze--Epps--Pulley), and \citet{HenzeZirkler1990} developed the parameterized class considered here. \citet{HenzeWagner1997} subsequently established a Gaussian-process representation of the limiting null distribution and studied contiguous alternatives. The resulting tests are consistent against every fixed nonnormal alternative.

To make explicit the sample standardization, the statistic, and its smoothing parameter, let $\bb{X}_1,\ldots,\bb{X}_n$ be independent and identically distributed random vectors in $\R^d$, and define
\[
\overline{\bb{X}}_n \leqdef \frac{1}{n}\sum_{j=1}^n\bb{X}_j, \qquad S_n \leqdef \frac{1}{n}\sum_{j=1}^n(\bb{X}_j - \overline{\bb{X}}_n)(\bb{X}_j - \overline{\bb{X}}_n)^{\top}.
\]
When $S_n$ is nonsingular, let $S_n^{-1/2}$ denote the symmetric positive definite square root of $S_n^{-1}$, and set $\bb{Y}_j \leqdef S_n^{-1/2}(\bb{X}_j - \overline{\bb{X}}_n)$ for $j = 1,\ldots,n$. If $n \geq d + 1$ and the common distribution of the observations assigns probability zero to every affine hyperplane in $\R^d$ (in particular, if it is absolutely continuous with respect to Lebesgue measure), then $S_n$ is nonsingular almost surely (see, e.g., \citet{EatonPerlman1973}).

Fix $\beta > 0$, let $\mu_{\beta}$ be the centered Gaussian probability measure on $\R^d$ with covariance matrix $\beta^2 I_d$, and set
\begin{equation}\label{eq:introduction.empirical.characteristic.function}
\Psi_n(\bb{t}) \leqdef \frac{1}{n}\sum_{j=1}^n\exp(\ii\bb{t}^{\top}\bb{Y}_j), \qquad \bb{t} \in \R^d,
\end{equation}
where $\ii = \sqrt{-1}$ denotes the imaginary unit.

The BHEP statistic studied here is
\begin{equation}\label{eq:introduction.BHEP.statistic}
T_{n,\beta} \leqdef n\int_{\R^d}\left|\Psi_n(\bb{t}) - \exp\left(-\frac{\|\bb{t}\|^2}{2}\right)\right|^2 \, \rd \mu_{\beta}(\bb{t}).
\end{equation}
Following \citet{HenzeWagner1997}, one may set $T_{n,\beta} = 4n$ when $S_n$ is singular; this convention does not affect the continuous normal null model once the sample size is large enough for nonsingularity to hold almost surely. When $S_n$ is nonsingular, Gaussian integration gives the exact computable representation
\[
T_{n,\beta} = \frac{1}{n}\sum_{j,k=1}^n\exp\left(-\frac{\beta^2}{2}\|\bb{Y}_j - \bb{Y}_k\|^2\right) - \frac{2}{(1 + \beta^2)^{d/2}} \sum_{j=1}^n \exp\left\{-\frac{\beta^2\|\bb{Y}_j\|^2}{2(1 + \beta^2)}\right\} + n(1 + 2 \beta^2)^{-d/2}.
\]
This representation avoids numerical integration, and the use of the sample mean and covariance makes the statistic invariant under every nonsingular affine transformation of the observations. The same calculation shows that $T_{n,\beta}$ is $n$ times the biased squared maximum mean discrepancy between the empirical law of the standardized residuals and the standard normal distribution under the Gaussian kernel
\[
k_{\beta}(\bb{x},\bb{y}) \leqdef \exp\left(-\frac{\beta^2}{2}\|\bb{x} - \bb{y}\|^2\right);
\]
see, e.g., \citet{Rustamov2021}.

The parameter $\beta$ determines the spread of the Gaussian frequency weight in \eqref{eq:introduction.BHEP.statistic} and can substantially affect finite-sample power. In the equivalent Gaussian kernel density formulation, the bandwidth is $1/(\beta\sqrt{2})$. Thus the boundary regimes $\beta \downarrow 0$ and $\beta \to \infty$ correspond, respectively, to bandwidths tending to infinity and zero. The limit $\beta \downarrow 0$ is the extreme-smoothing regime studied by \citet{Henze1997Extreme}; see Remark~\ref{rem:extreme.smoothing.limit} for details. The selection of $\beta$ for the BHEP test was studied by \citet{Tenreiro2009}. Under suitable regularity and integrability conditions, \citet{BaringhausEbnerHenze2017} obtained Gaussian limiting distributions for general weighted $L^2$ goodness-of-fit statistics under fixed alternatives, providing a basis for power approximation and inference on distance from the null model. Theorems~\ref{thm:complete.spectrum} and~\ref{thm:Henze.Zirkler.spectrum} determine the spectra for every $\beta > 0$. Their application in Section~\ref{sec:limiting.distribution} relies on the null convergence in \eqref{eq:background.limiting.distribution}, which takes $n \to \infty$ with $\beta$ fixed. This fixed-$\beta$ convergence therefore does not directly cover the commonly implemented Henze--Zirkler choice $\beta_n = 2^{-1/2}\{(2d + 1)n/4\}^{1/(d + 4)}$, which belongs to a different asymptotic regime; see \citet{EbnerHenze2020}.

Under the multivariate normality hypothesis $H_0$ that the common distribution of the observations is $N_d(\bb{\mu},\Sigma)$ for some $\bb{\mu} \in \R^d$ and some positive definite matrix $\Sigma$, \citet[Theorems~2.1 and~2.2]{HenzeWagner1997} obtained a centered real Gaussian process $Z$ with covariance kernel
\begin{equation}\label{eq:definitions.BHEP.kernel}
K(\bb{s},\bb{t}) \leqdef \exp\left(-\frac{\|\bb{s} - \bb{t}\|^2}{2}\right) - \left\{1 + \bb{s}^{\top}\bb{t} + \frac{(\bb{s}^{\top}\bb{t})^2}{2}\right\}\exp\left(-\frac{\|\bb{s}\|^2 + \|\bb{t}\|^2}{2}\right).
\end{equation}
The associated covariance operator on $L^2(\mu_{\beta})$ is
\begin{equation}\label{eq:definitions.BHEP.operator}
(A_{\beta,d}f)(\bb{s}) \leqdef \int_{\R^d}K(\bb{s},\bb{t})f(\bb{t}) \, \rd \mu_{\beta}(\bb{t}).
\end{equation}
They proved
\begin{equation}\label{eq:background.limiting.distribution}
T_{n,\beta} \xrightarrow{\mathrm{law}} T_{\beta}(d) \leqdef \int_{\R^d}Z(\bb{t})^2 \, \rd \mu_{\beta}(\bb{t}),
\end{equation}
and the spectral theorem gives
\begin{equation}\label{eq:background.weighted.chi.square}
T_{\beta}(d) \stackrel{\mathcal{D}}{=} \sum_{j=1}^{\infty}\lambda_j(\beta,d)N_j^2,
\end{equation}
where the $N_j$ are independent standard normal random variables and the $\lambda_j(\beta,d)$ are the positive eigenvalues of $A_{\beta,d}$, repeated according to multiplicity.

This is the classical spectral mechanism behind quadratic goodness-of-fit statistics; see, for example, \citet{AndersonDarling1952}. The complete spectrum determines the limiting distribution, and its power sums determine all cumulants through $\kappa_r\{T_{\beta}(d)\} = 2^{r - 1}(r - 1)!\sum_{j=1}^{\infty}\lambda_j(\beta,d)^r$, $r \in \N$. The largest eigenvalue governs the leading exponential rate of the upper tail and enters approximate Bahadur-efficiency calculations, as discussed by \citet{EbnerHenze2023}. Finite truncations of \eqref{eq:background.weighted.chi.square} can be evaluated by standard methods for quadratic forms, such as the inversion procedure of \citet{Imhof1961}; the omitted nonnegative tail after $m$ terms has mean $\sum_{j>m}\lambda_j(\beta,d)$. When exact weights are unavailable, calibration typically relies on simulation, numerical spectral approximation, or moment-based surrogates. In the univariate Epps--Pulley setting, \citet{Henze1990Approximation} used the first four moments for fitted approximations; \citet{HenzeWagner1997} later calculated the first three cumulants of the multivariate BHEP limit and used moment-matched approximations for critical values.

Two exact operator formulations of these weights have been available for several decades. \citet[Theorem~3.1]{HenzeZirkler1990} expressed the limiting null distribution of the same statistic in terms of the eigenvalues of an integral operator $\widetilde{A}_{\beta,d}$ on $L^2(\mu_1)$ with the explicit four-term kernel $h_{\beta,d}^*$ in \eqref{eq:definitions.HZ.kernel}, incorporating the corrections caused by estimating the normal mean vector and covariance matrix:
\[
(\widetilde{A}_{\beta,d}f)(\bb{x}) \leqdef \int_{\R^d}h_{\beta,d}^*(\bb{x},\bb{y})f(\bb{y}) \, \rd \mu_1(\bb{y}).
\]
Those authors stated that a closed form for its eigenvalues appeared to be out of reach. Seven years later, \citet[pp.~13--14]{HenzeWagner1997} obtained the operator $A_{\beta,d}$ with the substantially shorter kernel \eqref{eq:definitions.BHEP.kernel}, but explicitly stated that they had not succeeded in solving the equation $A_{\beta,d}f = \lambda f$. Since the two operators act on different Gaussian Hilbert spaces and their expanded kernels have very different forms, the equality of their positive spectra is not apparent from the original formulations.

The problem is subtle because the first term of $K$ has a simple Hermite spectrum, whereas the correction is not a sum of eigenprojections for the Gaussian measure $\mu_{\beta}$. The correction nevertheless has finite rank and respects rotations, so its effect can be isolated in three angular sectors.

There has been important recent progress on this spectral problem. In the univariate case, \citet{EbnerHenze2023} derived a Fredholm-determinant-based equation and used a numerically stable reciprocal-root procedure to tabulate the leading Henze--Wagner eigenvalues for several values of the smoothing parameter $\beta$, and explicitly left the multivariate extension open. The final assertion of their Theorem~2.2 excludes coincidences with the unperturbed Gaussian eigenvalues. Corollary~\ref{cor:exceptional.coincidences} below corrects that assertion by establishing the existence of exceptional pole coincidences, even when $d = 1$. \citet{EbnerJimenezGameroMilosevic2025} subsequently proved convergence of a general Rayleigh--Ritz approximation and applied it to leading BHEP eigenvalues and cumulants in dimensions $1$, $2$, and $3$.

Very recently, \citet{GkoumasPapadatosTrevezas2026} independently gave a detailed spectral analysis of the fully standardized univariate Epps–Pulley/BHEP problem. In particular, they established simplicity and strict parity alternation of the positive eigenvalues and obtained a complete classification of the exceptional even-pole coincidences. Their analysis is restricted to the univariate case. These works provided precise univariate and low-dimensional information, but not a complete exact characterization of the spectrum in arbitrary dimension, including multiplicities, eigenfunction reconstruction, and the relationship between the Henze–Zirkler and Henze–Wagner operators.

This paper resolves both spectral problems for every $d \in \N$ and every $\beta > 0$. For $\bb{x},\bb{t} \in \R^d$, define the real feature $\Phi_{\bb{t}}(\bb{x})$ and its centered version $\zeta(\bb{x},\bb{t})$ by
\begin{align}
\Phi_{\bb{t}}(\bb{x}) &\leqdef \cos(\bb{t}^{\top}\bb{x}) + \sin(\bb{t}^{\top}\bb{x}),\nonumber\\
\zeta(\bb{x},\bb{t}) &\leqdef \Phi_{\bb{t}}(\bb{x}) - \exp\left(-\frac{\|\bb{t}\|^2}{2}\right)\left[1 + \bb{t}^{\top}\bb{x} - \frac{1}{2}\{(\bb{t}^{\top}\bb{x})^2 - \|\bb{t}\|^2\}\right].\label{eq:definitions.centered.feature}
\end{align}
The common feature operator $\mathcal{X}_{\beta,d}:L^2(\mu_1)\to L^2(\mu_{\beta})$ is
\[
(\mathcal{X}_{\beta,d}f)(\bb{t}) \leqdef \int_{\R^d}\zeta(\bb{x},\bb{t})f(\bb{x}) \, \rd \mu_1(\bb{x}).
\]
The two historically different kernels admit the compact Gram representations
\begin{equation}\label{eq:introduction.kernel.Gram.representations}
h_{\beta,d}^*(\bb{x},\bb{y}) = \int_{\R^d}\zeta(\bb{x},\bb{t})\zeta(\bb{y},\bb{t}) \, \rd \mu_{\beta}(\bb{t}), \qquad K(\bb{s},\bb{t}) = \int_{\R^d}\zeta(\bb{x},\bb{s})\zeta(\bb{x},\bb{t}) \, \rd \mu_1(\bb{x}),
\end{equation}
and consequently
\begin{equation}\label{eq:introduction.operator.factorizations}
\widetilde{A}_{\beta,d} = \mathcal{X}_{\beta,d}^*\mathcal{X}_{\beta,d}, \qquad A_{\beta,d} = \mathcal{X}_{\beta,d}\mathcal{X}_{\beta,d}^*.
\end{equation}
These identities reveal the common structure hidden by the expanded kernels. Theorem~\ref{thm:Henze.Zirkler.spectrum}\textrm{(i)--(iii)} shows, in particular, that their positive eigenvalues agree, including multiplicities, and identifies the null space of $\widetilde{A}_{\beta,d}$ exactly as the space of polynomials of total degree at most two, generated by normalization together with the score directions for the normal mean and covariance parameters. By contrast, Theorem~\ref{thm:complete.spectrum}\textrm{(v)} shows that $A_{\beta,d}$ is injective. Theorems~\ref{thm:complete.spectrum} and~\ref{thm:Henze.Zirkler.spectrum} determine the complete spectra of both operators, including multiplicities, and Proposition~\ref{prop:eigenfunction.reconstruction} gives reconstruction formulas for the corresponding eigenfunctions. Together, these results determine every weight in \eqref{eq:background.weighted.chi.square}.

The rest of the paper is organized as follows. Section~\ref{sec:definitions} provides additional necessary definitions and notation. Section~\ref{sec:results} states the complete spectra of both operators. Section~\ref{sec:reduction} diagonalizes the Gaussian part and reduces the correction to explicit sequence-space blocks. Section~\ref{sec:proofs} contains the proofs of all results. Section~\ref{sec:limiting.distribution} applies the spectral results to the common limiting null distribution and records the resulting exact series representation. Section~\ref{sec:numerical.evaluations} presents numerical evaluations of the eigenvalues and cumulants.

%
\section{Definitions and notation}\label{sec:definitions}
%

Throughout, ``positive'' and ``negative'' mean strictly greater than and strictly less than zero, respectively, whereas ``nonnegative'' and ``nonpositive'' allow equality. Likewise, ``increasing'' and ``decreasing'' mean strictly increasing and strictly decreasing, whereas ``nondecreasing'' and ``nonincreasing'' allow equality. A self-adjoint operator $C$ on a real Hilbert space $H$ is called nonnegative when $\langle Cf,f\rangle_H \geq 0$ for every $f \in H$; the word positive is not used as a synonym for nonnegative in statements about operators.

Write $\N \leqdef \{1,2,\ldots\}$ and $\N_0 \leqdef \{0,1,2,\ldots\}$. Finite-dimensional vectors are written in bold; $\|\cdot\|$ is the Euclidean norm, $\bb{x}^{\top}$ is the transpose of $\bb{x}$, and $I_d$ is the $d \times d$ identity matrix. Let $\delta_{jk}$ denote the Kronecker delta, and let $\ell^2(\N_0)$ be the real Hilbert space of square-summable sequences indexed by $\N_0$.

Let $d \in \N$ and let $\beta > 0$. The Gaussian measure $\mu_{\beta}$ introduced in Section~\ref{sec:intro} has density
\[
\varphi_{\beta}(\bb{t}) \leqdef (2 \pi \beta^2)^{-d/2}\exp\left(-\frac{\|\bb{t}\|^2}{2 \beta^2}\right), \qquad \bb{t} \in \R^d.
\]
The real Hilbert space $L^2(\mu_{\beta})$ is equipped with the inner product and norm
\[
\langle f,g\rangle_{\beta} \leqdef \int_{\R^d}f(\bb{t})g(\bb{t})\varphi_{\beta}(\bb{t}) \, \rd \bb{t}, \qquad \|f\|_{\beta} \leqdef \langle f,f\rangle_{\beta}^{1/2}.
\]
If $H$ is a real Hilbert space and $v_1,v_2 \in H$, the rank-one operator $v_1 \otimes v_2$ on $H$ is defined by $(v_1 \otimes v_2)f \leqdef \langle f,v_2\rangle_H v_1$. The identity operator is denoted by $I$. For a trace-class operator $C$, $\det\nolimits_F(I - zC)$ denotes its Fredholm determinant. If the nonzero eigenvalues of a nonnegative trace-class operator $C$, repeated according to multiplicity, are $(\xi_j)_{j \in \N}$, then
\[
\det\nolimits_F(I - zC) = \prod_{j=1}^{\infty}(1 - z\xi_j), \qquad z \in \C.
\]
The product converges locally uniformly because $\sum_{j=1}^{\infty}\xi_j < \infty$.

For a multi-index $\bb{\nu} = (\nu_1,\ldots,\nu_d) \in \N_0^d$, write $|\bb{\nu}| \leqdef \nu_1 + \cdots + \nu_d$, $\bb{\nu}! \leqdef \nu_1!\cdots \nu_d!$, and $\bb{x}^{\bb{\nu}} \leqdef x_1^{\nu_1}\cdots x_d^{\nu_d}$. The rising factorial is
\[
(a)_0 \leqdef 1, \qquad (a)_k \leqdef a(a + 1)\cdots(a + k - 1), \qquad k \in \N.
\]
For $0 < \rho < 1$, the infinite Pochhammer symbol with base $\rho$ is
\[
(a;\rho)_{\infty} \leqdef \prod_{k=0}^{\infty}(1 - a\rho^k).
\]
The physicists' Hermite polynomials and the generalized Laguerre polynomials are defined by the Rodrigues formulas
\[
H_n(x) \leqdef (-1)^n e^{x^2}\frac{\rd ^n}{\rd x^n}e^{-x^2}, \qquad L_k^{(\alpha)}(u) \leqdef \frac{u^{-\alpha}e^u}{k!}\frac{\rd ^k}{\rd u^k}\{e^{-u}u^{k+\alpha}\}, \qquad \alpha > -1.
\]
The normalizations used below, as in \citet[Table~18.3.1]{KoornwinderEtAl2010}, are
\[
\int_{\R}H_m(x)H_n(x)e^{-x^2} \, \rd x = \sqrt{\pi}\,2^n n!\,\delta_{mn}, \qquad \int_0^{\infty}L_j^{(\alpha)}(u)L_k^{(\alpha)}(u)u^{\alpha}e^{-u} \, \rd u = \frac{\Gamma(k + \alpha + 1)}{k!}\delta_{jk}.
\]

For $d \geq 2$, let $\sigma_{d-1}$ be the surface area of the unit sphere $\mathbb{S}^{d-1}$, and let $\mathcal{H}_{d,\ell}$ be the space of restrictions to $\mathbb{S}^{d-1}$ of real homogeneous harmonic polynomials of degree $\ell$. Its elements are called spherical harmonics of degree $\ell$. Its dimension is
\[
h_{d,\ell} \leqdef \binom{d + \ell - 1}{\ell} - \binom{d + \ell - 3}{\ell - 2}, \qquad \ell \in \N_0,
\]
where the second binomial coefficient is understood to be zero when $\ell < 2$. For $d = 1$, set $h_{1,0} = h_{1,1} = 1$ and $h_{1,\ell} = 0$ for $\ell \geq 2$. Thus the two one-dimensional sectors are the even and odd subspaces.

The parameters used throughout the spectral calculation are
\[
\Omega \leqdef \sqrt{1 + 4 \beta^2}, \qquad q \leqdef \frac{\Omega - 1}{\Omega + 1}, \qquad \rho \leqdef q^2, \qquad c_{\beta} \leqdef (1 - q)^d.
\]
They satisfy
\begin{equation}\label{eq:definitions.parameter.identities}
0 < q < 1, \qquad \beta^2 = \frac{q}{(1 - q)^2}, \qquad \Omega = \frac{1 + q}{1 - q}, \qquad \omega \leqdef \frac{\Omega}{2 \beta^2} = \frac{1 - q^2}{2q}.
\end{equation}
For $\ell,k \in \N_0$, define
\begin{equation}\label{eq:definitions.Gaussian.levels.and.weights}
\Lambda_{\ell,k} \leqdef c_{\beta}q^{\ell}\rho^k, \qquad \pi_{\ell,k} \leqdef (1 - \rho)^{d/2+\ell}\frac{(\frac{d}{2} + \ell)_k}{k!}\rho^k.
\end{equation}
The binomial series gives $\sum_{k=0}^{\infty}\pi_{\ell,k} = 1$.

For $a > 0$, $j \in \{0,1,2\}$, and $x \in \C \setminus (\{0\} \cup \{\rho^m : m \in \N_0\})$, set
\[
\mathcal{Q}_{a,j}(x) \leqdef \sum_{k=0}^{\infty}\frac{k^j(a)_k}{k!}\frac{\rho^k}{\rho^k - x},
\]
where $k^0 \leqdef 1$, including when $k = 0$. The series converges locally uniformly on its stated domain. The point $x = 0$ is excluded because the defining series diverges there.

For the expanded form of the Henze--Zirkler kernel, let
\[
\tau_{\beta,1} \leqdef \frac{\beta^2}{1 + \beta^2}, \qquad \tau_{\beta,2} \leqdef \frac{\beta^2}{1 + 2 \beta^2}, \qquad \gamma_{\beta} \leqdef \frac{\tau_{\beta,1}}{2} = \frac{\beta^2}{2(1 + \beta^2)}.
\]
The first identity in \eqref{eq:introduction.kernel.Gram.representations} has the following expanded form, which is the kernel displayed in Theorem~3.1 of \citet{HenzeZirkler1990}, written with the normalized Gaussian measures used here:
\begin{align}\label{eq:definitions.HZ.kernel}
h_{\beta,d}^*(\bb{x},\bb{y})
&\leqdef \exp\left(-\frac{\beta^2}{2}\|\bb{x} - \bb{y}\|^2\right) \notag \\
&\quad - (1 + \beta^2)^{-d/2}\exp(-\gamma_{\beta}\|\bb{x}\|^2)\Bigg[1 + \gamma_{\beta}\Big\{\tau_{\beta,1}(\bb{x}^{\top}\bb{y})^2 - \tau_{\beta,1}\|\bb{x}\|^2 - \|\bb{y}\|^2 + 2\bb{x}^{\top}\bb{y} + d\Big\}\Bigg] \notag \\
&\quad - (1 + \beta^2)^{-d/2}\exp(-\gamma_{\beta}\|\bb{y}\|^2)\Bigg[1 + \gamma_{\beta}\Big\{\tau_{\beta,1}(\bb{x}^{\top}\bb{y})^2 - \|\bb{x}\|^2 - \tau_{\beta,1}\|\bb{y}\|^2 + 2\bb{x}^{\top}\bb{y} + d\Big\}\Bigg] \notag \\
&\quad + (1 + 2 \beta^2)^{-d/2}\Bigg[1 - \frac{\tau_{\beta,2}}{2}\Big\{\|\bb{x}\|^2 + \|\bb{y}\|^2 - 2d - 2\bb{x}^{\top}\bb{y}\Big\} \notag \\
&\hspace{35mm} \Bigg. + \frac{\tau_{\beta,2}^2}{4}\Big[\{\|\bb{x}\|^2 - d\}\{\|\bb{y}\|^2 - d\} + 2\Big\{(\bb{x}^{\top}\bb{y})^2 - \|\bb{x}\|^2 - \|\bb{y}\|^2 + d\Big\}\Big]\Bigg].
\end{align}
Let
\begin{equation}\label{eq:P.leq.2}
\mathcal{P}_{\leq 2} \leqdef \operatorname{span}\left(\{1\} \cup \{x_i : 1 \leq i \leq d\} \cup \{x_ix_j - \delta_{ij} : 1 \leq i \leq j \leq d\}\right) \subseteq L^2(\mu_1).
\end{equation}
Denote the orthogonal projection onto this space by $\Pi_{\leq 2}$ and let $\Pi_{\geq 3} \leqdef I - \Pi_{\leq 2}$. The standard-Gaussian kernel operator used below is
\[
(\widetilde{B}_{\beta,d}f)(\bb{x}) \leqdef \int_{\R^d}\exp\left(-\frac{\beta^2}{2}\|\bb{x} - \bb{y}\|^2\right)f(\bb{y})\varphi_1(\bb{y}) \, \rd \bb{y}.
\]
Its compression to $\mathcal{P}_{\leq 2}^{\perp}$ is the operator $\Pi_{\geq 3}\widetilde{B}_{\beta,d}\Pi_{\geq 3}$ restricted to that subspace, obtained by applying $\widetilde{B}_{\beta,d}$ and then projecting its output back onto $\mathcal{P}_{\leq 2}^{\perp}$.

%
\section{The two spectral problems and main results}\label{sec:results}
%

\subsection{The Henze--Wagner operator}\label{subsec:results.Henze.Wagner}

The problem considered in this paper is to determine all nonzero eigenvalues $\lambda$ in the equation
\[
\int_{\R^d}K(\bb{s},\bb{t})f(\bb{t})\varphi_{\beta}(\bb{t}) \, \rd \bb{t} = \lambda f(\bb{s}), \qquad \bb{s} \in \R^d,
\]
including their multiplicities and corresponding eigenfunctions. The next theorem gives a complete answer.

For $k \in \N_0$, define
\begin{equation}\label{eq:results.radial.y}
y_k \leqdef \frac{dq - 2(1 - \rho)k/q}{\sqrt{2d}}, \qquad \bb{v}_k \leqdef \begin{pmatrix}1\\y_k\end{pmatrix}.
\end{equation}
The radial problem requires a pole-safe formulation because a meromorphic resolvent equation can miss eigenvalues that coincide with its apparent poles. Fredholm determinants provide a standard spectral encoding for trace-class integral operators; see \citet{Bornemann2010} for a detailed account and numerical perspective. Define the \textit{pole-safe radial Fredholm determinant} by
\begin{equation}\label{eq:results.radial.Fredholm.determinant}
\mathcal{D}_0(z) \leqdef (zc_{\beta};\rho)_{\infty}\det\left[I_2 + z\sum_{k=0}^{\infty}\frac{c_{\beta}\rho^k\pi_{0,k}}{1 - zc_{\beta}\rho^k}\bb{v}_k\bb{v}_k^{\top}\right].
\end{equation}
At the points $z = (c_{\beta}\rho^m)^{-1}$, $m \in \N_0$, the right-hand side of \eqref{eq:results.radial.Fredholm.determinant} is understood by analytic continuation. Proposition~\ref{prop:radial.Fredholm.determinant} below proves that all apparent singularities are removable.

%
\begin{theorem}[Complete spectrum of the limiting BHEP
operator]\label{thm:complete.spectrum}
%
%

Let $d \in \N$ and $\beta > 0$.

\emph{\textrm{(i)} Unchanged sectors.} For every $\ell \geq 3$ such that $h_{d,\ell} > 0$, and every $k \in \N_0$,
\[
\lambda_{\ell,k}^{(\mathrm{u})} = c_{\beta}q^{\ell+2k}
\]
is an eigenvalue with angular multiplicity $h_{d,\ell}$.

\emph{\textrm{(ii)} Degree-$1$ sector.} For every $m \in \N_0$, the equation
\begin{equation}\label{eq:results.degree.1.equation}
\mathcal{Q}_{d/2+1,0}(x) = 0
\end{equation}
has a unique root $x_{1,m}$ in $(\rho^{m+1},\rho^m)$. The corresponding eigenvalue is
\[
\lambda_{1,m} = c_{\beta}q x_{1,m},
\]
and it has multiplicity $d$.

\emph{\textrm{(iii)} Degree-$2$ sector.} If $d \geq 2$, then, for every $m \in \N_0$, the equation
\begin{equation}\label{eq:results.degree.2.equation}
\mathcal{Q}_{d/2+2,0}(x) = 0
\end{equation}
has a unique root $x_{2,m}$ in $(\rho^{m+1},\rho^m)$. The corresponding eigenvalue is
\[
\lambda_{2,m} = c_{\beta}q^2x_{2,m},
\]
and it has multiplicity $h_{d,2} = (d - 1)(d + 2)/2$.

\emph{\textrm{(iv)} Radial sector.} The radial eigenvalues are precisely the numbers $\lambda = z^{-1}$ for which $z > 0$ and $\mathcal{D}_0(z) = 0$. Every radial eigenvalue is simple. Away from the unperturbed values $c_{\beta}\rho^m$, the radial eigenvalues are equivalently the numbers $\lambda = c_{\beta}x$, where $x \in (0,1) \setminus \{\rho^m : m \in \N_0\}$ satisfies
\begin{equation}\label{eq:results.radial.nonpole.equation}
\mathcal{Q}_{d/2,0}(x)\mathcal{Q}_{d/2,2}(x) - \mathcal{Q}_{d/2,1}(x)^2 = 0.
\end{equation}
For $m \in \N_0$, the unperturbed value $c_{\beta}\rho^m$ is a radial eigenvalue if and only if
\begin{equation}\label{eq:results.radial.pole.equation}
\sum_{\substack{k \in \N_0 \\ k\neq m}}\frac{(k - m)^2(d/2)_k}{k!}\frac{\rho^k}{\rho^k - \rho^m} = 0.
\end{equation}
If the radial eigenvalues are arranged in decreasing order as $(\lambda_{0,m})_{m \in \N_0}$, then
\begin{equation}\label{eq:results.radial.interlacing}
c_{\beta}\rho^m > \lambda_{0,m} > c_{\beta}\rho^{m+2}, \qquad m \in \N_0.
\end{equation}

\emph{\textrm{(v)} Exhaustiveness.} The multiset union of the eigenvalues in parts \textrm{(i)}--\textrm{(iv)} is the complete nonzero spectrum of $A_{\beta,d}$. The operator is nonnegative, trace class, and injective, so zero is an accumulation point but is not an eigenvalue. If values from different displayed families coincide, their multiplicities are added.
\end{theorem}

The unchanged contribution at a fixed total Hermite degree can be written without referring to angular sectors. This formulation is useful when the limiting distribution in \eqref{eq:background.weighted.chi.square} is assembled numerically.

%
\begin{corollary}[Baseline multiplicities of the geometric
values]\label{cor:baseline.multiplicities}
%
For $N \in \N_0$, the geometric value $c_{\beta}q^N$ has an unchanged spectral contribution of multiplicity
\[
m_N^{(\mathrm{u})} =
\begin{cases}
0, & N \in \{0,1,2\},\\
\displaystyle \binom{N + d - 1}{N} - d, & N \geq 3 \text{ is odd},\\
\displaystyle \binom{N + d - 1}{N} - \binom{d + 1}{2}, & N \geq 4 \text{ is even}.
\end{cases}
\]
Any coincident eigenvalues arising from the radial, degree-$1$, or degree-$2$ equations contribute additional multiplicity.
\end{corollary}

\begin{remark}[The one-dimensional case]\label{rem:one.dimensional.case}
When $d = 1$, only the even and odd sectors exist. The radial determinant in \eqref{eq:results.radial.Fredholm.determinant} gives the even spectrum, and \eqref{eq:results.degree.1.equation} gives the odd spectrum. There is no degree-$2$ traceless sector and no unchanged higher angular sector.
\end{remark}

\subsection{The Henze--Zirkler operator}\label{subsec:results.Henze.Zirkler}

The next theorem solves the eigenvalue problem for the longer kernel in \eqref{eq:definitions.HZ.kernel}. It also explains why the eigenvalues in the two earlier formulations of the limiting distribution are the same.

%
\begin{theorem}[Complete spectrum of the Henze--Zirkler operator]\label{thm:Henze.Zirkler.spectrum}
%
%
Let $d \in \N$ and $\beta > 0$. Recall the definition of $\Pi_{\geq 3}$ just below \eqref{eq:P.leq.2}.

\emph{\textrm{(i)} Exact factorizations.} The operator $\mathcal{X}_{\beta,d}$ is Hilbert--Schmidt, and
\begin{equation}\label{eq:results.HZ.factorizations}
\widetilde{A}_{\beta,d} = \mathcal{X}_{\beta,d}^*\mathcal{X}_{\beta,d} = \Pi_{\geq 3}\widetilde{B}_{\beta,d}\Pi_{\geq 3}, \qquad A_{\beta,d} = \mathcal{X}_{\beta,d}\mathcal{X}_{\beta,d}^*.
\end{equation}
Consequently, $\widetilde{A}_{\beta,d}$ is nonnegative, self-adjoint, and trace class.

\emph{\textrm{(ii)} Complete positive spectrum.} The positive eigenvalues of $\widetilde{A}_{\beta,d}$, including their multiplicities, are exactly the eigenvalues of $A_{\beta,d}$ listed in Theorem~\ref{thm:complete.spectrum}. Concretely, they consist of the unchanged values $c_{\beta}q^{\ell+2k}$ from part \textrm{(i)} of that theorem, the values $c_{\beta}q x_{1,m}$ and, when $d \geq 2$, $c_{\beta}q^2x_{2,m}$ obtained from the roots defined by \eqref{eq:results.degree.1.equation} and \eqref{eq:results.degree.2.equation}, respectively, and the radial values $z^{-1}$, where $z > 0$ and $\mathcal{D}_0(z) = 0$. Equivalently, the non-pole radial values are $c_{\beta}x$, where $x$ satisfies \eqref{eq:results.radial.nonpole.equation}, while every pole satisfying \eqref{eq:results.radial.pole.equation} contributes the value $c_{\beta}\rho^m$. The degree-$2$ family is absent when $d = 1$. All multiplicities and possible coincidences are exactly as stated there.

\emph{\textrm{(iii)} Zero eigenspace.} The null space is
\begin{equation}\label{eq:results.HZ.null.space}
\ker(\widetilde{A}_{\beta,d}) = \mathcal{P}_{\leq 2}, \qquad \dim\ker(\widetilde{A}_{\beta,d}) = 1 + d + \frac{d(d + 1)}{2} = \frac{(d + 1)(d + 2)}{2}.
\end{equation}
Thus zero is an eigenvalue with the displayed finite multiplicity and is also the accumulation point of the positive eigenvalues.

\emph{\textrm{(iv)} Unitary equivalence and eigenfunctions.} The restriction of $\widetilde{A}_{\beta,d}$ to $\mathcal{P}_{\leq 2}^{\perp}$ is unitarily equivalent to $A_{\beta,d}$. More explicitly, if $f \in \mathcal{P}_{\leq 2}^{\perp}$ is a normalized eigenfunction of $\widetilde{A}_{\beta,d}$ with eigenvalue $\lambda > 0$, then
\begin{equation}\label{eq:results.HZ.eigenfunction.forward}
\bb{t} \longmapsto \frac{1}{\sqrt{\lambda}}\int_{\R^d}\zeta(\bb{x},\bb{t})f(\bb{x})\varphi_1(\bb{x}) \, \rd \bb{x}
\end{equation}
is the corresponding normalized eigenfunction of $A_{\beta,d}$. Conversely, if $g \in L^2(\mu_{\beta})$ is a normalized eigenfunction of $A_{\beta,d}$ with eigenvalue $\lambda > 0$, then
\begin{equation}\label{eq:results.HZ.eigenfunction.backward}
\bb{x} \longmapsto \frac{1}{\sqrt{\lambda}}\int_{\R^d}\zeta(\bb{x},\bb{t})g(\bb{t})\varphi_{\beta}(\bb{t}) \, \rd \bb{t}
\end{equation}
is the corresponding normalized eigenfunction of $\widetilde{A}_{\beta,d}$. These two maps are mutually inverse on the corresponding eigenspaces.

\emph{\textrm{(v)} Trace and Fredholm determinant.} One has
\begin{equation}\label{eq:results.HZ.trace}
\tr(\widetilde{A}_{\beta,d}) = \tr(A_{\beta,d}) = 1 - (1 + 2 \beta^2)^{-d/2} - d\beta^2(1 + 2 \beta^2)^{-d/2-1} - \frac{d(d + 2)}{2}\beta^4(1 + 2 \beta^2)^{-d/2-2},
\end{equation}
and
\begin{equation}\label{eq:results.HZ.Fredholm.determinant}
\det\nolimits_F(I - z\widetilde{A}_{\beta,d}) = \det\nolimits_F(I - zA_{\beta,d}), \qquad z \in \C.
\end{equation}
In particular, the weights $\delta_k(\beta)$ in Theorem~3.1 of \citet{HenzeZirkler1990} are precisely the positive eigenvalues in Theorem~\ref{thm:complete.spectrum}.
\end{theorem}

%
\section{Diagonalization and finite-rank reduction}\label{sec:reduction}
%

The first step is to diagonalize the Gaussian integral operator $B_{\beta,d}$ on $L^2(\mu_{\beta})$, defined by
\[
(B_{\beta,d}f)(\bb{s}) \leqdef \int_{\R^d}\exp\left(-\frac{\|\bb{s} - \bb{t}\|^2}{2}\right)f(\bb{t}) \, \rd \mu_{\beta}(\bb{t}).
\]
For $\bb{\nu} \in \N_0^d$, define
\begin{equation}\label{eq:cartesian.eigenfunctions}
e_{\bb{\nu}}(\bb{x}) \leqdef \frac{\Omega^{d/4}}{\{2^{|\bb{\nu}|}\bb{\nu}!\}^{1/2}}\prod_{j=1}^d H_{\nu_j}(\sqrt{\omega}x_j)\exp\left(-\frac{\|\bb{x}\|^2}{\Omega + 1}\right), \qquad \bb{x}\in \R^d.
\end{equation}

%
\begin{proposition}[Cartesian Mercer expansion]\label{prop:Cartesian.Mercer.expansion}
%
The family $(e_{\bb{\nu}})_{\bb{\nu}\in\N_0^d}$ is an orthonormal basis of $L^2(\mu_{\beta})$, and
\begin{equation}\label{eq:reduction.Cartesian.Mercer.eigenvalues}
B_{\beta,d}e_{\bb{\nu}} = c_{\beta}q^{|\bb{\nu}|}e_{\bb{\nu}}, \qquad \bb{\nu} \in \N_0^d.
\end{equation}
Moreover, we have the Mercer expansion
\begin{equation}\label{eq:reduction.Cartesian.Mercer.kernel}
\exp\left(-\frac{\|\bb{s} - \bb{t}\|^2}{2}\right) = \sum_{\bb{\nu}\in\N_0^d}c_{\beta}q^{|\bb{\nu}|}e_{\bb{\nu}}(\bb{s})e_{\bb{\nu}}(\bb{t}),
\end{equation}
where the series converges pointwise and locally uniformly. In particular, $B_{\beta,d}$ is nonnegative and trace class, and $\tr(B_{\beta,d}) = 1$.
\end{proposition}

The one-dimensional version of \eqref{eq:reduction.Cartesian.Mercer.eigenvalues} gives the geometric spectrum found in \citet[Proposition~1]{ShiBelkinYu2009}, and the isotropic $d$-dimensional version follows by tensorization.

Passing from Cartesian Hermite functions to spherical harmonics and generalized Laguerre polynomials separates the angular and radial variables. For the rotational decomposition, fix an orthonormal basis $(Y_{\ell,j})_{1\leq j\leq h_{d,\ell}}$ of $\mathcal{H}_{d,\ell}$ with respect to ordinary surface measure. For $d \geq 2$, $\ell,k \in \N_0$, $1 \leq j \leq h_{d,\ell}$, $r \geq 0$, and $\bb{\theta} \in \mathbb{S}^{d-1}$, define
\begin{align}\label{eq:reduction.spherical.Laguerre.functions}
\psi_{\ell,k,j}(r\bb{\theta})
&\leqdef (2 \pi \beta^2)^{d/4}\left\{\frac{2\omega^{\ell+d/2}k!}{\Gamma(k + \ell + \frac{d}{2})}\right\}^{1/2}r^{\ell}L_k^{(\ell+d/2-1)}(\omega r^2) \exp\left(-\frac{1 - q}{2}r^2\right)Y_{\ell,j}(\bb{\theta}).
\end{align}
For $d = 1$, define $\psi_{0,k,1} \leqdef (-1)^k e_{2k}$ and $\psi_{1,k,1} \leqdef (-1)^k e_{2k+1}$, where the one-dimensional Cartesian eigenfunctions are defined in \eqref{eq:cartesian.eigenfunctions}.

%
\begin{proposition}[Spherical Mercer expansion]\label{prop:spherical.Mercer.expansion}
%
For $d \geq 2$, the functions in \eqref{eq:reduction.spherical.Laguerre.functions}, for all $\ell,k \in \N_0$ and $1 \leq j \leq h_{d,\ell}$, form an orthonormal basis of $L^2(\mu_{\beta})$ and satisfy
\begin{equation}\label{eq:reduction.spherical.Mercer.eigenvalues}
B_{\beta,d}\psi_{\ell,k,j} = \Lambda_{\ell,k}\psi_{\ell,k,j} = c_{\beta}q^{\ell+2k}\psi_{\ell,k,j}.
\end{equation}
For $d = 1$, the same assertion holds with the even and odd families just defined.
\end{proposition}

The finite-rank correction can now be placed exactly within this decomposition. Set
\begin{equation}\label{eq:def.g}
g(\bb{x}) \leqdef \exp(-\|\bb{x}\|^2/2), \qquad \bb{x}\in \R^d.
\end{equation}
If $d \geq 2$, let $(E_j)_{1\leq j\leq h_{d,2}}$ be a Frobenius-orthonormal basis of the real traceless symmetric $d \times d$ matrices, so that $\tr(E_iE_j) = \delta_{ij}$. The spherical harmonic bases of degrees $0$, $1$, and $2$ may and shall be chosen as
\begin{equation}\label{eq:reduction.aligned.harmonics}
Y_{0,1}(\bb{\theta}) = \sigma_{d-1}^{-1/2}, \qquad Y_{1,i}(\bb{\theta}) = \left(\frac{d}{\sigma_{d-1}}\right)^{1/2}\theta_i, \qquad Y_{2,j}(\bb{\theta}) = \left\{\frac{d(d + 2)}{2\sigma_{d-1}}\right\}^{1/2}\bb{\theta}^{\top}E_j\bb{\theta},
\end{equation}
respectively.

%
\begin{proposition}[Exact block reduction]\label{prop:exact.block.reduction}
%
The operator $A_{\beta,d}$ defined in \eqref{eq:definitions.BHEP.operator}, with kernel $K$ defined in \eqref{eq:definitions.BHEP.kernel}, has the finite-rank representation
\begin{equation}\label{eq:reduction.finite.rank.decomposition}
\begin{aligned}
A_{\beta,d}
&= B_{\beta,d} - g \otimes g - \sum_{i=1}^d(x_i g) \otimes (x_i g) \\
&\qquad- \left(\frac{\|\bb{x}\|^2g}{\sqrt{2d}}\right)\otimes\left(\frac{\|\bb{x}\|^2g}{\sqrt{2d}}\right) - \sum_{j=1}^{h_{d,2}}\left(\frac{\bb{x}^{\top}E_j\bb{x}}{\sqrt{2}}g\right)\otimes\left(\frac{\bb{x}^{\top}E_j\bb{x}}{\sqrt{2}}g\right),
\end{aligned}
\end{equation}
where the last sum is absent when $d = 1$. With the basis choices above, the coefficient formulas within the affected angular sectors are as follows. In all four identities, $k \in \N_0$; in the third identity, $1 \leq i,j \leq d$; and in the fourth identity, $d \geq 2$ and $1 \leq i,j \leq h_{d,2}$:
\begin{align}
\langle g,\psi_{0,k,1}\rangle_{\beta} &= \{\Lambda_{0,k}\pi_{0,k}\}^{1/2},\label{eq:reduction.coefficient.radial.first}\\
\left\langle\frac{\|\bb{x}\|^2g}{\sqrt{2d}},\psi_{0,k,1}\right\rangle_{\beta} &= y_k\{\Lambda_{0,k}\pi_{0,k}\}^{1/2},\label{eq:reduction.coefficient.radial.second}\\
\langle x_i g,\psi_{1,k,j}\rangle_{\beta} &= \delta_{ij}\{\Lambda_{1,k}\pi_{1,k}\}^{1/2},\label{eq:reduction.coefficient.linear}\\
\left\langle\frac{\bb{x}^{\top}E_i\bb{x}}{\sqrt{2}}g,\psi_{2,k,j}\right\rangle_{\beta} &= \delta_{ij}\{\Lambda_{2,k}\pi_{2,k}\}^{1/2}.\label{eq:reduction.coefficient.quadratic}
\end{align}
All coefficients of the removed functions $g$, $x_i g$, $\frac{\|\bb{x}\|^2g}{\sqrt{2d}}$, and $\frac{\bb{x}^{\top}E_j\bb{x}}{\sqrt{2}}g$ against basis functions of a different angular degree vanish. Also, \eqref{eq:reduction.coefficient.quadratic} does not appear when $d = 1$ because $h_{1,2} = 0$.

For $\ell \in \{1,2\}$ and $1 \leq j \leq h_{d,\ell}$, set
\[
\mathcal{V}_{\ell,j} \leqdef \overline{\operatorname{span}}\{\psi_{\ell,k,j} : k \in \N_0\}.
\]
Under the unitary identification $\sum_{k=0}^{\infty}a_k\psi_{\ell,k,j} \longleftrightarrow (a_k)_{k \in \N_0}$ between $\mathcal{V}_{\ell,j}$ and $\ell^2(\N_0)$, the restriction of $A_{\beta,d}$ to $\mathcal{V}_{\ell,j}$ is represented by
\begin{equation}\label{eq:reduction.rank.one.matrix}
D_{\ell}^{1/2}P_{\ell}D_{\ell}^{1/2}, \qquad D_{\ell} \leqdef \operatorname{diag}(\Lambda_{\ell,0},\Lambda_{\ell,1},\ldots), \qquad P_{\ell} \leqdef I - u_{\ell} \otimes u_{\ell},
\end{equation}
where $u_{\ell} \leqdef (\sqrt{\pi_{\ell,k}})_{k \in \N_0}$ is a unit sequence in $\ell^2(\N_0)$. The case $\ell = 2$ is present only when $d \geq 2$.

Likewise, set
\begin{equation}\label{eq:def.V.0}
\mathcal{V}_0 \leqdef \overline{\operatorname{span}}\{\psi_{0,k,1} : k \in \N_0\}.
\end{equation}
The subspace $\mathcal{V}_0$ is the closed subspace of radial functions when $d \geq 2$ and the even subspace when $d = 1$. Under the unitary identification $\sum_{k=0}^{\infty}a_k\psi_{0,k,1} \longleftrightarrow (a_k)_{k \in \N_0}$, the restriction of $A_{\beta,d}$ to $\mathcal{V}_0$, denoted $A_{\beta,d}^{(0)}$, is represented by
\begin{equation}\label{eq:reduction.radial.matrix}
A_{\beta,d}^{(0)} = D_0^{1/2}Q_0D_0^{1/2}, \qquad D_0 \leqdef \operatorname{diag}(c_{\beta},c_{\beta}\rho,\ldots), \qquad Q_0 \leqdef I - U_0U_0^*,
\end{equation}
where $U_0:\R^2\to\ell^2(\N_0)$ is defined by
\begin{equation}\label{eq:def.U.0}
(U_0\bb{a})_k \leqdef \sqrt{\pi_{0,k}}\bb{v}_k^{\top}\bb{a}, \qquad k \in \N_0.
\end{equation}
The map $U_0$ is an isometry. On every sector of degree $\ell \geq 3$, $A_{\beta,d}$ agrees with $B_{\beta,d}$.
\end{proposition}

Thus the finite-rank correction affects only spherical harmonic degrees $0$, $1$, and $2$. Each angular copy of degree $1$ or $2$ is a rank-one perturbation of a diagonal operator, whereas the radial sector is a rank-two perturbation.

%

\begin{proposition}[Nonnegativity, trace class, and injectivity]\label{prop:positivity.trace.injectivity}
%
The operator $A_{\beta,d}$ is nonnegative, self-adjoint, and trace class. Moreover,
\begin{equation}\label{eq:reduction.trace.formula}
\tr(A_{\beta,d}) = 1 - (1 + 2 \beta^2)^{-d/2} - d\beta^2(1 + 2 \beta^2)^{-d/2-1} - \frac{d(d + 2)}{2}\beta^4(1 + 2 \beta^2)^{-d/2-2}.
\end{equation}
It is also injective.
\end{proposition}

The product in \eqref{eq:results.radial.Fredholm.determinant} is already adequate for numerical evaluation away from its apparent poles, which occur when a denominator $1 - zc_{\beta}\rho^k$ vanishes. The following expansion removes these denominators entirely, exhibiting $\mathcal{D}_0$ as an entire function and thereby proving that the pole prescription (i.e., the value assigned to $\mathcal{D}_0$ at each apparent pole by analytic continuation) is unambiguous. Let
\[
d_k \leqdef c_{\beta}\rho^k, \qquad p_k(z) \leqdef 1 - zd_k, \qquad P(z) \leqdef \prod_{k=0}^{\infty}p_k(z),
\]
so that $p_k$ vanishes exactly at the $k$th apparent pole $z = d_k^{-1}$. Define $P_k(z) \leqdef \prod_{r\neq k}p_r(z)$ and $P_{i,j}(z) \leqdef \prod_{r\neq i,j}p_r(z)$, which are the products obtained by omitting the $k$th factor and the $i$th and $j$th factors, respectively.

%
\begin{proposition}[Entire radial Fredholm determinant]\label{prop:radial.Fredholm.determinant}
%
The function $\mathcal{D}_0$ defined in \eqref{eq:results.radial.Fredholm.determinant} has the manifestly entire expansion
\begin{equation}\label{eq:reduction.radial.entire.expansion}
\mathcal{D}_0(z) = P(z) + z\sum_{k=0}^{\infty}d_k\pi_{0,k}(1 + y_k^2)P_k(z) + z^2\sum_{0\leq i<j}d_id_j\pi_{0,i}\pi_{0,j}(y_i - y_j)^2P_{i,j}(z).
\end{equation}
All products and sums in \eqref{eq:reduction.radial.entire.expansion} converge locally uniformly on $\C$. Moreover,
\begin{equation}\label{eq:reduction.radial.Fredholm.identity}
\mathcal{D}_0(z) = \det\nolimits_F(I - zA_{\beta,d}^{(0)}), \qquad z \in \C,
\end{equation}
where $A_{\beta,d}^{(0)}$ is the restriction of $A_{\beta,d}$ to the subspace $\mathcal{V}_0$ defined in \eqref{eq:def.V.0}.
\end{proposition}

%
\begin{corollary}[Exceptional coincidences occur]\label{cor:exceptional.coincidences}
%
For every $d \in \N$ and every $m \in \N$, there is at least one $\beta > 0$ for which $c_{\beta}\rho^m$ is a radial eigenvalue of $A_{\beta,d}$. For $m=0$, the value $c_{\beta}$ itself is never a radial eigenvalue.
\end{corollary}

%
\begin{lemma}[One-coordinate compression]\label{lem:one.coordinate.compression}
%
Let $\delta_0 > \delta_1 > \cdots > 0$ with $\delta_k \to 0$, let $D = \operatorname{diag}(\delta_0,\delta_1,\ldots)$ on $\ell^2(\N_0)$, and let $u$ be a unit sequence all of whose coordinates are nonzero. If $P = I - u \otimes u$, then the nonzero spectra of $D^{1/2}PD^{1/2}$ and the compression $PDP$ on $u^{\perp}$ coincide, including multiplicity. All their positive eigenvalues are simple, and there is exactly one eigenvalue in every interval $(\delta_{m+1},\delta_m)$.
\end{lemma}

%
\begin{proposition}[Eigenfunction reconstruction]\label{prop:eigenfunction.reconstruction}
%
The eigenfunctions can be recovered explicitly from the roots in Theorem~\ref{thm:complete.spectrum}. In a degree-$\ell$ sector, $\ell \in \{1,2\}$, the eigenvalue $\lambda_{\ell,m} = c_{\beta}q^{\ell}x_{\ell,m}$ has, for each angular index $j$, an eigenfunction proportional to
\begin{equation}\label{eq:reduction.rank.one.eigenfunction}
\sum_{k=0}^{\infty}\frac{\{\Lambda_{\ell,k}\pi_{\ell,k}\}^{1/2}}{\Lambda_{\ell,k} - \lambda_{\ell,m}}\psi_{\ell,k,j}.
\end{equation}
If $\lambda = c_{\beta}x$ is a non-pole radial eigenvalue, choose a nonzero vector $\bb{\gamma}$ in the null space of
\[
\sum_{k=0}^{\infty}\frac{\pi_{0,k}}{\rho^k - x}\bb{v}_k\bb{v}_k^{\top}.
\]
An associated radial eigenfunction is proportional to
\begin{equation}\label{eq:reduction.radial.nonpole.eigenfunction}
\sum_{k=0}^{\infty}\frac{\{\Lambda_{0,k}\pi_{0,k}\}^{1/2}\bb{v}_k^{\top}\bb{\gamma}}{\Lambda_{0,k} - \lambda}\psi_{0,k,1}.
\end{equation}
At a radial pole eigenvalue $\lambda = d_m$, take $\bb{w}_m \leqdef (y_m,-1)^{\top}$ and set
\begin{equation}\label{eq:reduction.radial.pole.eigenfunction.coefficients}
a_k \leqdef \frac{\{d_k\pi_{0,k}\}^{1/2}\bb{v}_k^{\top}\bb{w}_m}{d_k - d_m}, \quad k \neq m, \qquad a_m \leqdef \frac{\bb{b}_m^{\top}\left[I_2 - \displaystyle\sum_{k\neq m}\frac{\bb{b}_k\bb{b}_k^{\top}}{d_k - d_m}\right]\bb{w}_m}{\|\bb{b}_m\|^2},
\end{equation}
where $\bb{b}_k \leqdef \{d_k\pi_{0,k}\}^{1/2}\bb{v}_k$. Then $\sum_{k=0}^{\infty}a_k\psi_{0,k,1}$ is an eigenfunction. Every series in this proposition converges in $L^2(\mu_{\beta})$ and may be normalized by dividing by its norm. The unchanged eigenfunctions are the $\psi_{\ell,k,j}$ with $\ell \geq 3$.
\end{proposition}

%
\section{Proofs}\label{sec:proofs}
%

\begin{figure}[H]
\centering
\resizebox{0.85\textwidth}{!}{%
\begin{tikzpicture}[
result/.style={draw, fill=white, rounded corners, align=center, text width=3.15cm, minimum height=1.00cm, inner sep=4pt, font=\small, execute at begin node={\hyphenpenalty=10000\relax}},
dependency/.style={-{Stealth[length=2.2mm,width=1.5mm]}, semithick, rounded corners=2pt, shorten <=1pt, shorten >=1pt}
]
\node[result] (cartesian) at (-6.25,9.6) {\textbf{Proposition~\ref{prop:Cartesian.Mercer.expansion}}\\Cartesian Mercer\\expansion};
\node[result] (spherical) at (-6.25,7.2) {\textbf{Proposition~\ref{prop:spherical.Mercer.expansion}}\\Spherical Mercer\\expansion};
\node[result] (block) at (-6.25,4.8) {\textbf{Proposition~\ref{prop:exact.block.reduction}}\\Exact block reduction};
\node[result] (fredholm) at (-6.25,2.4) {\textbf{Proposition~\ref{prop:radial.Fredholm.determinant}}\\Radial Fredholm\\determinant};
\node[result] (positivity) at (-2.0,8.0) {\textbf{Proposition~\ref{prop:positivity.trace.injectivity}}\\Nonnegativity, trace class,\\and injectivity};
\node[result] (compression) at (-2.4,3.65) {\textbf{Lemma~\ref{lem:one.coordinate.compression}}\\One-coordinate\\compression};
\node[result] (complete) at (2.2,4.8) {\textbf{Theorem~\ref{thm:complete.spectrum}}\\Complete BHEP\\spectrum};
\node[result] (reconstruction) at (6.6,4.8) {\textbf{Proposition~\ref{prop:eigenfunction.reconstruction}}\\Eigenfunction\\reconstruction};
\node[result] (baseline) at (-1.4,0.3) {\textbf{Corollary~\ref{cor:baseline.multiplicities}}\\Baseline\\multiplicities};
\node[result] (exceptional) at (2.2,0.3) {\textbf{Corollary~\ref{cor:exceptional.coincidences}}\\Exceptional\\coincidences};
\node[result] (HZ) at (6.2,0.3) {\textbf{Theorem~\ref{thm:Henze.Zirkler.spectrum}}\\Complete Henze--Zirkler\\spectrum};

\draw[dependency] (cartesian) -- (spherical);
\draw[dependency] (spherical) -- (block);
\draw[dependency] (block) -- (fredholm);
\draw[dependency] (spherical.east) -- (-4.20,7.2) -- (-4.20,6.45) -- (2.2,6.45) -- (complete.north);
\draw[dependency] (block.east) -- (complete.west);
\draw[dependency] (fredholm.east) -- (-0.50,2.4) -- ([xshift=4mm]complete.south west);
\draw[dependency] (positivity.south east) -- (complete.north west);
\draw[dependency] (compression.east) -- (complete.south west);

\draw[dependency] (complete.east) -- (reconstruction.west);
\draw[dependency] (spherical.west) -- (-8.20,7.2) -- (-8.20,10.75) -- (6.6,10.75) -- (6.6,6.55) -- (reconstruction.north);
\draw[dependency] (block.north) -- (-6.25,6.15) -- (4.30,6.15) -- (reconstruction.north west);

\draw[dependency] (complete) -- (baseline);
\draw[dependency] (complete) -- (exceptional);
\draw[dependency] (complete) -- (HZ);
\draw[dependency] (cartesian.west) -- (-8.55,9.6) -- (-8.55,0.3) -- (baseline.west);
\draw[dependency] (positivity.east) -- (8.85,8.0) -- (8.85,0.3) -- (HZ.east);
\end{tikzpicture}%
}
\caption{Proof-dependency graph.}
\label{fig:proof.dependencies}
\end{figure}
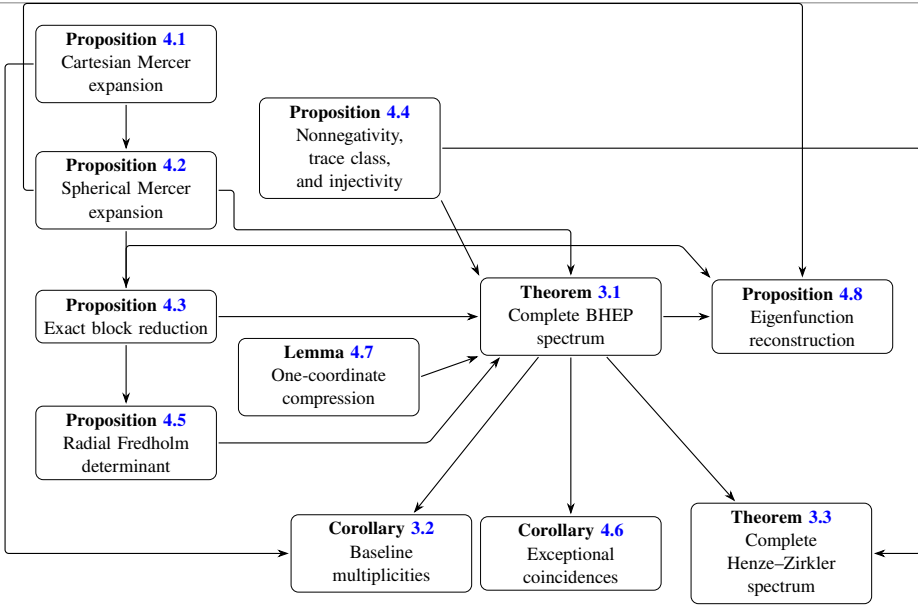

The order of the proofs in this section follows the dependencies among the theorem-like results displayed in Figure~\ref{fig:proof.dependencies}. An arrow $R_1 \longrightarrow R_2$ means that $R_1$ is used in the proof of $R_2$; only direct dependencies are shown.

\subsection{Proof of Proposition~\ref{prop:Cartesian.Mercer.expansion}}\label{subsec:proof.Cartesian.Mercer}

We begin with the one-dimensional case. Define
\[
e_n(x) \leqdef \frac{\Omega^{1/4}}{\sqrt{2^n n!}}H_n(\sqrt{\omega}x)\exp\left(-\frac{x^2}{\Omega + 1}\right), \qquad n \in \N_0.
\]
The identities in \eqref{eq:definitions.parameter.identities} give
\[
\frac{2}{\Omega + 1} + \frac{1}{2 \beta^2} = \omega, \qquad 2 \beta^2\omega = \Omega.
\]
Consequently, the change of variable $u = \sqrt{\omega}x$ and Hermite orthogonality yield
\[
\langle e_m,e_n\rangle_{\beta}
= \frac{\Omega^{1/2}}{\sqrt{2^{m+n}m!n!}}\frac{1}{\sqrt{2 \pi \beta^2}}\int_{\R}H_m(\sqrt{\omega}x)H_n(\sqrt{\omega}x)e^{-\omega x^2} \, \rd x
= \delta_{mn}.
\]
To prove completeness, suppose that $f$ is orthogonal to every $e_n$, and set
\[
h(u) \leqdef f(u/\sqrt{\omega})\exp\left\{\frac{u^2}{\omega(\Omega + 1)}\right\}.
\]
Since $\omega - 2/(\Omega + 1) = 1/(2\beta^2)$, the change of variable $u = \sqrt{\omega}x$ shows that $h \in L^2(\R,e^{-u^2} \, \rd u)$:
\[
\int_{\R} h(u)^2 e^{-u^2} \, \rd u
= \sqrt{\omega} \int_{\R} f(x)^2 \exp\left\{-\left(\omega - \frac{2}{\Omega + 1}\right)x^2\right\} \, \rd x
= \sqrt{2 \pi \beta^2 \omega} \, \|f\|_{\beta}^2 < \infty.
\]
The equations $\langle f,e_n\rangle_{\beta} = 0$ become
\[
\int_{\R}h(u)H_n(u)e^{-u^2} \, \rd u = 0, \qquad n \in \N_0.
\]
Completeness of the Hermite polynomials gives $h = 0$, and hence $f = 0$.

Mehler's formula \citep[Eq.~18.18.28]{KoornwinderEtAl2010}, written in the normalization of Section~\ref{sec:definitions}, is
\[
\sum_{n=0}^{\infty}\frac{q^n}{2^n n!}H_n(u)H_n(v) = \frac{1}{\sqrt{1 - q^2}}\exp\left\{\frac{2quv - q^2(u^2 + v^2)}{1 - q^2}\right\}, \qquad |q| < 1.
\]
With $u = \sqrt{\omega}s$ and $v = \sqrt{\omega}t$, the one-dimensional Mercer sum becomes
\begin{equation}\label{eq:proofs.Mehlers.formula}
\begin{aligned}
\sum_{n=0}^{\infty}(1 - q)q^ne_n(s)e_n(t)
&= (1 - q)\sqrt{\Omega} \exp\left(-\frac{s^2 + t^2}{\Omega + 1}\right) \sum_{n=0}^{\infty}\frac{q^n}{2^n n!}H_n(\sqrt{\omega}s)H_n(\sqrt{\omega}t) \\
&= \frac{(1 - q)\sqrt{\Omega}}{\sqrt{1 - q^2}}\exp\left\{-\frac{s^2 + t^2}{\Omega + 1} + \frac{2q\omega st - q^2\omega(s^2 + t^2)}{1 - q^2}\right\}.
\end{aligned}
\end{equation}
The constant prefactor reduces to one because
\[
\frac{(1 - q)\sqrt{\Omega}}{\sqrt{1 - q^2}} = \frac{(1 - q)\sqrt{(1 + q)/(1 - q)}}{\sqrt{(1 - q)(1 + q)}} = 1.
\]
Also, the last exponent in \eqref{eq:proofs.Mehlers.formula} simplifies to $-(s - t)^2/2$. Indeed, substituting $\Omega + 1 = 2 / (1 - q)$ and $\omega = (1 - q^2) / (2q)$ from \eqref{eq:definitions.parameter.identities} yields
\[
-\frac{s^2 + t^2}{\Omega + 1} + \frac{2q\omega st - q^2\omega(s^2 + t^2)}{1 - q^2}
= -\frac{1 - q}{2}(s^2 + t^2) + st - \frac{q}{2}(s^2 + t^2)
= -\frac{(s - t)^2}{2}.
\]
Hence
\begin{equation}\label{eq:1d.expansion}
\exp\left\{-\frac{(s - t)^2}{2}\right\} = \sum_{n=0}^{\infty}(1 - q)q^ne_n(s)e_n(t).
\end{equation}
On the diagonal, the nonnegative series in \eqref{eq:1d.expansion} sums to one. Hence the Cauchy--Schwarz inequality yields
\[
\sum_{m=0}^{\infty}(1 - q)q^m|e_m(s)e_m(t)| \leq \left\{\sum_{m=0}^{\infty}(1 - q)q^me_m(s)^2\right\}^{1/2}\left\{\sum_{m=0}^{\infty}(1 - q)q^me_m(t)^2\right\}^{1/2} = 1.
\]
The same inequality applied to the tails, together with Dini's theorem on compact sets for the nonnegative diagonal series, proves local uniform convergence. Since $\int_{\R}|e_n(t)| \, \rd \mu_{\beta}(t) \leq \|e_n\|_{\beta} = 1$, multiplication by $e_n(t)$ and dominated convergence justify termwise integration against $e_n(t) \, \rd \mu_{\beta}(t)$ and prove the one-dimensional case of the eigenvalue equation \eqref{eq:reduction.Cartesian.Mercer.eigenvalues}.

The $d$-dimensional kernel $\exp(-\|\bb{s} - \bb{t}\|^2/2) = \prod_{j=1}^d \exp\{-(s_j - t_j)^2/2\}$ and measure $\rd \mu_{\beta}(\bb{t}) = \prod_{j=1}^d \varphi_{\beta}(t_j) \, \rd t_j$ are tensor products of their one-dimensional counterparts. Taking products of the one-dimensional expansion \eqref{eq:1d.expansion} proves \eqref{eq:reduction.Cartesian.Mercer.kernel}. The same Cauchy--Schwarz domination, applied to the multi-index series, justifies integration against $e_{\bb{\nu}}(\bb{t}) \, \rd \mu_{\beta}(\bb{t})$ and proves \eqref{eq:reduction.Cartesian.Mercer.eigenvalues}. Finally,
\[
\sum_{\bb{\nu}\in\N_0^d}c_{\beta}q^{|\bb{\nu}|} = (1 - q)^d\left(\sum_{n=0}^{\infty}q^n\right)^d = 1,
\]
which proves both that $B_{\beta,d}$ is trace class and that $\tr(B_{\beta,d}) = 1$. This concludes the proof of Proposition~\ref{prop:Cartesian.Mercer.expansion}. \qed

\subsection{Proof of Proposition~\ref{prop:spherical.Mercer.expansion}}\label{subsec:proof.spherical.Mercer}

Let $d \geq 2$. By squaring \eqref{eq:reduction.spherical.Laguerre.functions}, multiplying by the density of $\mu_{\beta}$, and noting from \eqref{eq:definitions.parameter.identities} that $1 - q + 1/(2 \beta^2) = 2/(\Omega + 1) + 1/(2 \beta^2) = \omega$, we obtain the radial exponential
\[
\exp\left\{-\left(1 - q + \frac{1}{2 \beta^2}\right)r^2\right\} = e^{-\omega r^2}.
\]
The change of variable $u = \omega r^2$ therefore transforms the radial inner product into
\[
\frac{1}{2\omega^{\ell+d/2}}\int_0^{\infty}L_k^{(\ell+d/2-1)}(u)L_{k'}^{(\ell+d/2-1)}(u)u^{\ell+d/2-1}e^{-u} \, \rd u, \qquad k, k' \in \N_0.
\]
The normalization in \eqref{eq:reduction.spherical.Laguerre.functions}, the orthogonality of the Laguerre polynomials, and the orthogonality of the spherical harmonics $Y_{\ell,j}$ show that all the functions $\psi_{\ell,k,j}$ are orthonormal.

It remains to identify their eigenvalues. In the variable $\bb{z} = \sqrt{\omega}\bb{x}$, write $\bb{z} = s\bb{\theta}$ with $s = \|\bb{z}\| = \sqrt{\omega}\,r$, and introduce the polynomial differential operator
\[
\mathcal{L}P \leqdef -\Delta P + 2\bb{z}^{\top}\nabla P.
\]
The Hermite differential equation \citep[Table~18.8.1\#12]{KoornwinderEtAl2010} shows that
\begin{equation}\label{eq:product.Hermite}
\mathcal{L}\prod_{i=1}^dH_{\nu_i}(z_i) = 2|\bb{\nu}|\prod_{i=1}^dH_{\nu_i}(z_i).
\end{equation}
On the other hand, if $Y$ is a spherical harmonic of degree $\ell$ and $F$ is a twice differentiable function, then direct differentiation in polar coordinates yields
\[
\begin{aligned}
\mathcal{L}\{s^{\ell}Y(\bb{\theta})F(s^2)\}
&= - \left\{ F\Delta(s^{\ell}Y) + 2\nabla(s^{\ell}Y)^\top \nabla F + s^{\ell}Y\Delta F \right\} + 2F\bb{z}^\top \nabla(s^{\ell}Y) + 2s^{\ell}Y\bb{z}^\top \nabla F \\
&= - \left\{ 0 + 2\nabla(s^{\ell}Y)^\top (2\bb{z}F') + s^{\ell}Y(2d F' + 4s^2 F'') \right\} + 2\ell s^{\ell}Y F + 2s^{\ell}Y(2s^2 F') \\
&= - 4F'(\bb{z}^\top \nabla(s^{\ell}Y)) - s^{\ell}Y(2d F' + 4s^2 F'') + 2\ell s^{\ell}Y F + 4s^2 s^{\ell}Y F' \\
&= - 4\ell s^{\ell}Y F' - 2d s^{\ell}Y F' - 4s^2 s^{\ell}Y F'' + 2\ell s^{\ell}Y F + 4s^2 s^{\ell}Y F' \\
&= s^{\ell}Y(\bb{\theta})\left[ 2\ell F + 4 \{(s^2 - \ell - d/2)F' - s^2F''\} \right].
\end{aligned}
\]
The Laguerre differential equation \citep[Table~18.8.1\#8]{KoornwinderEtAl2010}
\[
u(L_k^{(\alpha)})''(u) + (\alpha + 1 - u)(L_k^{(\alpha)})'(u) + kL_k^{(\alpha)}(u) = 0
\]
with $u = s^2$, $\alpha = \ell + d/2 - 1$, and $F = L_k^{(\alpha)}$ then implies that
\[
\mathcal{L}\{s^{\ell}Y(\bb{\theta})L_k^{(\ell+d/2-1)}(s^2)\} = 2(\ell + 2k)\,s^{\ell}Y(\bb{\theta})L_k^{(\ell+d/2-1)}(s^2);
\]
the polynomial $s^{\ell}Y(\bb{\theta})L_k^{(\ell+d/2-1)}(s^2)$ is an eigenfunction of $\mathcal{L}$ with eigenvalue $2(\ell + 2k)$. Since
$$
H_n(z) = 2^nz^n + \ \textrm{terms of lower degree},
$$
the product Hermite polynomials form a basis of the polynomial algebra. Expanding the displayed polynomial in that basis and using \eqref{eq:product.Hermite} shows that only product Hermite polynomials of total index $\ell + 2k$ can occur. Since $s^2 = \omega r^2$, the displayed polynomial is, up to a nonzero constant, exactly the polynomial factor of $\psi_{\ell,k,j}(r\bb{\theta})$ in \eqref{eq:reduction.spherical.Laguerre.functions}. Moreover, $(1 - q)/2 = 1/(\Omega + 1)$, so the exponential factor of $\psi_{\ell,k,j}$ is the common exponential factor of the Cartesian Hermite functions. Hence $\psi_{\ell,k,j}$ belongs to the span of the product Hermite functions of total index $\ell + 2k$. To see that no dimension is missing, observe that, by a symmetry and telescoping argument,
\[
\sum_{\substack{\ell,k \in \N_0\\\ell+2k=N}}h_{d,\ell}
= \sum_{k=0}^{\lfloor N/2\rfloor}h_{d,N-2k}
= \sum_{k=0}^{\lfloor N/2\rfloor}\left\{\binom{d + N - 2k - 1}{N - 2k} - \binom{d + N - 2k - 3}{N - 2k - 2}\right\}
= \binom{d + N - 1}{N}.
\]
The last number is exactly the number of multi-indices $\bb{\nu}$ with $|\bb{\nu}| = N$. Thus, at every total degree $N$, the spherical functions $\{\psi_{\ell,k,j} : \ell + 2k = N,\ \ell\in \N_0,\ k\in \N_0,\ 1\leq j\leq h_{d,\ell}\}$ constitute an orthonormal basis of the Cartesian Hermite eigenspace. Proposition~\ref{prop:Cartesian.Mercer.expansion} now gives \eqref{eq:reduction.spherical.Mercer.eigenvalues} and completeness. In dimension one, the Cartesian Hermite functions $e_n$ defined in \eqref{eq:cartesian.eigenfunctions} split exactly into their even and odd subsequences $e_{2k}$ and $e_{2k+1}$, which we have identified with $\psi_{0,k,1}$ and $\psi_{1,k,1}$, respectively, proving the last assertion. This concludes the proof of Proposition~\ref{prop:spherical.Mercer.expansion}. \qed

\subsection{Proof of Proposition~\ref{prop:exact.block.reduction}}\label{subsec:proof.block.reduction}

Recall that $g(\bb{x}) \leqdef \exp(-\|\bb{x}\|^2/2)$, $\bb{x} \in \R^d$, from \eqref{eq:def.g}. Expanding the squared norm gives
\[
e^{-\|\bb{s} - \bb{t}\|^2/2} = g(\bb{s})g(\bb{t})e^{\bb{s}^{\top}\bb{t}}.
\]
The second term in the definition of the kernel $K$ in \eqref{eq:definitions.BHEP.kernel} removes the terms of Taylor degrees $0$, $1$, and $2$ from $e^{\bb{s}^{\top}\bb{t}}$. Since $f \otimes f$ has kernel $f(\bb{s})f(\bb{t})$, the constant term gives $g \otimes g$. The identity $\bb{s}^{\top}\bb{t} = \sum_{i=1}^d s_i t_i$ shows similarly that the linear term gives $\sum_{i=1}^d(x_i g) \otimes (x_i g)$.

Assume that $d \geq 2$. The matrices $I_d/\sqrt{d}, E_1,\ldots,E_{h_{d,2}}$ form a Frobenius-orthonormal basis of the real symmetric matrices. Applying Parseval's identity to $\bb{s}\bb{s}^{\top}$ and $\bb{t}\bb{t}^{\top}$ gives
\[
\begin{aligned}
(\bb{s}^{\top}\bb{t})^2
= \tr(\bb{s}\bb{s}^{\top}\bb{t}\bb{t}^{\top})
&= \tr\left(\bb{s}\bb{s}^{\top}\frac{I_d}{\sqrt{d}}\right)\tr\left(\bb{t}\bb{t}^{\top}\frac{I_d}{\sqrt{d}}\right) + \sum_{j=1}^{h_{d,2}}\tr(\bb{s}\bb{s}^{\top}E_j)\tr(\bb{t}\bb{t}^{\top}E_j) \\
&= \frac{\|\bb{s}\|^2\|\bb{t}\|^2}{d} + \sum_{j=1}^{h_{d,2}}(\bb{s}^{\top}E_j\bb{s})(\bb{t}^{\top}E_j\bb{t}).
\end{aligned}
\]
After multiplication by $g(\bb{s})g(\bb{t})/2$, the right-hand side represents the kernel corresponding to the quadratic terms, namely
\begin{equation}\label{eq:traceless.sum.2nd.term}
\left(\frac{\|\bb{x}\|^2g}{\sqrt{2d}}\right)\otimes\left(\frac{\|\bb{x}\|^2g}{\sqrt{2d}}\right) + \sum_{j=1}^{h_{d,2}}\left(\frac{\bb{x}^{\top}E_j\bb{x}}{\sqrt{2}}g\right)\otimes\left(\frac{\bb{x}^{\top}E_j\bb{x}}{\sqrt{2}}g\right).
\end{equation}
Combining the constant, linear, and two quadratic contributions proves \eqref{eq:reduction.finite.rank.decomposition}. When $d = 1$, the traceless sum (i.e., the summation term in \eqref{eq:traceless.sum.2nd.term}) is absent.

This also locates the correction in the angular decomposition. The functions $g$ and $\|\bb{x}\|^2g$ are radial and belong to angular degree $0$, while $x_i g$ has angular part $\theta_i$ and belongs to degree $1$. Moreover, $\bb{x}^{\top}E_j\bb{x}$ is homogeneous of degree $2$ and harmonic, since
\[
\Delta(\bb{x}^{\top}E_j\bb{x}) = 2 \, \tr(E_j) = 0,
\]
so its restriction to the sphere is a spherical harmonic of angular degree $2$. Hence the correction
$$
g \otimes g + \sum_{i=1}^d(x_i g) \otimes (x_i g) + \left(\frac{\|\bb{x}\|^2g}{\sqrt{2d}}\right)\otimes\left(\frac{\|\bb{x}\|^2g}{\sqrt{2d}}\right) + \sum_{j=1}^{h_{d,2}}\left(\frac{\bb{x}^{\top}E_j\bb{x}}{\sqrt{2}}g\right)\otimes\left(\frac{\bb{x}^{\top}E_j\bb{x}}{\sqrt{2}}g\right)
$$
vanishes on every sector of degree $\ell \geq 3$. When $d = 1$, the radial functions are even and $x_1g$ is odd.

For completeness, the projection coefficients are calculated explicitly. The Laguerre generating function gives, for $\alpha > -1$ and $s > 1$,
\begin{equation}\label{eq:proofs.Laguerre.Laplace.transform}
\int_0^{\infty}u^{\alpha}e^{-su}L_k^{(\alpha)}(u) \, \rd u = \frac{\Gamma(k + \alpha + 1)}{k!}\frac{(s - 1)^k}{s^{k+\alpha+1}}.
\end{equation}
This is a standard Laplace transform of the generalized Laguerre polynomial \citep[18.17.34]{KoornwinderEtAl2010}. Alternatively, multiplying
$$
\sum_{k=0}^{\infty}L_k^{(\alpha)}(u)z^k = (1 - z)^{-\alpha-1}\exp\{-uz/(1 - z)\}
$$
by $u^{\alpha}e^{-su}$, integrating, and comparing coefficients directly proves \eqref{eq:proofs.Laguerre.Laplace.transform}.

We next calculate the projection coefficients of the functions defining the finite-rank correction (namely $g$, $x_i g$, $\|\bb{x}\|^2 g/\sqrt{2d}$, and $\bb{x}^{\top}E_j\bb{x} g/\sqrt{2}$) with respect to the spherical basis when $d \geq 2$. Write $\bb{x} = r\bb{\theta}$ in polar coordinates, noting that $g(\bb{x}) = \exp(-r^2/2)$. The three exponential factors in each $L^2(\mu_{\beta})$ inner product (originating from the factor of $g$ present in each removed function; the basis function $\psi_{\ell,k,j}$; and the density of the Gaussian measure $\mu_{\beta}$) multiply to
\[
\exp\left(-\frac{r^2}{2}\right)\exp\left(-\frac{1 - q}{2}r^2\right)\exp\left(-\frac{r^2}{2\beta^2}\right) = \exp\left(-\frac{r^2}{2q}\right),
\]
because $\beta^{-2} = (1 - q)^2/q$ by \eqref{eq:definitions.parameter.identities}. For the coefficients of $g$, $x_i g$, and $(\bb{x}^{\top}E_i\bb{x})g/\sqrt{2}$, the two polynomial factors contribute $r^{2\ell}$ in angular degree $\ell \in \{0,1,2\}$. Let $s = 1/(2q\omega) = 1/(1 - \rho)$, so that $(s - 1)/s = \rho$. The substitution $u = \omega r^2$ and \eqref{eq:proofs.Laguerre.Laplace.transform} with $\alpha = \ell + d/2 - 1$ give
\begin{align}\label{eq:proofs.matching.radial.integral}
R_{\ell,k}
&\leqdef \int_0^{\infty}r^{2\ell+d-1}\exp\left(-\frac{r^2}{2q}\right)L_k^{(\ell+d/2-1)}(\omega r^2) \, \rd r \nonumber \\
&= \frac{1}{2\omega^{\ell+d/2}}\int_0^{\infty}u^{\ell+d/2-1}e^{-su}L_k^{(\ell+d/2-1)}(u) \, \rd u
= \frac{\Gamma(k + \ell + \frac{d}{2})}{2\omega^{\ell+d/2}k!}(1 - \rho)^{\ell+d/2}\rho^k.
\end{align}

The specific choice of basis in \eqref{eq:reduction.aligned.harmonics} gives the three matching angular factors
\begin{equation}\label{eq:3.identities}
\begin{aligned}
\int_{\mathbb{S}^{d-1}}Y_{0,1}(\bb{\theta}) \, \rd \bb{\theta} &= \sqrt{\sigma_{d-1}}, \\
\int_{\mathbb{S}^{d-1}}\theta_iY_{1,j}(\bb{\theta}) \, \rd \bb{\theta} &= \left(\frac{\sigma_{d-1}}{d}\right)^{1/2}\delta_{ij}, \\
\frac{1}{\sqrt{2}}\int_{\mathbb{S}^{d-1}}(\bb{\theta}^{\top}E_i\bb{\theta})Y_{2,j}(\bb{\theta}) \, \rd \bb{\theta} &= \left\{\frac{\sigma_{d-1}}{d(d + 2)}\right\}^{1/2}\delta_{ij}.
\end{aligned}
\end{equation}
Converting expectation with respect to the uniform probability measure on $\mathbb{S}^{d-1}$ into integration with respect to surface measure shows that the case $m = 1$ of \citet[Theorem~2, p.~2405]{VignatBhatnagar2008} gives the second identity in \eqref{eq:3.identities}. For the third identity in \eqref{eq:3.identities}, the case $m = 2$ of the same result gives the fourth-order spherical moment formula
\[
\int_{\mathbb{S}^{d-1}}\theta_a\theta_b\theta_c\theta_e \, \rd \bb{\theta}
= \frac{\sigma_{d-1}}{d(d + 2)}\left(\delta_{ab}\delta_{ce} + \delta_{ac}\delta_{be} + \delta_{ae}\delta_{bc}\right).
\]
Expanding both quadratic forms and applying this formula gives
\[
\int_{\mathbb{S}^{d-1}}(\bb{\theta}^{\top}E_i\bb{\theta})(\bb{\theta}^{\top}E_j\bb{\theta}) \, \rd \bb{\theta}
= \frac{\sigma_{d-1}}{d(d + 2)}\left\{\tr(E_i)\tr(E_j) + 2\tr(E_iE_j)\right\}
= \frac{2\sigma_{d-1}}{d(d + 2)}\delta_{ij},
\]
where the last equality uses $\tr(E_i) = \tr(E_j) = 0$ and $\tr(E_iE_j) = \delta_{ij}$. Multiplication by the factor $1/\sqrt{2}$ from the removed quadratic function and by the normalization factor $\{d(d + 2)/(2\sigma_{d-1})\}^{1/2}$ in $Y_{2,j}$ gives the third angular identity in \eqref{eq:3.identities}.

Define
\[
\alpha_0 \leqdef \int_{\mathbb{S}^{d-1}}Y_{0,1}(\bb{\theta}) \, \rd \bb{\theta}, \qquad
\alpha_1 \leqdef \int_{\mathbb{S}^{d-1}}\theta_iY_{1,i}(\bb{\theta}) \, \rd \bb{\theta} \quad (1 \leq i \leq d),
\]
and
\[
\alpha_2 \leqdef \frac{1}{\sqrt{2}}\int_{\mathbb{S}^{d-1}}(\bb{\theta}^{\top}E_i\bb{\theta})Y_{2,i}(\bb{\theta}) \, \rd \bb{\theta} \quad (1 \leq i \leq h_{d,2}).
\]
The values of $\alpha_1$ and $\alpha_2$ do not depend on the chosen index $i$. These three quantities are the positive coefficients of the matching angular terms in \eqref{eq:3.identities}. Using $\sigma_{d-1} = 2\pi^{d/2}/\Gamma(d/2)$ gives
\[
\alpha_{\ell}^2 = \frac{2^{1-\ell}\pi^{d/2}}{\Gamma(\frac{d}{2} + \ell)}, \qquad \ell \in \{0,1,2\}.
\]

For $k \in \N_0$, define the matching projection coefficients
\[
b_{0,k} \leqdef \langle g,\psi_{0,k,1}\rangle_{\beta}, \qquad
b_{1,k} \leqdef \langle x_i g,\psi_{1,k,i}\rangle_{\beta} \quad (1 \leq i \leq d),
\]
and
\[
b_{2,k} \leqdef \left\langle\frac{\bb{x}^{\top}E_i\bb{x}}{\sqrt{2}}g,\psi_{2,k,i}\right\rangle_{\beta} \quad (1 \leq i \leq h_{d,2}).
\]
The last two values do not depend on the chosen index $i$. The polar-coordinate formula, the density factor $(2\pi\beta^2)^{-d/2}$ from $\mu_{\beta}$, the normalization in \eqref{eq:reduction.spherical.Laguerre.functions}, the angular factor $\alpha_{\ell}$, and the radial integral $R_{\ell,k}$ in \eqref{eq:proofs.matching.radial.integral} give
\begin{equation}\label{eq:proofs.matching.projection.coefficient}
b_{\ell,k}
= (2\pi\beta^2)^{-d/4}\left\{\frac{2\omega^{\ell+d/2}k!}{\Gamma(k + \ell + \frac{d}{2})}\right\}^{1/2}\alpha_{\ell}R_{\ell,k}, \qquad \ell \in \{0,1,2\}.
\end{equation}
Squaring \eqref{eq:proofs.matching.projection.coefficient} and substituting the formulas for $\alpha_{\ell}^2$ and $R_{\ell,k}$ yields
\begin{equation}\label{eq:calc}
\frac{(\frac{d}{2} + \ell)_k}{k!}\frac{(1 - \rho)^{d+2\ell}\rho^{2k}}{2^{\ell}(2\beta^2)^{d/2}\omega^{d/2+\ell}}
= (1 - q)^dq^{\ell}(1 - \rho)^{d/2+\ell}\frac{(\frac{d}{2} + \ell)_k}{k!}\rho^{2k}
= \Lambda_{\ell,k}\pi_{\ell,k}.
\end{equation}
The first equality in \eqref{eq:calc} uses $\frac{1 - \rho}{2\beta^2\omega} = (1 - q)^2$ and $\frac{1 - \rho}{2\omega} = q$, which follow from \eqref{eq:definitions.parameter.identities}. The last equality in \eqref{eq:calc} follows directly from \eqref{eq:definitions.Gaussian.levels.and.weights}. Every factor on the right-hand side of \eqref{eq:proofs.matching.projection.coefficient} is strictly positive, because $\alpha_{\ell} > 0$ and the final expression for $R_{\ell,k}$ in \eqref{eq:proofs.matching.radial.integral} is positive. Hence each matching coefficient $b_{\ell,k}$ is the positive square root of $\Lambda_{\ell,k}\pi_{\ell,k}$. Restoring the Kronecker deltas for distinct angular indices gives
\begin{align*}
\langle g,\psi_{0,k,1}\rangle_{\beta} &= \{\Lambda_{0,k}\pi_{0,k}\}^{1/2}, \\
\langle x_i g,\psi_{1,k,j}\rangle_{\beta} &= \delta_{ij}\{\Lambda_{1,k}\pi_{1,k}\}^{1/2}, \\
\left\langle\frac{\bb{x}^{\top}E_i\bb{x}}{\sqrt{2}}g,\psi_{2,k,j}\right\rangle_{\beta} &= \delta_{ij}\{\Lambda_{2,k}\pi_{2,k}\}^{1/2}.
\end{align*}
This proves \eqref{eq:reduction.coefficient.radial.first}, \eqref{eq:reduction.coefficient.linear}, and \eqref{eq:reduction.coefficient.quadratic} when $d \geq 2$.

When $d = 1$, the quadratic transformations in \citet[Eqs.~18.7.19--18.7.20]{KoornwinderEtAl2010} give, for $z \in \C$,
\[
H_{2k}(z) = (-1)^k 2^{2k}k!L_k^{(-1/2)}(z^2), \qquad H_{2k+1}(z) = (-1)^k 2^{2k+1}k!zL_k^{(1/2)}(z^2).
\]
Substitution into \eqref{eq:cartesian.eigenfunctions}, together with the choices $\psi_{0,k,1} = (-1)^k e_{2k}$ and $\psi_{1,k,1} = (-1)^k e_{2k+1}$, gives
\[
\begin{aligned}
\psi_{0,k,1}(x)
&= \frac{2^k k! \Omega^{1/4}}{\sqrt{(2k)!}}L_k^{(-1/2)}(\omega x^2)\exp\left(-\frac{x^2}{\Omega + 1}\right), \\
\psi_{1,k,1}(x)
&= \frac{2^{k+1/2} k! \Omega^{1/4}\sqrt{\omega}}{\sqrt{(2k + 1)!}}xL_k^{(1/2)}(\omega x^2)\exp\left(-\frac{x^2}{\Omega + 1}\right).
\end{aligned}
\]
The first function is even and the second is odd. Splitting the corresponding inner products over the two half-lines, substituting $u = \omega x^2$, and applying \eqref{eq:proofs.Laguerre.Laplace.transform} with $\alpha = -1/2$ and $\alpha = 1/2$, respectively, gives
\[
\langle g,\psi_{0,k,1}\rangle_{\beta} = \{\Lambda_{0,k}\pi_{0,k}\}^{1/2}, \qquad
\langle xg,\psi_{1,k,1}\rangle_{\beta} = \{\Lambda_{1,k}\pi_{1,k}\}^{1/2}.
\]
Thus \eqref{eq:reduction.coefficient.radial.first} and \eqref{eq:reduction.coefficient.linear} also hold when $d = 1$, and \eqref{eq:reduction.coefficient.quadratic} is trivial because $h_{1,2} = 0$. The remaining coefficient \eqref{eq:reduction.coefficient.radial.second} is established next for every $d \in \N$.

For the second radial coefficient, \eqref{eq:reduction.coefficient.radial.second}, let $I_k(s)$ denote the left-hand side of \eqref{eq:proofs.Laguerre.Laplace.transform} with $\alpha = d/2 - 1$. Differentiation under the integral sign shows that the integral with one additional factor $u$ is $-I_k'(s)$. Logarithmic differentiation of the right-hand side of \eqref{eq:proofs.Laguerre.Laplace.transform} gives
\[
-\frac{I_k'(s)}{I_k(s)} = \frac{k + d/2}{s} - \frac{k}{s - 1}.
\]
Under $u = \omega r^2$, the additional factor $\|\bb{x}\|^2/\sqrt{2d}$ contributes $u/(\sqrt{2d}\omega)$. At $s = 1/(1 - \rho)$, the ratio of the second radial coefficient, \eqref{eq:reduction.coefficient.radial.second}, to the first, \eqref{eq:reduction.coefficient.radial.first}, is therefore
\[
\dfrac{\left\langle\dfrac{\|\bb{x}\|^2g}{\sqrt{2d}},\psi_{0,k,1}\right\rangle_{\beta}}{\langle g,\psi_{0,k,1}\rangle_{\beta}}
= \frac{1 - \rho}{\sqrt{2d}\omega}\left\{\frac{d}{2} - \frac{1 - \rho}{\rho}k\right\}
= \frac{dq - 2(1 - \rho)k/q}{\sqrt{2d}} = y_k,
\]
where $\rho = q^2$ and $\omega = (1 - \rho)/(2q)$ from \eqref{eq:definitions.parameter.identities} were used in the last equality. Combining this ratio with \eqref{eq:reduction.coefficient.radial.first} proves \eqref{eq:reduction.coefficient.radial.second}.

It remains to pass to the sequence-space representations in order to prove \eqref{eq:reduction.rank.one.matrix} and \eqref{eq:reduction.radial.matrix}. Fix an angular copy of degree $\ell \in \{1,2\}$, where $\ell = 2$ is present only for $d \geq 2$, and identify $\sum_{k=0}^{\infty}a_k\psi_{\ell,k,j}$ with $(a_k)_{k \in \N_0} \in \ell^2(\N_0)$. In these coordinates, $B_{\beta,d}$ is represented by $D_{\ell}$, defined in \eqref{eq:reduction.rank.one.matrix}, and the coefficient sequence of the unique removed function is
\[
b_{\ell} \leqdef (b_{\ell,k})_{k \in \N_0} = \left(\{\Lambda_{\ell,k}\pi_{\ell,k}\}^{1/2}\right)_{k \in \N_0} = D_{\ell}^{1/2}u_{\ell},
\]
where $u_{\ell} \leqdef (\sqrt{\pi_{\ell,k}})_{k \in \N_0}$. Since $\sum_{k=0}^{\infty}\pi_{\ell,k} = 1$, the sequence $u_{\ell}$ has unit norm. Restricting $A_{\beta,d}$ to the whole invariant subspace $\mathcal{V}_{\ell,j} = \overline{\operatorname{span}}\{\psi_{\ell,k,j} : k \in \N_0\}$, identified with $\ell^2(\N_0)$ as above, gives
\[
D_{\ell} - b_{\ell} \otimes b_{\ell} = D_{\ell}^{1/2}(I - u_{\ell} \otimes u_{\ell})D_{\ell}^{1/2} = D_{\ell}^{1/2}P_{\ell}D_{\ell}^{1/2},
\]
which proves \eqref{eq:reduction.rank.one.matrix}.

In the radial sector, first consider the probability generating function of $(\pi_{0,k})_{k \in \N_0}$, defined in \eqref{eq:definitions.Gaussian.levels.and.weights} with $\ell = 0$:
\[
G_0(z) \leqdef \sum_{k=0}^{\infty}\pi_{0,k}z^k = (1-\rho)^{d/2}\sum_{k=0}^{\infty}\frac{(d/2)_k}{k!}(\rho z)^k = \left(\frac{1 - \rho}{1 - \rho z}\right)^{d/2};
\]
see, e.g., \citet[Eq.~15.4.6]{DLMF15}. Its first two derivatives at $z = 1$ give
\[
G_0'(1) = \frac{(d/2) \rho}{1 - \rho}, \qquad G_0''(1) = \frac{(d/2)(d/2 + 1)\rho^2}{(1 - \rho)^2}.
\]
The mean and variance of $k$ under these weights are, respectively,
\[
\sum_{k=0}^{\infty}k\,\pi_{0,k} = G_0'(1), \qquad \sum_{k=0}^{\infty}\left(k - G_0'(1)\right)^2\pi_{0,k} = G_0''(1) + G_0'(1) - \{G_0'(1)\}^2 = \frac{(d/2) \rho}{(1 - \rho)^2}.
\]
Since $\rho = q^2$, the quantity $y_k$ in \eqref{eq:results.radial.y} can be written as the negative of the standardized value
\[
y_k = -\frac{k - \frac{(d/2) \rho}{1 - \rho}}{\left\{\frac{(d/2) \rho}{(1 - \rho)^2}\right\}^{1/2}}.
\]
The total mass, mean, and variance identities are therefore equivalent to
\begin{equation}\label{eq:proofs.radial.isometry.moments}
\sum_{k=0}^{\infty}\pi_{0,k} = 1, \qquad \sum_{k=0}^{\infty}\pi_{0,k}y_k = 0, \qquad \sum_{k=0}^{\infty}\pi_{0,k}y_k^2 = 1.
\end{equation}
Since $\bb{v}_k = (1,y_k)^{\top}$, these identities give, by the definition of $U_0$ in \eqref{eq:def.U.0},
\[
U_0^*U_0 = \sum_{k=0}^{\infty}\pi_{0,k}\bb{v}_k\bb{v}_k^{\top} = I_2.
\]
Thus $U_0$ is an isometry. Under the radial identification with $\ell^2(\N_0)$, the coefficient sequences of $g$ and $\|\bb{x}\|^2g/\sqrt{2d}$ are the two columns of $D_0^{1/2}U_0$. Their two rank-one corrections therefore sum to $D_0^{1/2}U_0U_0^*D_0^{1/2}$, and the radial restriction of $A_{\beta,d}$ is
\[
D_0 - D_0^{1/2}U_0U_0^*D_0^{1/2} = D_0^{1/2}(I - U_0U_0^*)D_0^{1/2} = D_0^{1/2}Q_0D_0^{1/2}.
\]
This proves \eqref{eq:reduction.radial.matrix}.

Finally, all removed functions have angular degree $0$, $1$, or $2$, so they are orthogonal to every sector of degree $\ell \geq 3$. On each such sector, $A_{\beta,d}$ therefore agrees with $B_{\beta,d}$. This completes the proof of Proposition~\ref{prop:exact.block.reduction}. \qed

\subsection{Proof of Proposition~\ref{prop:positivity.trace.injectivity}}\label{subsec:proof.positivity}

The exponential series gives the pointwise feature expansion of $K$, defined in \eqref{eq:definitions.BHEP.kernel}:
\begin{equation}\label{eq:proofs.positive.feature.expansion}
K(\bb{s},\bb{t}) = \sum_{\substack{\bb{\nu}\in\N_0^d\\|\bb{\nu}|\geq 3}}\phi_{\bb{\nu}}(\bb{s})\phi_{\bb{\nu}}(\bb{t}), \qquad \phi_{\bb{\nu}}(\bb{x}) \leqdef \frac{\bb{x}^{\bb{\nu}}}{\sqrt{\bb{\nu}!}}e^{-\|\bb{x}\|^2/2}.
\end{equation}
On the diagonal,
\[
K(\bb{t},\bb{t}) = 1 - \left(1 + \|\bb{t}\|^2 + \frac{\|\bb{t}\|^4}{2}\right)e^{-\|\bb{t}\|^2},
\]
which is nonnegative and integrable with respect to $\mu_{\beta}$. Tonelli's theorem therefore gives
\begin{equation}\label{eq:trace.norm.finite}
\sum_{|\bb{\nu}|\geq 3}\|\phi_{\bb{\nu}}\|_{\beta}^2 = \int_{\R^d}K(\bb{t},\bb{t}) \, \rd \mu_{\beta}(\bb{t}) < \infty.
\end{equation}
Let $K_N$ be the truncation of \eqref{eq:proofs.positive.feature.expansion} to $3 \leq |\bb{\nu}| \leq N$, and set $\varepsilon_N(\bb{t}) \leqdef K(\bb{t},\bb{t}) - K_N(\bb{t},\bb{t})$. Cauchy--Schwarz in the feature index gives
\[
|K(\bb{s},\bb{t}) - K_N(\bb{s},\bb{t})|^2 \leq \varepsilon_N(\bb{s})\varepsilon_N(\bb{t}).
\]
Since $\int_{\R^d} \varepsilon_N(\bb{t}) \, \rd \mu_{\beta}(\bb{t}) \to 0$, it follows that $K_N \to K$ in $L^2(\mu_{\beta}\otimes\mu_{\beta})$. The series $\sum_{|\bb{\nu}|\geq 3}\phi_{\bb{\nu}}\otimes\phi_{\bb{\nu}}$ also converges absolutely in trace norm (i.e., $\sum_{|\bb{\nu}|\geq 3}\|\phi_{\bb{\nu}}\otimes\phi_{\bb{\nu}}\|_1 < \infty$) because $\|\phi_{\bb{\nu}}\otimes\phi_{\bb{\nu}}\|_1 = \tr(\phi_{\bb{\nu}}\otimes\phi_{\bb{\nu}}) = \langle\phi_{\bb{\nu}},\phi_{\bb{\nu}}\rangle_{\beta} = \|\phi_{\bb{\nu}}\|_{\beta}^2$, and the sum of these trace norms is finite by \eqref{eq:trace.norm.finite}. Trace-norm convergence implies Hilbert--Schmidt convergence, and the Hilbert--Schmidt norm of an integral operator equals the $L^2$ norm of its kernel. Uniqueness of the Hilbert--Schmidt limit therefore shows that the trace-norm limit has kernel $K$ and is exactly $A_{\beta,d}$. Hence $A_{\beta,d}$ is nonnegative, self-adjoint, and trace class.

If $\bb{X}$ has distribution $\mu_{\beta}$, the gamma integral \citep[Eq.~5.9.1]{DLMF5}, applied to the chi-square density, gives $\EE(e^{-u\|\bb{X}\|^2}) = (1 + 2\beta^2u)^{-d/2}$ for $u \geq 0$. Differentiating twice with respect to $u$ and evaluating at $u = 1$ gives
\begin{align*}
\EE(e^{-\|\bb{X}\|^2}) &= (1 + 2\beta^2)^{-d/2},\\
\EE(\|\bb{X}\|^2e^{-\|\bb{X}\|^2}) &= d\beta^2(1 + 2\beta^2)^{-d/2-1},\\
\EE(\|\bb{X}\|^4e^{-\|\bb{X}\|^2}) &= d(d + 2)\beta^4(1 + 2\beta^2)^{-d/2-2}.
\end{align*}
Integrating the diagonal of $K$, i.e., $K(\bb{t},\bb{t}) = 1 - \left(1 + \|\bb{t}\|^2 + \|\bb{t}\|^4/2\right)e^{-\|\bb{t}\|^2}$, proves \eqref{eq:reduction.trace.formula}.

It remains to prove injectivity. If $A_{\beta,d}f = 0$, then in particular $\langle A_{\beta,d}f,f\rangle_{\beta} = 0$. Substituting the feature expansion \eqref{eq:proofs.positive.feature.expansion} for $K$ and using its absolute trace-norm convergence to justify interchanging sum and integral gives
\[
0 = \int_{\R^d}\int_{\R^d}K(\bb{s},\bb{t})f(\bb{t})f(\bb{s}) \, \rd \mu_{\beta}(\bb{t}) \, \rd \mu_{\beta}(\bb{s}) = \sum_{|\bb{\nu}|\geq 3}c_{\bb{\nu}}^2, \qquad c_{\bb{\nu}} \leqdef \int_{\R^d}f(\bb{t})\phi_{\bb{\nu}}(\bb{t}) \, \rd \mu_{\beta}(\bb{t}).
\]
This sum of nonnegative terms equals zero, so every term vanishes and $c_{\bb{\nu}} = 0$ for every $\bb{\nu}$ with $|\bb{\nu}| \geq 3$; this is equivalent to
\begin{equation}\label{eq:proofs.high.moments.zero}
\int_{\R^d}f(\bb{t})\bb{t}^{\bb{\nu}}e^{-\|\bb{t}\|^2/2} \, \rd \mu_{\beta}(\bb{t}) = 0, \qquad |\bb{\nu}| \geq 3.
\end{equation}
For $\bb{z} \in \C^d$, define
\[
F(\bb{z}) \leqdef \int_{\R^d}f(\bb{t})e^{-\|\bb{t}\|^2/2}e^{\bb{z}^{\top}\bb{t}} \, \rd \mu_{\beta}(\bb{t}).
\]
Applying the Cauchy--Schwarz inequality and completing the square show that the integrand is dominated by an integrable function uniformly for $\bb{z}$ in each compact subset of $\C^d$, so $F$ is entire and differentiation may be carried out under the integral sign. Equation \eqref{eq:proofs.high.moments.zero} implies that every Taylor coefficient of total degree at least three vanishes, so $F$ is a polynomial of degree at most two. For $\bb{z} = \ii\bb{\xi}$, $F$ is the Fourier transform of the integrable Lebesgue density $f(\bb{t})e^{-\|\bb{t}\|^2/2}\varphi_{\beta}(\bb{t})$. By the Riemann--Lebesgue lemma \citep[Theorem~8.22(f), p.~249]{Folland1999}, $F(\ii\bb{\xi}) \to 0$ as $\|\bb{\xi}\| \to \infty$. The only polynomial with this property is the zero polynomial. The uniqueness of the Fourier transform \citep[Corollary~8.27, p.~252]{Folland1999} then yields $f(\bb{t})e^{-\|\bb{t}\|^2/2}\varphi_{\beta}(\bb{t}) = 0$ almost everywhere, and therefore $f = 0$ in $L^2(\mu_{\beta})$. Hence $A_{\beta,d}$ is injective, completing the proof of Proposition~\ref{prop:positivity.trace.injectivity}. \qed

\subsection{Proof of Proposition~\ref{prop:radial.Fredholm.determinant}}\label{subsec:proof.radial.Fredholm}

In the radial sequence model of \eqref{eq:reduction.radial.matrix}, let $D \leqdef \operatorname{diag}(d_0,d_1,\ldots)$ and $W \leqdef D^{1/2}U_0$. Then
\[
A_{\beta,d}^{(0)} = D - WW^*.
\]
For $z \notin \{d_k^{-1} : k \in \N_0\}$, factorization and multiplicativity of the Fredholm determinant give
\[
\det\nolimits_F(I - zA_{\beta,d}^{(0)}) = \det\nolimits_F(I - zD)\det\{I_2 + zW^*(I - zD)^{-1}W\}.
\]
The matrix $(I-zD)^{-1}$ is diagonal with $k$th entry $(1-zd_k)^{-1}$, so $\det\nolimits_F(I-zD) = \prod_{k=0}^{\infty}(1-zd_k)$. Recall from \eqref{eq:def.U.0} that $W = D^{1/2}U_0$ has $k$th row $\sqrt{d_k}\cdot\sqrt{\pi_{0,k}}\bb{v}_k^{\top} = \{d_k\pi_{0,k}\}^{1/2}\bb{v}_k^{\top}$. Hence
\[
W^*(I-zD)^{-1}W = \sum_{k=0}^{\infty}\frac{1}{1-zd_k}\left(\{d_k\pi_{0,k}\}^{1/2}\bb{v}_k\right)\left(\{d_k\pi_{0,k}\}^{1/2}\bb{v}_k\right)^{\top} = \sum_{k=0}^{\infty}\frac{d_k\pi_{0,k}}{1-zd_k}\bb{v}_k\bb{v}_k^{\top},
\]
and substituting these two expressions gives
\[
\det\nolimits_F(I - zA_{\beta,d}^{(0)})
= \prod_{k=0}^{\infty}(1 - zd_k)\det\left\{I_2 + z\sum_{k=0}^{\infty}\frac{d_k\pi_{0,k}}{1 - zd_k}\bb{v}_k\bb{v}_k^{\top}\right\}.
\]
By \eqref{eq:results.radial.Fredholm.determinant}, the right-hand side is exactly $\mathcal{D}_0(z)$, so $\det\nolimits_F(I - zA_{\beta,d}^{(0)}) = \mathcal{D}_0(z)$ for every $z \notin \{d_k^{-1} : k \in \N_0\}$.

For vectors $\bb{r}_0,\ldots,\bb{r}_N \in \R^2$ and scalars $a_0,\ldots,a_N$, expansion of a $2 \times 2$ determinant gives
\[
\det\left(I_2 + \sum_{k=0}^N a_k\bb{r}_k\bb{r}_k^{\top}\right) = 1 + \sum_{k=0}^N a_k\|\bb{r}_k\|^2 + \sum_{0 \leq i < j \leq N}a_ia_j\det(\bb{r}_i,\bb{r}_j)^2.
\]
For fixed $z \notin \{d_k^{-1} : k \in \N_0\}$, the denominators $1 - zd_k$ are bounded away from zero. Apply the identity to the first $N + 1$ terms with $\bb{r}_k = \bb{v}_k$ and $a_k = zd_k\pi_{0,k}/(1 - zd_k)$, and then multiply by $P(z)$. Since $\det(\bb{v}_i,\bb{v}_j) = y_j - y_i$, passage to the limit as $N \to \infty$, justified by the absolute convergence established below, gives exactly \eqref{eq:reduction.radial.entire.expansion}.

For every $R < \infty$, all products over subsets of the factors $1 - zd_k$ are uniformly bounded on the ball $\{z: |z| \leq R\}$ by $\exp(R\sum_{k=0}^{\infty}d_k)$, since $d_k = c_{\beta}\rho^k$ decays geometrically. The standard criterion for infinite products also shows that $P$, $P_k$, and $P_{i,j}$ converge locally uniformly. Because $(d/2)_k/k!$ grows at most polynomially in $k$, $d_k\pi_{0,k}$ decays as a polynomial factor multiplied by $\rho^{2k}$, and $y_k$ grows linearly in $k$, one has $\sum_{k=0}^{\infty}d_k\pi_{0,k}(1 + y_k^2) < \infty$. Since $(y_i - y_j)^2 \leq 2(1 + y_i^2)(1 + y_j^2)$, the terms of the double series in \eqref{eq:reduction.radial.entire.expansion} are bounded by products of this same summable sequence, so the double series is also absolutely summable. The Weierstrass test therefore proves local uniform convergence of both series in \eqref{eq:reduction.radial.entire.expansion}; hence that expression is entire. Moreover, the expression \eqref{eq:reduction.radial.entire.expansion} agrees with the Fredholm determinant away from the discrete set $\{d_k^{-1} : k \in \N_0\}$. The identity theorem proves \eqref{eq:reduction.radial.Fredholm.identity} everywhere and, in particular, proves that every apparent pole in \eqref{eq:results.radial.Fredholm.determinant} is removable. This concludes the proof of Proposition~\ref{prop:radial.Fredholm.determinant}. \qed

\subsection{Proof of Lemma~\ref{lem:one.coordinate.compression}}\label{subsec:proof.compression}

Let $X = D^{1/2}P$. Then $XX^* = D^{1/2}PD^{1/2}$ and $X^*X = PDP$. The standard correspondence $f \mapsto X^*f$ and $g \mapsto Xg$ between eigenvectors shows directly that $XX^*$ and $X^*X$ have the same nonzero eigenvalues with the same multiplicities; see also \citet[proof of Corollary~3.4.3, p.~42]{Kostenko2019}. Since $PDP$ vanishes on the span of $u$ and leaves $u^{\perp}$ invariant, it remains to analyze its restriction to $u^{\perp}$.

Suppose that $PDPx = \lambda x$ with $x \perp u$ and $\lambda > 0$. Then $Px = x$ and
\[
(D - \lambda I)x = au
\]
for some scalar $a$. If $\lambda = \delta_m$, the $m$-th coordinate forces $a = 0$ because $u_m \neq 0$, and then $x$ is a multiple of the $m$-th coordinate vector, contrary to $x \perp u$. Thus $\lambda \neq \delta_m$ for every $m$, $a \neq 0$, and
\[
x_k = \frac{au_k}{\delta_k - \lambda}, \qquad \sum_{k=0}^{\infty}\frac{u_k^2}{\delta_k - \lambda} = 0.
\]
Conversely, every positive solution $\lambda$ of the scalar equation defines an element of $\ell^2(\N_0)$ by the displayed coordinate formula, since the numbers $|\delta_k - \lambda|$ are bounded away from zero. The scalar equation gives $x \perp u$, and the coordinate equation then gives $PDPx = \lambda x$. Set
\[
F(\lambda) \leqdef \sum_{k=0}^{\infty}\frac{u_k^2}{\delta_k - \lambda}.
\]
On every compact subset of $(0,\infty) \setminus \{\delta_k : k \in \N_0\}$, the series defining $F$ and its derivative converge uniformly since $\sum_{k=0}^{\infty}u_k^2 = 1$ and the corresponding denominators are uniformly bounded away from zero for all sufficiently large $k$. Hence
\[
F'(\lambda) = \sum_{k=0}^{\infty}\frac{u_k^2}{(\delta_k - \lambda)^2} > 0.
\]
On $(\delta_{m+1},\delta_m)$, separating the singular summand at each endpoint shows that $F$ tends to negative infinity at the left endpoint and to positive infinity at the right endpoint. It therefore has exactly one zero in that interval. For $\lambda > \delta_0$, every summand is negative, so there is no zero. The displayed intervals cover every positive number below $\delta_0$ except the points $\delta_m$, which have already been excluded. Finally, $PDP$ is compact and self-adjoint, so every nonzero spectral value is an eigenvalue and the preceding analysis exhausts the nonzero spectrum. The coordinate formula also shows that each eigenspace is one-dimensional. This proves the lemma. \qed

\subsection{Proof of Theorem~\ref{thm:complete.spectrum}}\label{subsec:proof.main.theorem}

The orthogonal decomposition in Proposition~\ref{prop:spherical.Mercer.expansion} reduces the proof to the individual angular sectors. For every $\ell \geq 3$ such that $h_{d,\ell} > 0$, Proposition~\ref{prop:exact.block.reduction} and \eqref{eq:reduction.spherical.Mercer.eigenvalues} give
\[
A_{\beta,d}\psi_{\ell,k,j} = B_{\beta,d}\psi_{\ell,k,j} = c_{\beta}q^{\ell+2k}\psi_{\ell,k,j}, \qquad k \in \N_0, \quad 1 \leq j \leq h_{d,\ell}.
\]
The $h_{d,\ell}$ functions obtained by varying $j$ are orthogonal. This proves part \textrm{(i)}, including the stated angular multiplicity.

Fix $\ell = 1$ or, when $d \geq 2$, $\ell = 2$, and an angular index $1 \leq j \leq h_{d,\ell}$. This is possible because $h_{d,1} = d$ and $h_{d,2} = (d - 1)(d + 2)/2$, so a degree-$2$ sector exists exactly when $d \geq 2$. Consider the restriction of $A_{\beta,d}$ to $\mathcal{V}_{\ell,j}$, as defined in Proposition~\ref{prop:exact.block.reduction}. Applying Lemma~\ref{lem:one.coordinate.compression} to \eqref{eq:reduction.rank.one.matrix} with $\delta_k = \Lambda_{\ell,k}$ shows that none of the numbers $\Lambda_{\ell,k}$ is an eigenvalue. For the corresponding compression $P_{\ell}D_{\ell}P_{\ell}$ on $u_{\ell}^{\perp}$, the eigenvector equation has the form $(D_{\ell} - \lambda I)a = \alpha u_{\ell}$ with $\alpha \neq 0$. Substituting $a_k = \alpha\sqrt{\pi_{\ell,k}}/(\Lambda_{\ell,k} - \lambda)$ into $\langle a,u_{\ell}\rangle = 0$ shows that $\lambda = c_{\beta}q^{\ell}x$ is an eigenvalue if and only if
\[
0 = \sum_{k=0}^{\infty}\frac{\pi_{\ell,k}}{\Lambda_{\ell,k} - \lambda} = \frac{(1 - \rho)^{d/2+\ell}}{c_{\beta}q^{\ell}}\mathcal{Q}_{d/2+\ell,0}(x).
\]
Thus the eigenvalue equation is precisely \eqref{eq:results.degree.1.equation} when $\ell = 1$ and \eqref{eq:results.degree.2.equation} when $\ell = 2$. Termwise differentiation on every interval $(\rho^{m+1},\rho^m)$ gives
\[
\frac{\rd }{\rd x}\mathcal{Q}_{d/2+\ell,0}(x) = \sum_{k=0}^{\infty}\frac{(\frac{1}{2}d + \ell)_k}{k!}\frac{\rho^k}{(\rho^k - x)^2} > 0.
\]
The defining series and its derivative converge locally uniformly away from the poles. Moreover,
\[
\lim_{x\downarrow\rho^{m+1}}\mathcal{Q}_{d/2+\ell,0}(x) = -\infty, \qquad \lim_{x\uparrow\rho^m}\mathcal{Q}_{d/2+\ell,0}(x) = +\infty,
\]
because the terms with $k = m + 1$ and $k = m$ are the only singular terms at the left and right endpoints, respectively. Hence there is exactly one root in each interval $(\rho^{m+1},\rho^m)$. Lemma~\ref{lem:one.coordinate.compression} also shows that the corresponding eigenvalue is simple within each angular copy. For $\ell = 1$, the identical copies are indexed by $1 \leq j \leq h_{d,1} = d$, which proves part \textrm{(ii)}, including its multiplicity. For $d \geq 2$ and $\ell = 2$, they are indexed by $1 \leq j \leq h_{d,2} = (d - 1)(d + 2)/2$, which proves part \textrm{(iii)}. When $d = 1$, one has $h_{1,2} = 0$, so there is no degree-$2$ sector.

It remains to analyze the radial sector. Write
\[
D \leqdef D_0 = \operatorname{diag}(d_0,d_1,\ldots), \qquad \bb{b}_k \leqdef \{d_k\pi_{0,k}\}^{1/2}\bb{v}_k, \qquad (W\bb{\gamma})_k \leqdef \bb{b}_k^{\top}\bb{\gamma}.
\]
Since $W = D^{1/2}U_0$, \eqref{eq:reduction.radial.matrix} gives
\[
D - WW^* = D^{1/2}(I - U_0U_0^*)D^{1/2} = A_{\beta,d}^{(0)}.
\]
Moreover, $0 \leq A_{\beta,d}^{(0)} \leq D \leq c_{\beta}I$. Thus every nonzero radial eigenvalue satisfies $0 < \lambda \leq c_{\beta}$. Suppose first that $\lambda = c_{\beta}x$ and $x \neq \rho^k$ for every $k$. Since $\rho^0 = 1$, one has $0 < x < 1$, so $x$ belongs to a unique interval $(\rho^{m+1},\rho^m)$. If $a$ is an eigenvector and $\bb{\gamma} \leqdef W^*a$, the coordinate equations give
\begin{equation}\label{eq:proofs.radial.nonpole.coordinates}
a_k = \frac{\bb{b}_k^{\top}\bb{\gamma}}{d_k - \lambda}, \qquad \left[I_2 - \sum_{k=0}^{\infty}\frac{d_k\pi_{0,k}}{d_k - \lambda}\bb{v}_k\bb{v}_k^{\top}\right]\bb{\gamma} = \bb{0}.
\end{equation}
The vector $\bb{\gamma}$ cannot be zero, because otherwise $(D - \lambda I)a = 0$ at a non-pole value. By \eqref{eq:proofs.radial.isometry.moments} and the identity $\rho^k/(\rho^k - x) = 1 + x/(\rho^k - x)$, the matrix in brackets in \eqref{eq:proofs.radial.nonpole.coordinates} is
\begin{equation}\label{eq:proofs.radial.nonpole.matrix}
\begin{aligned}
I_2 - \sum_{k=0}^{\infty}\frac{d_k\pi_{0,k}}{d_k - \lambda}\bb{v}_k\bb{v}_k^{\top}
&= I_2 - \sum_{k=0}^{\infty}\frac{\rho^k\pi_{0,k}}{\rho^k - x}\bb{v}_k\bb{v}_k^{\top} \\
&= I_2 - \sum_{k=0}^{\infty}\pi_{0,k}\bb{v}_k\bb{v}_k^{\top} - x\sum_{k=0}^{\infty}\frac{\pi_{0,k}}{\rho^k - x}\bb{v}_k\bb{v}_k^{\top} \\
&= -x\sum_{k=0}^{\infty}\frac{\pi_{0,k}}{\rho^k - x}\bb{v}_k\bb{v}_k^{\top}.
\end{aligned}
\end{equation}
Set
\[
a_0 \leqdef \frac{dq}{\sqrt{2d}}, \qquad b_0 \leqdef \frac{2(1 - \rho)}{q\sqrt{2d}},
\]
so that $y_k = a_0 - b_0k$. With
\[
\Upsilon \leqdef \begin{pmatrix}1&0\\a_0&-b_0\end{pmatrix},
\]
the matrix without the factor $-x$ in \eqref{eq:proofs.radial.nonpole.matrix} equals
\[
(1 - \rho)^{d/2}\Upsilon\begin{pmatrix}\mathcal{Q}_{d/2,0}(x)&\mathcal{Q}_{d/2,1}(x)\\\mathcal{Q}_{d/2,1}(x)&\mathcal{Q}_{d/2,2}(x)\end{pmatrix}\Upsilon^{\top}.
\]
Since $\det(\Upsilon) = -b_0 \neq 0$, the determinant of this matrix vanishes if and only if \eqref{eq:results.radial.nonpole.equation} holds. Conversely, suppose that $x \in (0,1) \setminus \{\rho^k : k \in \N_0\}$ satisfies \eqref{eq:results.radial.nonpole.equation}, and choose a nonzero null vector $\bb{\gamma}$ for the matrix in \eqref{eq:proofs.radial.nonpole.matrix}. Define $a$ by the first equation in \eqref{eq:proofs.radial.nonpole.coordinates}. Since $\lambda > 0$, $d_k \to 0$, and $\lambda \neq d_k$ for every $k$, one has $\inf_{k \in \N_0}|d_k - \lambda| > 0$. Also, \eqref{eq:proofs.radial.isometry.moments} gives
\[
\sum_{k=0}^{\infty}\|\bb{b}_k\|^2 = \sum_{k=0}^{\infty}d_k\pi_{0,k}(1 + y_k^2) \leq 2d_0.
\]
It follows that $a \in \ell^2(\N_0)$. The null-vector equation and the definition of $a$ give
\[
W^*a = \sum_{k=0}^{\infty}\frac{\bb{b}_k\bb{b}_k^{\top}}{d_k - \lambda}\bb{\gamma} = \bb{\gamma}.
\]
Thus $(D - \lambda I)a = W\bb{\gamma} = WW^*a$, so $(D - WW^*)a = \lambda a$. This proves necessity and sufficiency away from the poles.

It remains to show that the null space of the matrix in \eqref{eq:proofs.radial.nonpole.matrix} is one-dimensional. Let $m$ be the unique index such that $\rho^{m+1} < x < \rho^m$. If the matrix in \eqref{eq:proofs.radial.nonpole.matrix} were zero, then the invertibility of $\Upsilon$ would imply $\mathcal{Q}_{d/2,0}(x) = \mathcal{Q}_{d/2,1}(x) = 0$. Let
\[
\varrho_k \leqdef \frac{(d/2)_k\rho^k}{k! \, |\rho^k - x|}.
\]
The first equation, $\mathcal{Q}_{d/2,0}(x) = 0$, means that $\sum_{k\leq m}\varrho_k = \sum_{k\geq m+1}\varrho_k \leqdef R > 0$. The second equality would then yield
\[
0 = \mathcal{Q}_{d/2,1}(x) = \sum_{k\leq m}k\varrho_k - \sum_{k\geq m+1}k\varrho_k \leq mR - (m + 1)R < 0,
\]
which is contradictory. Hence the singular $2 \times 2$ matrix in \eqref{eq:proofs.radial.nonpole.matrix} is not the zero matrix, so its null space is one-dimensional.

Now let $\lambda = d_m = c_{\beta}\rho^m$. From the $m$-th coordinate of the eigenvalue equation,
\[
\bb{b}_m^{\top}\bb{\gamma} = 0, \qquad \bb{\gamma} \leqdef W^*a.
\]
The vector $\bb{\gamma}$ cannot vanish: otherwise $a$ would be a multiple of the $m$-th coordinate vector, whose image under $W^*$ is nonzero. Hence $\bb{\gamma}$ is a nonzero multiple of $\bb{w}_m = (y_m,-1)^{\top}$. For $k \neq m$,
\[
a_k = \frac{\bb{b}_k^{\top}\bb{\gamma}}{d_k - d_m}.
\]
Define
\[
S_m \leqdef I_2 - \sum_{k\neq m}\frac{\bb{b}_k\bb{b}_k^{\top}}{d_k - d_m}.
\]
Then the remaining equation $\bb{\gamma} = W^*a$ is
\begin{equation}\label{eq:proofs.radial.pole.vector.equation}
S_m \bb{\gamma} = \bb{b}_ma_m.
\end{equation}
The matrix $S_m$ is well defined because $\inf_{k\neq m}|d_k - d_m| > 0$ and $\sum_{k=0}^{\infty}\|\bb{b}_k\|^2 < \infty$. Taking the inner product of \eqref{eq:proofs.radial.pole.vector.equation} with $\bb{w}_m$ and using $\bb{w}_m^{\top}\bb{b}_m = 0$ shows that the necessary scalar condition is $\bb{w}_m^{\top}S_m\bb{w}_m = 0$. Moreover, \eqref{eq:proofs.radial.isometry.moments} gives $\|\bb{w}_m\|^2 = \sum_{k=0}^{\infty}\pi_{0,k}(\bb{w}_m^{\top}\bb{v}_k)^2$. Since $\bb{w}_m^{\top}\bb{v}_k = y_m - y_k$ and the term with $k = m$ vanishes, direct subtraction gives
\begin{equation}\label{eq:proofs.radial.pole.scalar}
\bb{w}_m^{\top}S_m\bb{w}_m
= \sum_{k\neq m}\pi_{0,k}(y_m - y_k)^2 - \sum_{k\neq m}\frac{\rho^k\pi_{0,k}(y_m - y_k)^2}{\rho^k - \rho^m}
= -\rho^m\sum_{k\neq m}\frac{\pi_{0,k}(y_m - y_k)^2}{\rho^k - \rho^m}.
\end{equation}
The identity
\[
y_m - y_k = \frac{2(1 - \rho)(k - m)}{q\sqrt{2d}},
\]
together with the formula for $\pi_{0,k}$ in \eqref{eq:definitions.Gaussian.levels.and.weights}, shows that the right-hand side of \eqref{eq:proofs.radial.pole.scalar} vanishes exactly when \eqref{eq:results.radial.pole.equation} holds. This proves necessity. Conversely, suppose that \eqref{eq:results.radial.pole.equation} holds and take $\bb{\gamma} = \bb{w}_m$. Then $S_m\bb{w}_m$ is orthogonal to $\bb{w}_m$, so
\[
S_m\bb{w}_m \in \bb{w}_m^{\perp} = \operatorname{span}\{\bb{v}_m\} = \operatorname{span}\{\bb{b}_m\}.
\]
Consequently, the unique choice
\[
a_m = \frac{\bb{b}_m^{\top}S_m\bb{w}_m}{\|\bb{b}_m\|^2}
\]
satisfies $S_m\bb{w}_m = \bb{b}_ma_m$. For $k \neq m$, define $a_k = \bb{b}_k^{\top}\bb{w}_m/(d_k - d_m)$. The two bounds used to define $S_m$ give
\[
\sum_{k\neq m}|a_k|^2 \leq \frac{\|\bb{w}_m\|^2}{\inf_{k\neq m}|d_k - d_m|^2}\sum_{k\neq m}\|\bb{b}_k\|^2 < \infty.
\]
The identity $S_m\bb{w}_m = \bb{b}_ma_m$ is equivalent to $W^*a = \bb{w}_m$. The coordinate equations therefore give $(D - WW^*)a = d_ma$, which proves sufficiency. Every pole eigenvector has $\bb{\gamma}$ in the one-dimensional space $\operatorname{span}\{\bb{w}_m\}$, after which all coordinates of $a$ are determined by \eqref{eq:proofs.radial.pole.vector.equation}. Hence every pole eigenspace is one-dimensional.

By \eqref{eq:reduction.radial.Fredholm.identity}, $\mathcal{D}_0(z) = \det\nolimits_F(I - zA_{\beta,d}^{(0)})$. The invertibility and multiplicity theorem for Fredholm determinants \citep[Theorem~3.4.6, pp.~40--41]{Kostenko2019} therefore shows that $z \neq 0$ is a zero of $\mathcal{D}_0$ if and only if $z^{-1}$ is a nonzero radial eigenvalue, and the order of the zero equals the eigenvalue multiplicity. Since the radial restriction is nonnegative and self-adjoint, its nonzero eigenvalues, and therefore the zeros of $\mathcal{D}_0$, are positive. This proves all the radial characterizations in part \textrm{(iv)} except the strict inequalities, which are established next.

Define
\[
u_0 \leqdef (\sqrt{\pi_{0,k}})_{k \in \N_0}, \qquad u_1 \leqdef (\sqrt{\pi_{0,k}}y_k)_{k \in \N_0};
\]
then it follows from \eqref{eq:proofs.radial.isometry.moments} that $u_0$ and $u_1$ are orthonormal. Let $P_0 \leqdef I - u_0 \otimes u_0$, and let $C_0$ be the restriction of $P_0DP_0$ to $u_0^{\perp}$. By Lemma~\ref{lem:one.coordinate.compression}, $C_0$ has simple eigenvalues $(\eta_m)_{m \in \N_0}$ satisfying
\[
d_{m+1} < \eta_m < d_m.
\]
Every corresponding eigenvector $x^{(m)}$ has coordinates proportional to $\sqrt{\pi_{0,k}}/(d_k - \eta_m)$, and its inner product with $u_1$ is nonzero. Indeed, orthogonality to $u_0$ gives
\[
\sum_{k=0}^{\infty}\frac{\pi_{0,k}}{d_k - \eta_m} = 0.
\]
If the inner product with $u_1$ also vanished, then the affine term $y_k = a_0 - b_0k$ would yield
\[
0 = \sum_{k=0}^{\infty}\frac{\pi_{0,k}y_k}{d_k - \eta_m}
= a_0\sum_{k=0}^{\infty}\frac{\pi_{0,k}}{d_k - \eta_m} - b_0\sum_{k=0}^{\infty}\frac{k\pi_{0,k}}{d_k - \eta_m}.
\]
As we showed earlier, the first sum equals zero, so the second sum would also vanish. Since $d_{m+1} < \eta_m < d_m$, the terms $\varepsilon_k \leqdef \pi_{0,k}/(d_k - \eta_m)$ are positive for $k \leq m$ and negative for $k \geq m + 1$. Their total sum is zero, so the total positive mass and the absolute negative mass are equal to some $R > 0$. It follows that
\[
\sum_{k=0}^{\infty}k\varepsilon_k \leq mR - (m + 1)R < 0,
\]
which is a contradiction.

The operator $C_0$ is nonnegative and compact. It is also injective, because for $x \in u_0^{\perp}$,
\[
\langle C_0x,x\rangle = \langle Dx,x\rangle = \sum_{k=0}^{\infty}d_k|x_k|^2,
\]
which is positive when $x \neq 0$. The compact self-adjoint spectral theorem \citep[Theorem~2.3.2, p.~12]{Kostenko2019} therefore shows that the eigenvectors of $C_0$ form an orthonormal basis of $u_0^{\perp}$. In this basis, $C_0$ is represented by $\operatorname{diag}(\eta_0,\eta_1,\ldots)$, and $u_1 \in u_0^{\perp}$ has a nonzero coordinate along every eigenvector. A second application of Lemma~\ref{lem:one.coordinate.compression} therefore shows that the compression of $C_0$ to $\{u_0,u_1\}^{\perp}$ has simple eigenvalues satisfying
\[
\eta_{m+1} < \lambda_{0,m} < \eta_m.
\]
Since $u_0$ and $u_1$ are orthonormal,
\[
Q_0 = I - u_0 \otimes u_0 - u_1 \otimes u_1
\]
is the projection onto $\{u_0,u_1\}^{\perp}$. Since $Q_0P_0 = P_0Q_0 = Q_0$, the compression of $C_0$ to this space is $Q_0P_0DP_0Q_0 = Q_0DQ_0$. Taking $X = D^{1/2}Q_0$ gives $XX^* = D^{1/2}Q_0D^{1/2}$ and $X^*X = Q_0DQ_0$. The eigenvector correspondence used in Lemma~\ref{lem:one.coordinate.compression} therefore proves equality of their nonzero spectra. Moreover, $\lambda_{0,m+1} < \eta_{m+1} < \lambda_{0,m}$, so $(\lambda_{0,m})_{m \in \N_0}$ is in decreasing order. Combining the bounds for $\eta_m$ and $\lambda_{0,m}$ gives
\[
c_{\beta}\rho^{m+2} = d_{m+2} < \lambda_{0,m} < d_m = c_{\beta}\rho^m,
\]
which proves \eqref{eq:results.radial.interlacing} and radial simplicity. Together with the Fredholm determinant, non-pole, and pole characterizations proved above, this completes the proof of part \textrm{(iv)}.

It remains to prove part \textrm{(v)}. Proposition~\ref{prop:positivity.trace.injectivity} shows that $A_{\beta,d}$ is nonnegative, trace class, and injective. In particular, zero is not an eigenvalue, and trace class implies compactness. The compact self-adjoint spectral theorem \citep[Theorem~2.3.2, p.~12]{Kostenko2019}, together with the orthogonal angular decomposition, shows that the eigenvalues found in all the sectors exhaust the nonzero spectrum. Moreover, part \textrm{(ii)} gives
\[
0 < \lambda_{1,m} < c_{\beta}q\rho^m \to 0,
\]
so zero is an accumulation point of the nonzero spectrum. The angular multiplicities have already been accounted for, and the multiplicities of eigenvalues belonging to orthogonal sectors add when their numerical values coincide. This proves every assertion in part \textrm{(v)} and completes the proof of Theorem~\ref{thm:complete.spectrum}. \qed

\subsection{Proof of Theorem~\ref{thm:Henze.Zirkler.spectrum}}\label{subsec:proof.Henze.Zirkler.spectrum}

We first calculate the orthogonal projection that appears in the centered feature. The functions
\[
1, \qquad x_i \quad (1 \leq i \leq d), \qquad x_ix_j \quad (1 \leq i < j \leq d), \qquad \frac{x_i^2 - 1}{\sqrt{2}} \quad (1 \leq i \leq d)
\]
form an orthonormal basis of $\mathcal{P}_{\leq 2}$ in $L^2(\mu_1)$ by the second and fourth centered Gaussian moment formulas; see \citet[Theorem~1, p.~2404]{VignatBhatnagar2008}. Therefore the integral kernel of $\Pi_{\leq 2}$ is
\[
p_{\leq 2}(\bb{x},\bb{y}) = 1 + \bb{x}^{\top}\bb{y} + \frac{1}{2}\{(\bb{x}^{\top}\bb{y})^2 - \|\bb{x}\|^2 - \|\bb{y}\|^2 + d\}.
\]
Indeed, the contribution from the quadratic basis functions is
\[
\frac{1}{2}\sum_{i=1}^d(x_i^2 - 1)(y_i^2 - 1) + \sum_{1 \leq i<j \leq d}x_ix_jy_iy_j = \frac{1}{2}\{(\bb{x}^{\top}\bb{y})^2 - \|\bb{x}\|^2 - \|\bb{y}\|^2 + d\}.
\]

Let $\bb{X}$ have distribution $\mu_1$. The Gaussian Fourier transform \citep[Proposition~8.24, p.~251]{Folland1999}, after rescaling to the present convention, gives
\begin{equation}\label{eq:proofs.HZ.standard.Gaussian.Fourier}
\EE\{\exp(\ii\bb{t}^{\top}\bb{X})\} = \exp\left(-\frac{\|\bb{t}\|^2}{2}\right).
\end{equation}
Differentiating with respect to the coordinates of $\bb{t}$, and then taking real and imaginary parts, gives
\begin{align}
\langle \Phi_{\bb{t}},1\rangle_1 &= \exp\left(-\frac{\|\bb{t}\|^2}{2}\right),\label{eq:proofs.HZ.projection.constant}\\
\langle \Phi_{\bb{t}},x_i\rangle_1 &= t_i\exp\left(-\frac{\|\bb{t}\|^2}{2}\right),\label{eq:proofs.HZ.projection.linear}\\
\langle \Phi_{\bb{t}},x_ix_j - \delta_{ij}\rangle_1 &= -t_it_j\exp\left(-\frac{\|\bb{t}\|^2}{2}\right).\label{eq:proofs.HZ.projection.quadratic}
\end{align}
Using the orthonormal basis above in \eqref{eq:proofs.HZ.projection.constant}--\eqref{eq:proofs.HZ.projection.quadratic} yields
\begin{equation}\label{eq:proofs.HZ.projected.feature}
(\Pi_{\leq 2}\Phi_{\bb{t}})(\bb{x}) = \exp\left(-\frac{\|\bb{t}\|^2}{2}\right)\left[1 + \bb{t}^{\top}\bb{x} - \frac{1}{2}\{(\bb{t}^{\top}\bb{x})^2 - \|\bb{t}\|^2\}\right].
\end{equation}
Comparison with \eqref{eq:definitions.centered.feature} proves
\begin{equation}\label{eq:proofs.HZ.feature.is.projection}
\zeta(\bb{x},\bb{t}) = (\Pi_{\geq 3}\Phi_{\bb{t}})(\bb{x}).
\end{equation}
Since $(\cos u + \sin u)^2 = 1 + \sin(2u)$ and the distribution of $\bb{t}^{\top}\bb{X}$ is symmetric,
\[
\|\Phi_{\bb{t}}\|_1^2 = \EE\{1 + \sin(2\bb{t}^{\top}\bb{X})\} = 1.
\]
An orthogonal projection does not enlarge the norm, so \eqref{eq:proofs.HZ.feature.is.projection} gives
\begin{equation}\label{eq:proofs.HZ.feature.square.integrability}
\int_{\R^d}\int_{\R^d}|\zeta(\bb{x},\bb{t})|^2\varphi_1(\bb{x})\varphi_{\beta}(\bb{t}) \, \rd \bb{x} \, \rd \bb{t} \leq 1.
\end{equation}
Thus $\mathcal{X}_{\beta,d}$ is Hilbert--Schmidt by the kernel characterization \citep[Eq.~3.3.9, p.~28]{Kostenko2019}. Its adjoint, with respect to $\langle\cdot,\cdot\rangle_1$ in the domain and $\langle\cdot,\cdot\rangle_{\beta}$ in the codomain, is
\[
(\mathcal{X}_{\beta,d}^*g)(\bb{x}) = \int_{\R^d}\zeta(\bb{x},\bb{t})g(\bb{t})\varphi_{\beta}(\bb{t}) \, \rd \bb{t}.
\]

We next verify the long kernel $h_{\beta,d}^*$ term by term. Let $\bb{G}$ have distribution $\mu_{\beta}$, and let
\begin{equation}\label{eq:proofs.HZ.polynomial.factor}
r(\bb{x},\bb{t}) \leqdef 1 + \bb{t}^{\top}\bb{x} - \frac{1}{2}\{(\bb{t}^{\top}\bb{x})^2 - \|\bb{t}\|^2\}.
\end{equation}
For $a \geq 0$, multiplying the Gaussian density by $\exp(-a\|\bb{t}\|^2/2)$ gives
\[
\exp\left(-\frac{a}{2}\|\bb{t}\|^2\right)\varphi_{\beta}(\bb{t}) = (1 + a\beta^2)^{-d/2}\varphi_{\beta/\sqrt{1+a\beta^2}}(\bb{t}).
\]
Multiplying by $H(\bb{t})$ and integrating proves that, for every measurable function $H$ for which the expectations below are absolutely finite,
\begin{equation}\label{eq:proofs.HZ.tilting.identity}
\EE\left\{\exp\left(-\frac{a}{2}\|\bb{G}\|^2\right)H(\bb{G})\right\} = (1 + a\beta^2)^{-d/2}\EE\{H(\bb{Z}_a)\}, \qquad \bb{Z}_a \sim N_d\left(\bb{0},\frac{\beta^2}{1 + a\beta^2}I_d\right).
\end{equation}

Let $\bb{Z}$ be centered Gaussian with covariance matrix $\tau I_d$. Differentiation of the scaled form of \eqref{eq:proofs.HZ.standard.Gaussian.Fourier}, now with covariance matrix $\tau I_d$, gives
\begin{align}
\EE\{\cos(\bb{x}^{\top}\bb{Z})\} &= \exp\left(-\frac{\tau}{2}\|\bb{x}\|^2\right),\label{eq:proofs.HZ.Gaussian.moment.1}\\
\EE\{(\bb{y}^{\top}\bb{Z})\sin(\bb{x}^{\top}\bb{Z})\} &= \tau\bb{x}^{\top}\bb{y}\exp\left(-\frac{\tau}{2}\|\bb{x}\|^2\right),\label{eq:proofs.HZ.Gaussian.moment.2}\\
\EE\{(\bb{y}^{\top}\bb{Z})^2\cos(\bb{x}^{\top}\bb{Z})\} &= \{\tau\|\bb{y}\|^2 - \tau^2(\bb{x}^{\top}\bb{y})^2\}\exp\left(-\frac{\tau}{2}\|\bb{x}\|^2\right),\label{eq:proofs.HZ.Gaussian.moment.3}\\
\EE\{\|\bb{Z}\|^2\cos(\bb{x}^{\top}\bb{Z})\} &= \{d\tau - \tau^2\|\bb{x}\|^2\}\exp\left(-\frac{\tau}{2}\|\bb{x}\|^2\right).\label{eq:proofs.HZ.Gaussian.moment.4}
\end{align}
All terms with an odd integrand have expectation zero. Apply \eqref{eq:proofs.HZ.tilting.identity} with $a = 1$, so that the covariance of $\bb{Z}_1$ is $\tau_{\beta,1}I_d$ and $\tau_{\beta,1}/2 = \gamma_{\beta}$. Equations \eqref{eq:proofs.HZ.Gaussian.moment.1}--\eqref{eq:proofs.HZ.Gaussian.moment.4} give
\begin{align}
J_{01}(\bb{x},\bb{y}) &\leqdef \EE\left[\Phi_{\bb{G}}(\bb{x})\exp\left(-\frac{\|\bb{G}\|^2}{2}\right)r(\bb{y},\bb{G})\right]\nonumber\\
&= (1 + \beta^2)^{-d/2}\exp(-\gamma_{\beta}\|\bb{x}\|^2)\left[1 + \tau_{\beta,1}\bb{x}^{\top}\bb{y} \right. \nonumber\\
&\qquad \left. - \frac{1}{2}\{\tau_{\beta,1}\|\bb{y}\|^2 - \tau_{\beta,1}^2(\bb{x}^{\top}\bb{y})^2 - d\tau_{\beta,1} + \tau_{\beta,1}^2\|\bb{x}\|^2\}\right]\nonumber\\
&= (1 + \beta^2)^{-d/2}\exp(-\gamma_{\beta}\|\bb{x}\|^2)\left[1 + \gamma_{\beta}\left\{\tau_{\beta,1}(\bb{x}^{\top}\bb{y})^2 - \tau_{\beta,1}\|\bb{x}\|^2 - \|\bb{y}\|^2 + 2\bb{x}^{\top}\bb{y} + d\right\}\right].\label{eq:proofs.HZ.first.cross.term}
\end{align}
The same calculation with $\bb{x}$ and $\bb{y}$ interchanged gives $J_{10}(\bb{x},\bb{y}) = J_{01}(\bb{y},\bb{x})$, which is the second subtracted term in \eqref{eq:definitions.HZ.kernel}.

For the last term, let $\bb{Z}$ have covariance matrix $\tau_{\beta,2}I_d$. The fourth centered Gaussian moment identity, which is the fourth-order case of Wick's theorem \citep[Theorem~1, p.~2404]{VignatBhatnagar2008}, is
\begin{equation}\label{eq:proofs.HZ.fourth.moment.identity}
\EE(Z_iZ_jZ_kZ_l) = \tau_{\beta,2}^2(\delta_{ij}\delta_{kl} + \delta_{ik}\delta_{jl} + \delta_{il}\delta_{jk}).
\end{equation}
Summing \eqref{eq:proofs.HZ.fourth.moment.identity} over the coordinates gives
\begin{align}
\EE\{(\bb{x}^{\top}\bb{Z})(\bb{y}^{\top}\bb{Z})\} &= \tau_{\beta,2}\bb{x}^{\top}\bb{y},\label{eq:proofs.HZ.Wick.1}\\
\EE\{(\bb{x}^{\top}\bb{Z})^2(\bb{y}^{\top}\bb{Z})^2\} &= \tau_{\beta,2}^2\{\|\bb{x}\|^2\|\bb{y}\|^2 + 2(\bb{x}^{\top}\bb{y})^2\},\label{eq:proofs.HZ.Wick.2}\\
\EE\{(\bb{x}^{\top}\bb{Z})^2\|\bb{Z}\|^2\} &= \tau_{\beta,2}^2(d + 2)\|\bb{x}\|^2,\label{eq:proofs.HZ.Wick.3}\\
\EE\{\|\bb{Z}\|^4\} &= \tau_{\beta,2}^2d(d + 2).\label{eq:proofs.HZ.Wick.4}
\end{align}
It follows from \eqref{eq:proofs.HZ.Wick.2}--\eqref{eq:proofs.HZ.Wick.4} that
\begin{equation}\label{eq:proofs.HZ.quadratic.product.moment}
\begin{aligned}
&\EE\left[\{(\bb{x}^{\top}\bb{Z})^2 - \|\bb{Z}\|^2\}\{(\bb{y}^{\top}\bb{Z})^2 - \|\bb{Z}\|^2\}\right] \\
&\qquad= \tau_{\beta,2}^2\left[\{\|\bb{x}\|^2 - d\}\{\|\bb{y}\|^2 - d\} + 2\{(\bb{x}^{\top}\bb{y})^2 - \|\bb{x}\|^2 - \|\bb{y}\|^2 + d\}\right].
\end{aligned}
\end{equation}
Apply \eqref{eq:proofs.HZ.tilting.identity} with $a = 2$ and expand the two factors $r$ in \eqref{eq:proofs.HZ.polynomial.factor}. The products of a linear term and a constant or quadratic term have zero expectation. Equations \eqref{eq:proofs.HZ.Wick.1} and \eqref{eq:proofs.HZ.quadratic.product.moment} therefore give
\begin{align}
J_{11}(\bb{x},\bb{y}) &\leqdef \EE\{e^{-\|\bb{G}\|^2}r(\bb{x},\bb{G})r(\bb{y},\bb{G})\}\nonumber\\
&= (1 + 2 \beta^2)^{-d/2}\left[1 - \frac{\tau_{\beta,2}}{2}\{\|\bb{x}\|^2 + \|\bb{y}\|^2 - 2d - 2\bb{x}^{\top}\bb{y}\} \right.\nonumber\\
&\qquad \left. + \frac{\tau_{\beta,2}^2}{4}\left[\{\|\bb{x}\|^2 - d\}\{\|\bb{y}\|^2 - d\} + 2\{(\bb{x}^{\top}\bb{y})^2 - \|\bb{x}\|^2 - \|\bb{y}\|^2 + d\}\right]\right].\label{eq:proofs.HZ.last.term}
\end{align}

The unprojected term is simpler. The elementary identity
\[
\Phi_{\bb{G}}(\bb{x})\Phi_{\bb{G}}(\bb{y}) = \cos\{\bb{G}^{\top}(\bb{x} - \bb{y})\} + \sin\{\bb{G}^{\top}(\bb{x} + \bb{y})\}
\]
and the symmetry of $\mu_{\beta}$, together with the scaled form of \eqref{eq:proofs.HZ.standard.Gaussian.Fourier}, show that
\begin{equation}\label{eq:proofs.HZ.unprojected.feature.product}
J_{00}(\bb{x},\bb{y}) \leqdef \EE\{\Phi_{\bb{G}}(\bb{x})\Phi_{\bb{G}}(\bb{y})\} = \exp\left(-\frac{\beta^2}{2}\|\bb{x} - \bb{y}\|^2\right).
\end{equation}
By \eqref{eq:definitions.centered.feature},
\[
\EE\{\zeta(\bb{x},\bb{G})\zeta(\bb{y},\bb{G})\} = J_{00}(\bb{x},\bb{y}) - J_{01}(\bb{x},\bb{y}) - J_{10}(\bb{x},\bb{y}) + J_{11}(\bb{x},\bb{y}).
\]
Comparison of \eqref{eq:proofs.HZ.first.cross.term}, its transposed version, \eqref{eq:proofs.HZ.last.term}, and \eqref{eq:proofs.HZ.unprojected.feature.product} with \eqref{eq:definitions.HZ.kernel} proves the exact identity
\begin{equation}\label{eq:proofs.HZ.kernel.feature.identity}
h_{\beta,d}^*(\bb{x},\bb{y}) = \int_{\R^d}\zeta(\bb{x},\bb{t})\zeta(\bb{y},\bb{t})\varphi_{\beta}(\bb{t}) \, \rd \bb{t}.
\end{equation}
For $f,g \in L^2(\mu_1)$, Cauchy--Schwarz and \eqref{eq:proofs.HZ.feature.square.integrability} give
\[
\begin{aligned}
&\int_{\R^d}\left\{\int_{\R^d}|f(\bb{x})\zeta(\bb{x},\bb{t})|\varphi_1(\bb{x}) \, \rd \bb{x}\right\}
\left\{\int_{\R^d}|g(\bb{y})\zeta(\bb{y},\bb{t})|\varphi_1(\bb{y}) \, \rd \bb{y}\right\}\varphi_{\beta}(\bb{t}) \, \rd \bb{t} \\
&\qquad\leq \|f\|_1\|g\|_1\int_{\R^d}\|\zeta(\cdot,\bb{t})\|_1^2\varphi_{\beta}(\bb{t}) \, \rd \bb{t} < \infty.
\end{aligned}
\]
Thus Fubini's theorem applies to the associated bilinear forms, and \eqref{eq:proofs.HZ.kernel.feature.identity} gives
\begin{equation}\label{eq:proofs.HZ.XstarX}
\widetilde{A}_{\beta,d} = \mathcal{X}_{\beta,d}^*\mathcal{X}_{\beta,d}.
\end{equation}

We next identify the operator in the other order. From \eqref{eq:proofs.HZ.feature.is.projection},
\begin{equation}\label{eq:proofs.HZ.BHEP.kernel.start}
\int_{\R^d}\zeta(\bb{x},\bb{s})\zeta(\bb{x},\bb{t})\varphi_1(\bb{x}) \, \rd \bb{x}
= \langle\Pi_{\geq 3}\Phi_{\bb{s}},\Pi_{\geq 3}\Phi_{\bb{t}}\rangle_1
= \langle \Phi_{\bb{s}},\Phi_{\bb{t}}\rangle_1 - \langle\Pi_{\leq 2}\Phi_{\bb{s}},\Pi_{\leq 2}\Phi_{\bb{t}}\rangle_1.
\end{equation}
The first inner product on the right is $\exp(-\|\bb{s} - \bb{t}\|^2/2)$ by \eqref{eq:proofs.HZ.standard.Gaussian.Fourier}. The constant, linear, and quadratic parts of \eqref{eq:proofs.HZ.projected.feature} are mutually orthogonal. Moreover, the Wick formula in \eqref{eq:proofs.HZ.fourth.moment.identity}, with covariance matrix $I_d$, gives
\[
\EE[\{(\bb{s}^{\top}\bb{X})^2 - \|\bb{s}\|^2\}\{(\bb{t}^{\top}\bb{X})^2 - \|\bb{t}\|^2\}] = 2(\bb{s}^{\top}\bb{t})^2.
\]
It follows that
\[
\langle\Pi_{\leq 2}\Phi_{\bb{s}},\Pi_{\leq 2}\Phi_{\bb{t}}\rangle_1 = \left\{1 + \bb{s}^{\top}\bb{t} + \frac{(\bb{s}^{\top}\bb{t})^2}{2}\right\}\exp\left(-\frac{\|\bb{s}\|^2 + \|\bb{t}\|^2}{2}\right).
\]
The difference in \eqref{eq:proofs.HZ.BHEP.kernel.start} is exactly $K(\bb{s},\bb{t})$ in \eqref{eq:definitions.BHEP.kernel}. The same Cauchy--Schwarz argument as above, now applied to $f,g \in L^2(\mu_{\beta})$ with integration first with respect to $\mu_1$, shows by \eqref{eq:proofs.HZ.feature.square.integrability} that the corresponding integral is finite. Hence Fubini's theorem gives
\begin{equation}\label{eq:proofs.HZ.XXstar}
A_{\beta,d} = \mathcal{X}_{\beta,d}\mathcal{X}_{\beta,d}^*.
\end{equation}

There is also a direct compression interpretation of the first factorization. For $f,g \in L^2(\mu_1)$, \eqref{eq:proofs.HZ.feature.is.projection} and the self-adjointness of $\Pi_{\geq 3}$ give
\[
(\mathcal{X}_{\beta,d}f)(\bb{t}) = \langle f,\Pi_{\geq 3}\Phi_{\bb{t}}\rangle_1 = \langle\Pi_{\geq 3}f,\Phi_{\bb{t}}\rangle_1.
\]
Using \eqref{eq:proofs.HZ.unprojected.feature.product} and the Fubini estimate above, one obtains
\[
\langle f,\widetilde{A}_{\beta,d}g\rangle_1
= \langle\mathcal{X}_{\beta,d}f,\mathcal{X}_{\beta,d}g\rangle_{\beta}
= \langle\Pi_{\geq 3}f,\widetilde{B}_{\beta,d}\Pi_{\geq 3}g\rangle_1
= \langle f,\Pi_{\geq 3}\widetilde{B}_{\beta,d}\Pi_{\geq 3}g\rangle_1.
\]
Since this identity holds for every $f,g \in L^2(\mu_1)$,
\begin{equation}\label{eq:proofs.HZ.compression}
\widetilde{A}_{\beta,d} = \Pi_{\geq 3}\widetilde{B}_{\beta,d}\Pi_{\geq 3}.
\end{equation}
Together, \eqref{eq:proofs.HZ.XstarX}, \eqref{eq:proofs.HZ.XXstar}, and \eqref{eq:proofs.HZ.compression} prove \eqref{eq:results.HZ.factorizations}. They also prove that $\widetilde{A}_{\beta,d}$ is nonnegative and self-adjoint. Since $\mathcal{X}_{\beta,d}$ is Hilbert--Schmidt, the Hilbert--Schmidt product theorem \citep[Theorem~3.3.1\textrm{(vi)}, pp.~27--28]{Kostenko2019} shows that $\mathcal{X}_{\beta,d}^*\mathcal{X}_{\beta,d}$ is trace class. This proves every assertion in part \textrm{(i)}.

We now determine the null space. If $f \in \mathcal{P}_{\leq 2}$, then $\Pi_{\geq 3}f = 0$, so \eqref{eq:proofs.HZ.compression} gives $\widetilde{A}_{\beta,d}f = 0$. Conversely, suppose that $\widetilde{A}_{\beta,d}f = 0$ and let $h \leqdef \Pi_{\geq 3}f$. The Gaussian Fourier identity \citep[Proposition~8.24, p.~251]{Folland1999}, after rescaling to the present convention, is
\[
\exp\left(-\frac{\beta^2}{2}\|\bb{x} - \bb{y}\|^2\right) = \int_{\R^d}\exp\{\ii\bb{t}^{\top}(\bb{x} - \bb{y})\}\varphi_{\beta}(\bb{t}) \, \rd \bb{t}.
\]
Since $h \in L^2(\mu_1)$ and $\mu_1$ is a probability measure, $h \in L^1(\mu_1)$. The absolute triple integral arising from the Fourier identity is bounded by $\|h\|_{L^1(\mu_1)}^2$, so Fubini's theorem gives
\begin{equation}\label{eq:proofs.HZ.strict.positivity}
0
= \langle f,\widetilde{A}_{\beta,d}f\rangle_1
= \langle h,\widetilde{B}_{\beta,d}h\rangle_1
= \int_{\R^d}\left|\int_{\R^d}\exp(\ii\bb{t}^{\top}\bb{x})h(\bb{x})\varphi_1(\bb{x}) \, \rd \bb{x}\right|^2\varphi_{\beta}(\bb{t}) \, \rd \bb{t}.
\end{equation}
The function $h\varphi_1$ is integrable, and its Fourier transform in \eqref{eq:proofs.HZ.strict.positivity} is continuous by \citet[Theorem~8.22\textrm{(f)}, p.~249]{Folland1999}. Since $\varphi_{\beta}$ is strictly positive everywhere, the last integral in \eqref{eq:proofs.HZ.strict.positivity} can vanish only if this Fourier transform vanishes Lebesgue almost everywhere. Continuity then implies that it vanishes for every $\bb{t} \in \R^d$. Fourier-transform uniqueness \citep[Corollary~8.27, p.~252]{Folland1999} gives $h\varphi_1 = 0$ almost everywhere. Since $\varphi_1$ is strictly positive, $h = 0$ in $L^2(\mu_1)$. Thus $f \in \mathcal{P}_{\leq 2}$, and $\ker(\widetilde{A}_{\beta,d}) = \mathcal{P}_{\leq 2}$. The orthonormal basis at the beginning of the proof contains $1 + d + d(d + 1)/2$ functions, which proves \eqref{eq:results.HZ.null.space}. Since $\widetilde{A}_{\beta,d}$ is trace class, it is compact. Its restriction to the infinite-dimensional space $\mathcal{P}_{\leq 2}^{\perp}$ is injective by the null-space identity just proved, so $\widetilde{A}_{\beta,d}$ has infinite rank. The compact self-adjoint spectral theorem \citep[Theorem~2.3.2, p.~12]{Kostenko2019}, together with nonnegativity, therefore shows that $\widetilde{A}_{\beta,d}$ has infinitely many positive eigenvalues and that they converge to zero. This proves every assertion in part \textrm{(iii)}.

Let $\lambda > 0$. If $\widetilde{A}_{\beta,d}f = \lambda f$, then \eqref{eq:proofs.HZ.XstarX} and \eqref{eq:proofs.HZ.XXstar} give
\[
\begin{aligned}
A_{\beta,d}(\mathcal{X}_{\beta,d}f)
&= \mathcal{X}_{\beta,d}\mathcal{X}_{\beta,d}^*\mathcal{X}_{\beta,d}f
= \mathcal{X}_{\beta,d}\widetilde{A}_{\beta,d}f
= \lambda\mathcal{X}_{\beta,d}f, \\
\|\mathcal{X}_{\beta,d}f\|_{\beta}^2
&= \langle f,\mathcal{X}_{\beta,d}^*\mathcal{X}_{\beta,d}f\rangle_1
= \langle f,\widetilde{A}_{\beta,d}f\rangle_1
= \lambda\|f\|_1^2.
\end{aligned}
\]
Similarly, if $A_{\beta,d}g = \lambda g$, then
\[
\begin{aligned}
\widetilde{A}_{\beta,d}(\mathcal{X}_{\beta,d}^*g)
&= \mathcal{X}_{\beta,d}^*\mathcal{X}_{\beta,d}\mathcal{X}_{\beta,d}^*g
= \mathcal{X}_{\beta,d}^*A_{\beta,d}g
= \lambda\mathcal{X}_{\beta,d}^*g, \\
\|\mathcal{X}_{\beta,d}^*g\|_1^2
&= \langle g,\mathcal{X}_{\beta,d}\mathcal{X}_{\beta,d}^*g\rangle_{\beta}
= \langle g,A_{\beta,d}g\rangle_{\beta}
= \lambda\|g\|_{\beta}^2.
\end{aligned}
\]
On the two $\lambda$-eigenspaces, the maps $\lambda^{-1/2}\mathcal{X}_{\beta,d}$ and $\lambda^{-1/2}\mathcal{X}_{\beta,d}^*$ are mutually inverse because
\[
\frac{1}{\lambda}\mathcal{X}_{\beta,d}^*\mathcal{X}_{\beta,d}f = \frac{1}{\lambda}\widetilde{A}_{\beta,d}f = f, \qquad \frac{1}{\lambda}\mathcal{X}_{\beta,d}\mathcal{X}_{\beta,d}^*g = \frac{1}{\lambda}A_{\beta,d}g = g.
\]
Thus the two operators have the same positive eigenvalues with the same multiplicities. The complete list, including the absence of the degree-$2$ family when $d = 1$, the exceptional pole cases, and the addition of multiplicities when values from different families coincide, is therefore exactly the list in Theorem~\ref{thm:complete.spectrum}. This proves part \textrm{(ii)}. The integral definitions of $\mathcal{X}_{\beta,d}$ and $\mathcal{X}_{\beta,d}^*$ show that the two normalized maps are precisely \eqref{eq:results.HZ.eigenfunction.forward} and \eqref{eq:results.HZ.eigenfunction.backward}. This proves the eigenfunction assertions in part \textrm{(iv)}.

For clarity, we also construct the single unitary operator that gives the asserted equivalence, rather than only separate maps on the eigenspaces. Define $\mathcal{U}_{\beta,d}$ first on $\operatorname{Ran}(\widetilde{A}_{\beta,d}^{1/2})$ by
\begin{equation}\label{eq:proofs.HZ.unitary.definition}
\mathcal{U}_{\beta,d}(\widetilde{A}_{\beta,d}^{1/2}f) \leqdef \mathcal{X}_{\beta,d}f.
\end{equation}
This definition is independent of the choice of $f$. Indeed, if $\widetilde{A}_{\beta,d}^{1/2}f_1 = \widetilde{A}_{\beta,d}^{1/2}f_2$, then
\[
\|\mathcal{X}_{\beta,d}(f_1 - f_2)\|_{\beta}^2 = \|\widetilde{A}_{\beta,d}^{1/2}(f_1 - f_2)\|_1^2 = 0.
\]
Moreover,
\[
\|\mathcal{U}_{\beta,d}(\widetilde{A}_{\beta,d}^{1/2}f)\|_{\beta}^2 = \|\mathcal{X}_{\beta,d}f\|_{\beta}^2 = \|\widetilde{A}_{\beta,d}^{1/2}f\|_1^2,
\]
so $\mathcal{U}_{\beta,d}$ is an isometry on its initial domain. The square-root and adjoint range identities \citep[Theorems~3.3.7 and~3.4.3, pp.~73--75]{HsingEubank2015} give $\overline{\operatorname{Ran}(\widetilde{A}_{\beta,d}^{1/2})} = \ker(\widetilde{A}_{\beta,d})^{\perp} = \mathcal{P}_{\leq 2}^{\perp}$. By definition \eqref{eq:proofs.HZ.unitary.definition}, the image of $\operatorname{Ran}(\widetilde{A}_{\beta,d}^{1/2})$ is $\operatorname{Ran}(\mathcal{X}_{\beta,d})$. For $g \in L^2(\mu_{\beta})$, \eqref{eq:proofs.HZ.XXstar} gives
\[
\langle A_{\beta,d}g,g\rangle_{\beta} = \langle\mathcal{X}_{\beta,d}\mathcal{X}_{\beta,d}^*g,g\rangle_{\beta} = \|\mathcal{X}_{\beta,d}^*g\|_1^2.
\]
Proposition~\ref{prop:positivity.trace.injectivity} therefore implies $\ker(\mathcal{X}_{\beta,d}^*) = \ker(A_{\beta,d}) = \{0\}$. The adjoint range identity \citep[Theorem~3.3.7\textrm{(part 4)}, pp.~73--74]{HsingEubank2015} gives $\operatorname{Ran}(\mathcal{X}_{\beta,d})^{\perp} = \ker(\mathcal{X}_{\beta,d}^*)$, so the range is dense in $L^2(\mu_{\beta})$. Hence \eqref{eq:proofs.HZ.unitary.definition} extends uniquely to a unitary map from $\mathcal{P}_{\leq 2}^{\perp}$ onto $L^2(\mu_{\beta})$. For $f \in L^2(\mu_1)$, the commutation of $\widetilde{A}_{\beta,d}$ with its square root and the two factorizations give
\begin{align*}
A_{\beta,d}\mathcal{U}_{\beta,d}(\widetilde{A}_{\beta,d}^{1/2}f)
&= A_{\beta,d}\mathcal{X}_{\beta,d}f
= \mathcal{X}_{\beta,d}\mathcal{X}_{\beta,d}^*\mathcal{X}_{\beta,d}f
= \mathcal{X}_{\beta,d}\widetilde{A}_{\beta,d}f, \\
\mathcal{U}_{\beta,d}\widetilde{A}_{\beta,d}(\widetilde{A}_{\beta,d}^{1/2}f)
&= \mathcal{U}_{\beta,d}(\widetilde{A}_{\beta,d}^{1/2}\widetilde{A}_{\beta,d}f)
= \mathcal{X}_{\beta,d}\widetilde{A}_{\beta,d}f.
\end{align*}
Hence $A_{\beta,d}\mathcal{U}_{\beta,d} = \mathcal{U}_{\beta,d}\widetilde{A}_{\beta,d}$ on $\operatorname{Ran}(\widetilde{A}_{\beta,d}^{1/2})$, and continuity extends this identity to $\mathcal{P}_{\leq 2}^{\perp}$. Together with the eigenfunction correspondence proved above, this proves every assertion in part \textrm{(iv)}.

It remains to prove part \textrm{(v)}. The Hilbert--Schmidt trace identity \citep[Theorem~3.3.1\textrm{(iv)}, p.~27]{Kostenko2019} and the trace equality for the two products \citep[Corollary~3.4.3, p.~42]{Kostenko2019} give
\[
\begin{aligned}
\tr(\widetilde{A}_{\beta,d})
&= \tr(\mathcal{X}_{\beta,d}^*\mathcal{X}_{\beta,d})
= \|\mathcal{X}_{\beta,d}\|_2^2
= \tr(\mathcal{X}_{\beta,d}\mathcal{X}_{\beta,d}^*)
= \tr(A_{\beta,d}) \\
&= 1 - (1 + 2 \beta^2)^{-d/2} - d\beta^2(1 + 2 \beta^2)^{-d/2-1} - \frac{d(d + 2)}{2}\beta^4(1 + 2 \beta^2)^{-d/2-2},
\end{aligned}
\]
where the last equality is \eqref{eq:reduction.trace.formula}, so this proves \eqref{eq:results.HZ.trace}.

Let $(\lambda_j(\beta,d))_{j \in \N}$ be the common positive eigenvalues, repeated according to multiplicity. Then the Fredholm determinant product formula \citep[Eq.~3.4.25, p.~41]{Kostenko2019} yields, for every $z \in \C$,
\[
\det\nolimits_F(I - z\widetilde{A}_{\beta,d}) = \prod_{j=1}^{\infty}(1 - z\lambda_j(\beta,d)) = \det\nolimits_F(I - zA_{\beta,d}).
\]
The zero eigenspace in \eqref{eq:results.HZ.null.space} contributes only factors equal to one. Finally, Theorem~3.1 of \citet{HenzeZirkler1990} identifies the weights $\delta_k(\beta)$ with the positive eigenvalues of $\widetilde{A}_{\beta,d}$, and part \textrm{(ii)} identifies those eigenvalues with the list in Theorem~\ref{thm:complete.spectrum}. This proves every assertion in part \textrm{(v)} and completes the proof of Theorem~\ref{thm:Henze.Zirkler.spectrum}. \qed

\subsection{Proof of Corollary~\ref{cor:baseline.multiplicities}}\label{subsec:proof.baseline.multiplicities}

By Proposition~\ref{prop:Cartesian.Mercer.expansion}, the eigenspace of $B_{\beta,d}$ associated with $c_{\beta}q^N$ has an orthonormal basis $\{e_{\bb{\nu}} : \bb{\nu} \in \N_0^d,\ |\bb{\nu}| = N\}$. There are $\binom{N + d - 1}{N}$ such multi-indices, so this eigenvalue has that multiplicity. The angular degrees at total degree $N$ satisfy $\ell + 2k = N$ for some $k \in \N_0$. If $N$ is odd, the only affected degree is $\ell = 1$, of dimension $d$. If $N \geq 2$ is even, the affected degrees are $\ell = 0$ and $\ell = 2$, of total dimension
\[
h_{d,0} + h_{d,2} = 1 + \frac{(d - 1)(d + 2)}{2} = \binom{d + 1}{2}.
\]
For $N = 0,1,2$, all available sectors are affected. Subtracting these dimensions proves the three cases in the corollary. The formulas remain valid when $d = 1$, in which case every unchanged multiplicity is zero. This concludes the proof of Corollary~\ref{cor:baseline.multiplicities}. \qed

\subsection{Proof of Corollary~\ref{cor:exceptional.coincidences}}\label{subsec:proof.exceptional.coincidences}

By Theorem~\ref{thm:complete.spectrum}\textrm{(iv)}, $c_{\beta}\rho^m$ is a radial eigenvalue of $A_{\beta,d}$ if and only if \eqref{eq:results.radial.pole.equation} holds, that is, if and only if $F_m(\rho) = 0$, where $\rho = \rho(\beta)$ and
\[
F_m(\rho) \leqdef \sum_{k\neq m}(k - m)^2\frac{(d/2)_k}{k!}\frac{\rho^k}{\rho^k - \rho^m}, \qquad 0 < \rho < 1.
\]
Fix $m \in \N$. The series converges locally uniformly, so $F_m$ is continuous. For $0 < \rho \leq 1/2$ and $k = m + r > m$,
\[
\left|\frac{\rho^k}{\rho^k - \rho^m}\right| = \frac{\rho^r}{1 - \rho^r} \leq 2\rho^r.
\]
Consequently,
\[
\sum_{k=m+1}^{\infty}(k - m)^2\frac{(d/2)_k}{k!}\left|\frac{\rho^k}{\rho^k - \rho^m}\right| \leq 2\sum_{r=1}^{\infty}r^2\frac{(d/2)_{m+r}}{(m + r)!}\rho^r = O(\rho), \qquad \rho \downarrow 0,
\]
where the last power series has radius one. The sum over $k < m$ is finite, and each of its terms tends to $(m - k)^2(d/2)_k/k!$. Hence
\[
\lim_{\rho\downarrow 0}F_m(\rho) = \sum_{k=0}^{m-1}(m - k)^2\frac{(d/2)_k}{k!} > 0.
\]
For each fixed $k \neq m$,
\[
\lim_{\rho\uparrow 1}(1 - \rho)(k - m)^2\frac{(d/2)_k}{k!}\frac{\rho^k}{\rho^k - \rho^m} =
\begin{cases}
(m - k)(d/2)_k/k!, & k < m,\\
-(k - m)(d/2)_k/k!, & k > m.
\end{cases}
\]
The gamma-ratio asymptotic \citep[Eq.~5.11.12]{DLMF5} gives
\[
\frac{(k - m)(d/2)_k}{k!} = \frac{(k - m)}{\Gamma(d/2)}\frac{\Gamma(k + d/2)}{\Gamma(k + 1)} \sim \frac{k^{d/2}}{\Gamma(d/2)}, \qquad k \to \infty.
\]
Thus the positive sum over $k < m$ is finite, whereas $\sum_{k>m}(k - m)(d/2)_k/k!$ diverges. Choose a finite $M > m$ for which the sum of the absolute values of the negative coefficients with $m < k \leq M$ exceeds the sum of the positive coefficients with $k < m$.  Then all remaining terms are negative, and it follows that $F_m(\rho) < 0$ for $\rho$ sufficiently close to one. By the intermediate value theorem, there exists $\rho^* \in (0,1)$ such that $F_m(\rho^*) = 0$, and by the equivalence above, $c_{\beta}\rho^m$ is a radial eigenvalue of $A_{\beta,d}$ for the corresponding $\beta$. Taking $q = \sqrt{\rho^*}$ and $\beta^2 = q/(1 - q)^2$ produces the required smoothing parameter. When $m = 0$, every term in $F_0(\rho)$ is strictly negative, so $F_0$ never vanishes, and by the same equivalence $c_{\beta}$ is never a radial eigenvalue of $A_{\beta,d}$. This concludes the proof of Corollary~\ref{cor:exceptional.coincidences}. \qed

\subsection{Proof of Proposition~\ref{prop:eigenfunction.reconstruction}}\label{subsec:proof.eigenfunction.reconstruction}

Fix $\ell \in \{1,2\}$, an angular index $j$, and one of the eigenvalues $\lambda = \lambda_{\ell,m}$ from Theorem~\ref{thm:complete.spectrum}. Under the unitary identification of $\mathcal{V}_{\ell,j}$ with $\ell^2(\N_0)$ used in Proposition~\ref{prop:exact.block.reduction}, the affected operator is
\[
D_{\ell} - b_{\ell} \otimes b_{\ell}, \qquad b_{\ell} \leqdef (b_{\ell,k})_{k \in \N_0}, \qquad b_{\ell,k} \leqdef \{\Lambda_{\ell,k}\pi_{\ell,k}\}^{1/2}.
\]
Since $x_{\ell,m} \in (\rho^{m+1},\rho^m)$, the value $\lambda = c_{\beta}q^{\ell}x_{\ell,m}$ differs from every diagonal entry $\Lambda_{\ell,k} = c_{\beta}q^{\ell}\rho^k$. If $a = (a_k)_{k \in \N_0}$ is an eigenvector, then the $k$th coordinate of $(D_{\ell} - b_{\ell} \otimes b_{\ell})a = \lambda a$ is
\[
(\Lambda_{\ell,k} - \lambda)a_k = b_{\ell,k}\langle a,b_{\ell}\rangle.
\]
The scalar $\langle a,b_{\ell}\rangle$ cannot vanish. Indeed, if it vanished, then $(\Lambda_{\ell,k} - \lambda)a_k = 0$ for every $k$, and the fact that $\lambda \neq \Lambda_{\ell,k}$ would force $a_k = 0$ for every $k$, contrary to the choice of an eigenvector. Hence
\[
a_k = \langle a,b_{\ell}\rangle\frac{\{\Lambda_{\ell,k}\pi_{\ell,k}\}^{1/2}}{\Lambda_{\ell,k} - \lambda}, \qquad k \in \N_0,
\]
which gives \eqref{eq:reduction.rank.one.eigenfunction} up to the nonzero overall factor $\langle a,b_{\ell}\rangle$.

For completeness, the displayed coefficient sequence can also be verified directly. Let
\[
\widetilde{a}_k \leqdef \frac{b_{\ell,k}}{\Lambda_{\ell,k} - \lambda}.
\]
The denominators are bounded away from zero: they are nonzero for every $k$, and $\Lambda_{\ell,k} \to 0$ while $\lambda > 0$. Moreover,
\[
\sum_{k=0}^{\infty}b_{\ell,k}^2 = \sum_{k=0}^{\infty}\Lambda_{\ell,k}\pi_{\ell,k} \leq \Lambda_{\ell,0}\sum_{k=0}^{\infty}\pi_{\ell,k} = \Lambda_{\ell,0},
\]
so $(\widetilde{a}_k)_{k \in \N_0} \in \ell^2(\N_0)$. The scalar compression equation used in the proof of Theorem~\ref{thm:complete.spectrum} is $\sum_{k=0}^{\infty} \pi_{\ell,k}/(\Lambda_{\ell,k} - \lambda) = 0$. Consequently,
\[
\langle\widetilde{a},b_{\ell}\rangle = \sum_{k=0}^{\infty}\frac{\Lambda_{\ell,k}\pi_{\ell,k}}{\Lambda_{\ell,k} - \lambda} = \sum_{k=0}^{\infty}\pi_{\ell,k} + \lambda\sum_{k=0}^{\infty}\frac{\pi_{\ell,k}}{\Lambda_{\ell,k} - \lambda} = 1.
\]
It follows coordinatewise that
\[
\{(D_{\ell} - b_{\ell} \otimes b_{\ell})\widetilde{a}\}_k = \frac{\Lambda_{\ell,k}b_{\ell,k}}{\Lambda_{\ell,k} - \lambda} - b_{\ell,k} = \lambda\frac{b_{\ell,k}}{\Lambda_{\ell,k} - \lambda} = \lambda\widetilde{a}_k,
\]
which confirms directly that the series in \eqref{eq:reduction.rank.one.eigenfunction} is an eigenfunction.

Next consider a non-pole radial eigenvalue $\lambda = c_{\beta}x$. Use the notation from the proof of Theorem~\ref{thm:complete.spectrum}:
\[
D = \operatorname{diag}(d_0,d_1,\ldots), \qquad \bb{b}_k = \{d_k\pi_{0,k}\}^{1/2}\bb{v}_k, \qquad (W\bb{\gamma})_k = \bb{b}_k^{\top}\bb{\gamma}.
\]
Thus the radial sequence operator is $D - WW^*$ and $W = D^{1/2}U_0$. Let $\bb{\gamma} \neq \bb{0}$ be a null vector chosen as in the proposition and define
\[
a_k \leqdef \frac{\bb{b}_k^{\top}\bb{\gamma}}{d_k - \lambda}, \qquad k \in \N_0.
\]
Since $\lambda$ is not a pole, the numbers $d_k - \lambda$ are bounded away from zero. The isometry of $U_0$ gives
\[
\sum_{k=0}^{\infty}|\bb{b}_k^{\top}\bb{\gamma}|^2 = \|W\bb{\gamma}\|_{\ell^2}^2 = \|D^{1/2}U_0\bb{\gamma}\|_{\ell^2}^2 \leq d_0\|\bb{\gamma}\|^2,
\]
so $a \in \ell^2(\N_0)$. By \eqref{eq:proofs.radial.nonpole.matrix}, the null-vector condition in the proposition is equivalent to
\[
\left[I_2 - \sum_{k=0}^{\infty}\frac{\bb{b}_k\bb{b}_k^{\top}}{d_k - \lambda}\right]\bb{\gamma} = \bb{0}.
\]
Therefore
\[
W^*a = \sum_{k=0}^{\infty}\bb{b}_ka_k = \sum_{k=0}^{\infty}\frac{\bb{b}_k\bb{b}_k^{\top}}{d_k - \lambda}\bb{\gamma} = \bb{\gamma}.
\]
The definition of the coefficients also gives $(D - \lambda I)a = W\bb{\gamma}$. Combining the last two identities yields
\[
(D - WW^*)a = Da - W\bb{\gamma} = \lambda a.
\]
Moreover, $a$ is nonzero because $W^*a = \bb{\gamma} \neq \bb{0}$. Since $d_k = \Lambda_{0,k}$, these coefficients are exactly those in \eqref{eq:reduction.radial.nonpole.eigenfunction}.

Finally, suppose that the pole value $\lambda = d_m$ is a radial eigenvalue. Set $\bb{w}_m = (y_m,-1)^{\top}$ and
\[
S_m \leqdef I_2 - \sum_{k\neq m}\frac{\bb{b}_k\bb{b}_k^{\top}}{d_k - d_m}.
\]
This matrix is well defined. Indeed, the numbers $|d_k - d_m|$, $k \neq m$, are bounded away from zero, and \eqref{eq:proofs.radial.isometry.moments} gives
\[
\sum_{k=0}^{\infty}\|\bb{b}_k\|^2 = \sum_{k=0}^{\infty}d_k\pi_{0,k}(1 + y_k^2) \leq d_0\sum_{k=0}^{\infty}\pi_{0,k}(1 + y_k^2) = 2d_0.
\]
Because the pole criterion \eqref{eq:results.radial.pole.equation} holds, \eqref{eq:proofs.radial.pole.scalar} gives $\bb{w}_m^{\top}S_m\bb{w}_m = 0$. Hence $S_m\bb{w}_m$ is orthogonal to $\bb{w}_m$. On the other hand,
\[
\bb{w}_m^{\top}\bb{v}_m = (y_m,-1)\begin{pmatrix}1\\y_m\end{pmatrix} = 0,
\]
so the one-dimensional orthogonal complement of $\bb{w}_m$ is $\operatorname{span}\{\bb{v}_m\} = \operatorname{span}\{\bb{b}_m\}$. It follows that there is a unique scalar $a_m$ such that $S_m\bb{w}_m = \bb{b}_m a_m$. Taking the inner product with $\bb{b}_m$ gives
\[
a_m = \frac{\bb{b}_m^{\top}S_m\bb{w}_m}{\|\bb{b}_m\|^2},
\]
which is the separately defined coefficient in \eqref{eq:reduction.radial.pole.eigenfunction.coefficients}. For $k \neq m$, define
\[
a_k \leqdef \frac{\bb{b}_k^{\top}\bb{w}_m}{d_k - d_m}.
\]
The same bounded-denominator argument used above shows that $(a_k)_{k\neq m}$ is square-summable, and adjoining the single finite coordinate $a_m$ gives $a \in \ell^2(\N_0)$. By the definition of $S_m$ and the identity $S_m\bb{w}_m = \bb{b}_ma_m$,
\[
W^*a = \sum_{k\neq m}\frac{\bb{b}_k\bb{b}_k^{\top}}{d_k - d_m}\bb{w}_m + \bb{b}_ma_m = \bb{w}_m.
\]
In particular, $a \neq 0$. For $k \neq m$, the coefficient definition gives $(d_k - d_m)a_k = \bb{b}_k^{\top}\bb{w}_m$, while for $k = m$ both sides vanish because $\bb{b}_m^{\top}\bb{w}_m = 0$. Thus $(D - d_m I)a = W\bb{w}_m = WW^*a$, and hence
\[
(D - WW^*)a = d_m a.
\]
This proves the pole reconstruction formula.

Under each of the unitary identifications above, square-summability of the coefficient sequence is equivalent to convergence of the corresponding orthonormal series in $L^2(\mu_{\beta})$. Each constructed sequence is nonzero and may therefore be divided by its $\ell^2$ norm to obtain a normalized eigenfunction. For every sector of degree $\ell \geq 3$, Proposition~\ref{prop:exact.block.reduction} gives $A_{\beta,d} = B_{\beta,d}$ on that sector, while Proposition~\ref{prop:spherical.Mercer.expansion} gives $B_{\beta,d}\psi_{\ell,k,j} = \Lambda_{\ell,k}\psi_{\ell,k,j}$. Hence each stated function $\psi_{\ell,k,j}$ is an unchanged eigenfunction, as claimed. This completes the proof of Proposition~\ref{prop:eigenfunction.reconstruction}. \qed

%
\section{The common limiting null distribution}\label{sec:limiting.distribution}
%

Theorem~\ref{thm:Henze.Zirkler.spectrum} shows that the weights $\delta_k(\beta)$ in Theorem~3.1 of \citet{HenzeZirkler1990} and the weights $\lambda_j(\beta,d)$ in \eqref{eq:background.weighted.chi.square} are the same positive eigenvalues, repeated according to multiplicity. The finite-dimensional zero eigenspace of $\smash{\widetilde{A}_{\beta,d}}$ does not contribute to the limiting quadratic form. Tonelli's theorem and Proposition~\ref{prop:positivity.trace.injectivity} give
\[
\EE\left\{\int_{\R^d}Z(\bb{t})^2 \, \rd \mu_{\beta}(\bb{t})\right\} = \int_{\R^d}K(\bb{t},\bb{t}) \, \rd \mu_{\beta}(\bb{t}) = \tr(A_{\beta,d}) < \infty.
\]
Consequently, the jointly measurable version used in \eqref{eq:background.limiting.distribution} has paths in $L^2(\mu_{\beta})$ almost surely. After redefining this version to be zero on the exceptional null set, joint measurability and Fubini's theorem show that $\langle Z,f\rangle_{\beta}$ is measurable for every $f \in L^2(\mu_{\beta})$. Since $L^2(\mu_{\beta})$ is separable, \citet[Theorem~7.1.2, p.~177]{HsingEubank2015} implies that $Z$ is an $L^2(\mu_{\beta})$-valued random element. Its covariance operator is $A_{\beta,d}$, and the covariance trace identity \citep[Theorem~7.2.5, pp.~180--181]{HsingEubank2015} gives
\[
\EE\|Z\|_{\beta}^2 = \tr(A_{\beta,d}).
\]
Since $A_{\beta,d}$ is compact, self-adjoint, and injective, the spectral theorem shows that the orthonormal eigenfunctions corresponding to its positive eigenvalues form a complete orthonormal basis of $L^2(\mu_{\beta})$. Let $(e_j)_{j \in \N}$ be such a basis, indexed so that $A_{\beta,d}e_j = \lambda_j(\beta,d)e_j$, and let $N_j \leqdef \lambda_j(\beta,d)^{-1/2}\langle Z,e_j\rangle_{\beta}$. Because $Z$ is a centered Gaussian process with continuous covariance kernel, its $L^2(\mu_{\beta})$-valued version is a centered Gaussian random element. The Hilbert-space Karhunen--Lo\`{e}ve expansion \citep[Theorem~7.2.7, pp.~181--182]{HsingEubank2015} and Parseval's identity yield
\begin{equation}\label{eq:limiting.distribution.exact.series}
T_{\beta}(d) = \sum_{j=1}^{\infty}\lambda_j(\beta,d)N_j^2
\end{equation}
almost surely. Every finite subvector of $(N_j)_{j \in \N}$ is centered Gaussian with identity covariance matrix, so the $N_j$ are independent standard normal random variables. Moreover,
\[
\EE\left\{\sum_{j>m}\lambda_j(\beta,d)N_j^2\right\} = \sum_{j>m}\lambda_j(\beta,d) \to 0, \qquad m \to \infty,
\]
which proves convergence in $L^1$. The common complete list of positive weights and multiplicities is supplied by Theorems~\ref{thm:complete.spectrum} and~\ref{thm:Henze.Zirkler.spectrum}. In particular,
\begin{equation}\label{eq:limiting.distribution.mean}
\EE\{T_{\beta}(d)\}
= \sum_{j=1}^{\infty}\lambda_j(\beta,d)
= 1 - (1 + 2 \beta^2)^{-d/2} - d\beta^2(1 + 2 \beta^2)^{-d/2-1} - \frac{d(d + 2)}{2}\beta^4(1 + 2 \beta^2)^{-d/2-2}.
\end{equation}
This identity provides an exact check on any numerical root calculation. The degree-$1$ and degree-$2$ roots can be bracketed in the adjacent intervals $(\rho^{m+1},\rho^m)$. Radial roots should be obtained from the entire determinant \eqref{eq:reduction.radial.entire.expansion}, and not only from the meromorphic equation \eqref{eq:results.radial.nonpole.equation} because the latter omits the exceptional pole eigenvalues. Truncating \eqref{eq:limiting.distribution.exact.series} once the retained eigenvalues account for a desired proportion of the trace permits numerical approximation of asymptotic critical values and probabilities. This trace proportion controls the mean of the omitted nonnegative tail, but does not by itself give an exact bound on the error in the resulting critical values.

\begin{remark}[Extreme-smoothing limit]\label{rem:extreme.smoothing.limit}
Since the bandwidth of the Gaussian kernel density estimator in the equivalent formulation is $1/(\beta\sqrt{2})$, the limit $\beta \downarrow 0$ is the extreme-smoothing limit considered by \citet{Henze1997Extreme}. The formulas in Theorem~\ref{thm:complete.spectrum} recover the corresponding finite-rank spectrum. Indeed, $\frac{q}{\beta^2} \to 1$, $\rho = q^2$, and $c_{\beta} \to 1$. Hence
\[
\frac{c_{\beta}q^3}{\beta^6} \to 1.
\]
When $d \geq 2$, the value $c_{\beta}q^3 = \lambda_{3,0}^{(\mathrm{u})}$ has an unchanged spectral contribution of multiplicity
\[
h_{d,3} = \binom{d + 2}{3} - d = \frac{d(d - 1)(d + 4)}{6}.
\]
When $d = 1$, this unchanged family is absent and $h_{1,3} = 0$.

For every compact set $K \subset (1,\infty)$, separating the terms with $k = 0$ and $k = 1$ in the defining series gives, uniformly for $y \in K$,
\[
\begin{aligned}
\mathcal{Q}_{d/2+1,0}(\rho y)
&= \frac{1}{1 - \rho y} + \frac{d/2 + 1}{1 - y} + \sum_{k=2}^{\infty}\frac{(\frac{d}{2} + 1)_k}{k!}\frac{\rho^{k-1}}{\rho^{k-1} - y} \\
&= 1 + \frac{d/2 + 1}{1 - y} + O_K(\rho), \qquad \rho \downarrow 0,
\end{aligned}
\]
where the last estimate follows from the binomial series and the fact that the denominators in the remaining sum are uniformly bounded away from zero. For each $\varepsilon \in (0,d/2 + 1)$, the limiting function is negative at $(d + 4)/2 - \varepsilon$ and positive at $(d + 4)/2 + \varepsilon$. The uniform estimate and the strict increase of $\mathcal{Q}_{d/2+1,0}(\rho y)$ in $y$ therefore imply that its unique root $x_{1,0}/\rho$ lies between these two points for all sufficiently small $\rho$. Letting $\varepsilon$ decrease to zero gives
\[
\frac{x_{1,0}}{\rho} \to \frac{d + 4}{2}, \qquad \frac{\lambda_{1,0}}{\beta^6} \to \frac{d + 4}{2},
\]
and the degree-$1$ contribution has multiplicity $d$. Moreover, \eqref{eq:limiting.distribution.mean} yields
\[
\tr(A_{\beta,d}) = \frac{d(d + 2)(d + 4)}{6}\beta^6 + O(\beta^8).
\]
Let $R_{\beta}$ denote the contribution to \eqref{eq:limiting.distribution.exact.series} of all terms other than the $d$ copies of $\lambda_{1,0}$ and, when $d \geq 2$, the $h_{d,3}$ copies of $\lambda_{3,0}^{(\mathrm{u})}$. Since all eigenvalues are nonnegative, $R_{\beta} \geq 0$, and
\[
\frac{\EE(R_{\beta})}{\beta^6} = \frac{\tr(A_{\beta,d}) - d\lambda_{1,0} - h_{d,3}c_{\beta}q^3}{\beta^6} \to 0.
\]
Thus $R_{\beta}/\beta^6$ tends to zero in $L^1$ and hence in probability. By independence of the Gaussian coordinates in \eqref{eq:limiting.distribution.exact.series}, the two retained contributions are independent and have distributions $\lambda_{1,0}\chi_d^2$ and $c_{\beta}q^3\chi_{h_{d,3}}^2$, respectively, with the second contribution absent when $d = 1$. Slutsky's theorem now gives
\[
\frac{T_{\beta}(d)}{\beta^6} \xrightarrow{\mathrm{law}} \frac{d + 4}{2}\chi_d^2 + \chi_{d(d - 1)(d + 4)/6}^2, \qquad \beta \downarrow 0,
\]
where the two chi-square variables are independent and the second term is absent when $d = 1$. This is exactly the finite-rank null law obtained in \citet[Theorem~2.2]{Henze1997Extreme} for the extreme-smoothing statistic identified in Theorem~2.1 of that paper. There, $\beta \downarrow 0$ is taken first for fixed $n$, followed by $n \to \infty$, whereas the present result takes these limits in the opposite order. Thus the two iterated null limits agree.
\end{remark}

\smallskip

\begin{remark}[Vanishing-bandwidth limit]\label{rem:large_beta}
The opposite boundary regime $\beta \to \infty$ also admits a simple interpretation in terms of the spectral representation~\eqref{eq:limiting.distribution.exact.series}. Recall that the bandwidth in the equivalent Gaussian kernel density formulation is $1/(\beta\sqrt{2})$, so that $\beta \to \infty$ corresponds to a bandwidth tending to zero. Let $m_{\beta} \leqdef \EE\{T_{\beta}(d)\} = \tr(A_{\beta,d})$. It follows from~\eqref{eq:limiting.distribution.mean} that
\[
m_{\beta} = 1 - 2^{-d/2 - 3}(d + 2)(d + 4)\beta^{-d} + o(\beta^{-d}), \qquad \beta \to \infty.
\]

\noindent We next determine the asymptotic variance. Write
\[
K(\bb{s},\bb{t}) = G(\bb{s},\bb{t}) - R(\bb{s},\bb{t}),
\]
where
\[
G(\bb{s},\bb{t}) \leqdef \exp\left(-\frac{\|\bb{s} - \bb{t}\|^2}{2}\right), \qquad
R(\bb{s},\bb{t}) \leqdef \left\{1 + \bb{s}^{\top}\bb{t} + \frac{(\bb{s}^{\top}\bb{t})^2}{2}\right\}\exp\left(-\frac{\|\bb{s}\|^2 + \|\bb{t}\|^2}{2}\right).
\]
Since $A_{\beta,d}$ is self-adjoint and Hilbert--Schmidt,
\[
\tr(A_{\beta,d}^2) = \iint_{\R^d\times\R^d}K(\bb{s},\bb{t})^2 \, \rd \mu_{\beta}(\bb{s}) \rd \mu_{\beta}(\bb{t}).
\]
The contribution of $G^2$ is
\[
\iint_{\R^d\times\R^d}\exp\{-\|\bb{s} - \bb{t}\|^2\} \, \rd \mu_{\beta}(\bb{s}) \rd \mu_{\beta}(\bb{t}) = (1 + 4\beta^2)^{-d/2}.
\]
On the other hand,
\[
\begin{aligned}
\iint_{\R^d\times\R^d} |G(\bb{s},\bb{t})R(\bb{s},\bb{t})| \, \rd \mu_{\beta}(\bb{s}) \rd \mu_{\beta}(\bb{t}) &= O(\beta^{-2d}), \\
\iint_{\R^d\times\R^d}R(\bb{s},\bb{t})^2 \, \rd \mu_{\beta}(\bb{s}) \rd \mu_{\beta}(\bb{t}) &= O(\beta^{-2d}).
\end{aligned}
\]
Indeed, after insertion of the Gaussian densities, each of these integrals is $(2\pi\beta^2)^{-d}$ times an integral over $\R^{2d}$ of a fixed polynomial multiplied by an integrable Gaussian function; the remaining factor involving $\beta$ is bounded by one. Consequently, $\tr(A_{\beta,d}^2) = (1 + 4\beta^2)^{-d/2} + O(\beta^{-2d}) \sim 2^{-d}\beta^{-d}$, and hence
\[
\operatorname{Var}\{T_{\beta}(d)\} = 2\tr(A_{\beta,d}^2) \sim 2^{1 - d}\beta^{-d}.
\]

It remains to identify the limiting distribution. By~\eqref{eq:limiting.distribution.exact.series},
\[
T_{\beta}(d) - m_{\beta} = \sum_{j=1}^{\infty}\lambda_j(\beta,d)(N_j^2 - 1).
\]
The nonnegativity in Proposition~\ref{prop:positivity.trace.injectivity} and the finite-rank representation in Proposition~\ref{prop:exact.block.reduction} give $0 \leq A_{\beta,d} \leq B_{\beta,d}$ in the operator order. Hence, by Proposition~\ref{prop:Cartesian.Mercer.expansion},
\[
\max_{j\geq1}\lambda_j(\beta,d) \leq \|B_{\beta,d}\| = c_{\beta} = (1 - q)^d \sim \beta^{-d}.
\]
Together with $\sum_{j=1}^{\infty}\lambda_j(\beta,d)^2 = \tr(A_{\beta,d}^2) \sim 2^{-d}\beta^{-d}$, this yields
\[
\frac{\max_{j\geq1}\lambda_j(\beta,d)}{\left\{\sum_{j=1}^{\infty}\lambda_j(\beta,d)^2\right\}^{1/2}} \to 0.
\]
Since $\EE |N_1^2 - 1|^3 < \infty$, Lyapunov's condition follows from
\[
\frac{\sum_{j=1}^{\infty}\EE\left|\lambda_j(\beta,d)(N_j^2 - 1)\right|^3}{\left\{2\sum_{j=1}^{\infty}\lambda_j(\beta,d)^2\right\}^{3/2}} \leq C \, \frac{\max_{j\geq1}\lambda_j(\beta,d)}{\left\{\sum_{j=1}^{\infty}\lambda_j(\beta,d)^2\right\}^{1/2}} \to 0,
\]
where $C$ is a finite constant independent of $\beta$. Lyapunov's theorem, applied to the infinite series by $L^2$ truncation, therefore gives
\[
\frac{T_{\beta}(d) - m_{\beta}}{\{2\tr(A_{\beta,d}^2)\}^{1/2}} \xrightarrow{\mathrm{law}} \ N(0,1), \qquad \beta \to \infty.
\]
Since $\beta^{d/2}(m_{\beta} - 1) \to 0$, we obtain
\[
\beta^{d/2}\{T_{\beta}(d) - 1\} \ \xrightarrow{\mathrm{law}} \ N\bigl(0,2^{1 - d}\bigr), \qquad \beta \to \infty.
\]

This result complements Remark~\ref{rem:extreme.smoothing.limit}. As $\beta \downarrow 0$, only finitely many eigenvalues contribute on the relevant scale, and the limit is a finite weighted sum of independent chi-square variables. By contrast, as $\beta \to \infty$, no individual eigenvalue contributes appreciably to the total variance, and the accumulation of increasingly many small spectral contributions produces a Gaussian limit.

There is also a connection with the other extreme-smoothing result in \citet{Henze1997Extreme}. For fixed $n \geq d + 1$, suppose that $S_n$ is nonsingular and the observations are pairwise distinct. These conditions hold almost surely under the normal null hypothesis. Then Theorem~3.1 of that paper gives
\[
\beta^d\{T_{n,\beta} - 1\} \to n2^{-d/2} - 2\sum_{j=1}^n\exp\left(-\frac{\|\bb{Y}_j\|^2}{2}\right), \qquad \beta \to \infty.
\]
Thus the limiting statistic depends only on the squared Mahalanobis distances and is, in this sense, similar to Mardia's measure of multivariate kurtosis. Under the normal null hypothesis, Corollary~3.3 of \citet{Henze1997Extreme} shows that
\[
\sqrt{n}\left\{\frac{1}{n}\sum_{j=1}^n\exp\left(-\frac{\|\bb{Y}_j\|^2}{2}\right) - 2^{-d/2}\right\} \ \ \xrightarrow{\mathrm{law}} \ N(0,\sigma_d^2),
\]
where $\sigma_d^2 = 3^{-d/2} - 2^{-d} - d \, 2^{-(d + 3)}$. In contrast to the situation in Remark~\ref{rem:extreme.smoothing.limit}, however, the natural normalizations in the two iterated limits are different: the fixed-$n$ limit uses the scale $\beta^d$, whereas when the limit $n \to \infty$ is taken first, the fluctuations of the limiting quadratic form are of order $\beta^{-d/2}$.
\end{remark}

%
\section{Numerical evaluations}\label{sec:numerical.evaluations}
%

The numerical validation uses two complementary finite-dimensional approximations of the compact, self-adjoint integral BHEP covariance operator in \eqref{eq:definitions.BHEP.operator}. The Rayleigh--Ritz method follows \citet{EbnerJimenezGameroMilosevic2025}: the operator is projected onto the space of multivariate Hermite polynomials whose total degree does not exceed a prescribed cutoff, and its eigenvalues are approximated by those of the resulting symmetric Galerkin matrix. Whereas \citet{EbnerJimenezGameroMilosevic2025} represent this space using Cartesian tensor-product Hermite polynomials, the implementation here makes an additional change to an equivalent rotation-adapted basis consisting of spherical harmonics and generalized Laguerre polynomials (see \citet{DunklXu2014,DaiXu2013}). Since the BHEP kernel and the Gaussian weight are rotationally invariant, this representation decomposes the Galerkin matrix into smaller radial blocks indexed by the angular degree, without changing the Rayleigh--Ritz eigenvalues. The second numerical method is the so-called Nystr\"om method. Here, the same angular decomposition is used, but each remaining radial integral operator is discretized directly by generalized Gauss--Laguerre quadrature. Its eigenvalues are then approximated by those of the corresponding symmetrically weighted kernel matrix, as described by \citet{Bornemann2010}. This Gauss--Laguerre Nystr\"om approximation is denoted by GLN. For the GLN calculations below, where $\beta = 1$, the substitution $u = r^2 / 2$ gives the probability weight $u^{d/2 - 1}e^{-u} / \Gamma(d/2)$. The same radial nodes are used in every angular sector, with angular degrees $0 \leq \ell \leq 80$; only degrees $0$ and $1$ occur when $d = 1$. Agreement between the two approximations provides a numerical validation of the eigenvalues characterized in Theorem~\ref{thm:complete.spectrum}.

The results are summarized in Tables~\ref{tab:bhep-d1}--\ref{tab:bhep-d10}. Write $\nu_1 > \nu_2 > \cdots > 0$ for the distinct positive eigenvalue levels, and let $m(\nu_j)$ denote the full exact eigenspace multiplicity of $\nu_j$. Each of the ten largest levels is recorded once. The last four rows provide the cumulants of the limit distribution in \eqref{eq:limiting.distribution.exact.series}, calculated using $\kappa_r = 2^{r - 1}(r - 1)!\sum_{j=1}^{\infty} m(\nu_j)\nu_j^r$, $r = 1,\ldots,4$. Their numerical approximations use all eigenvalues of each finite numerical spectrum, counted with multiplicity. In the reference column, $\kappa_1,\kappa_2,\kappa_3$ are the closed forms in Theorem~2.3 of \citet{HenzeWagner1997}. The reference fourth cumulant is calculated from the affected-sector matrices in \eqref{eq:reduction.rank.one.matrix} and \eqref{eq:reduction.radial.matrix}, compressed to $k = 0,\ldots,39$, together with the exact unchanged-sector fourth-power sum
\[
\sum_{\ell\geq3}\sum_{k=0}^{\infty}h_{d,\ell}\Lambda_{\ell,k}^4 = c_{\beta}^4\left\{(1 - q^4)^{-d} - \frac{1 + dq^4 + h_{d,2}q^8}{1 - \rho^4}\right\}.
\]
The sum of these fourth-power contributions, including the angular multiplicities, is multiplied by $48$ to obtain $\kappa_4$. Increasing the reference compression to $k = 0,\ldots,49$, comparing GLN angular cutoffs $60$ and $80$, and increasing each stated GLN quadrature order by four changed the corresponding displayed eigenvalues and cumulants by less than $10^{-9}$. Both the displayed eigenvalues and cumulants agree with the reference values with a maximum absolute error of at most $10^{-6}$. For the unnormalized $w_{\gamma}$ convention in \citet{EbnerJimenezGameroMilosevic2025}, with $\gamma = 1/(2\beta^2)$, multiply $\kappa_r$ by $(2\pi\beta^2)^{rd/2}$. For $d = 1,2,3$ and $r = 1,2,3$, this gives the corresponding reference values in Tables~12--14 of that article; the present choice $\beta = 1$ corresponds to $\gamma = 1/2$.

\begin{table}[htbp]
  \centering
  \small
  \caption{The 10 largest distinct BHEP eigenvalues, their multiplicities, and the first four cumulants of the limiting distribution for $d = 1$ and $\beta = 1$. Rayleigh--Ritz uses a total-degree cutoff of 38; GLN uses 24 radial nodes per angular block. Reference cumulants 1--3 use Henze--Wagner's closed forms; all numerical cumulants and the reference fourth cumulant use the spectral truncations described in the text, not only the displayed eigenvalues.}\label{tab:bhep-d1}
  \begin{tabular*}{\linewidth}{@{\extracolsep{\fill}}rrccc@{}}
    \toprule
    Quantity & Multiplicity $m(\nu)$ & Reference
    & Rayleigh--Ritz & GLN \\
    \midrule
    $\nu_{1}$  & $1$ & $7.4274839\times 10^{-2}$
      & $7.4274839\times 10^{-2}$ & $7.4274839\times 10^{-2}$ \\
    $\nu_{2}$  & $1$ & $4.4810433\times 10^{-2}$
      & $4.4810433\times 10^{-2}$ & $4.4810433\times 10^{-2}$ \\
    $\nu_{3}$  & $1$ & $8.4190709\times 10^{-3}$
      & $8.4190709\times 10^{-3}$ & $8.4190710\times 10^{-3}$ \\
    $\nu_{4}$  & $1$ & $4.5868443\times 10^{-3}$
      & $4.5868443\times 10^{-3}$ & $4.5868443\times 10^{-3}$ \\
    $\nu_{5}$  & $1$ & $1.0799752\times 10^{-3}$
      & $1.0799751\times 10^{-3}$ & $1.0799752\times 10^{-3}$ \\
    $\nu_{6}$  & $1$ & $5.5193868\times 10^{-4}$
      & $5.5193850\times 10^{-4}$ & $5.5193867\times 10^{-4}$ \\
    $\nu_{7}$  & $1$ & $1.4573879\times 10^{-4}$
      & $1.4573840\times 10^{-4}$ & $1.4573879\times 10^{-4}$ \\
    $\nu_{8}$  & $1$ & $7.1210966\times 10^{-5}$
      & $7.1210395\times 10^{-5}$ & $7.1210932\times 10^{-5}$ \\
    $\nu_{9}$  & $1$ & $2.0182080\times 10^{-5}$
      & $2.0180819\times 10^{-5}$ & $2.0182115\times 10^{-5}$ \\
    $\nu_{10}$ & $1$ & $9.5383947\times 10^{-6}$
      & $9.5370573\times 10^{-6}$ & $9.5382962\times 10^{-6}$ \\
    \midrule
    $\kappa_{1}$ & -- & $1.3397460\times 10^{-1}$
      & $1.3397458\times 10^{-1}$ & $1.3397460\times 10^{-1}$ \\
    $\kappa_{2}$ & -- & $1.5236289\times 10^{-2}$
      & $1.5236289\times 10^{-2}$ & $1.5236289\times 10^{-2}$ \\
    $\kappa_{3}$ & -- & $4.0034300\times 10^{-3}$
      & $4.0034300\times 10^{-3}$ & $4.0034300\times 10^{-3}$ \\
    $\kappa_{4}$ & -- & $1.6546551\times 10^{-3}$
      & $1.6546551\times 10^{-3}$ & $1.6546551\times 10^{-3}$ \\
    \bottomrule
  \end{tabular*}
\end{table}

\begin{table}[htbp]
  \centering
  \small
  \caption{The 10 largest distinct BHEP eigenvalues, their multiplicities, and the first four cumulants of the limiting distribution for $d = 2$ and $\beta = 1$. Rayleigh--Ritz uses a total-degree cutoff of 42; GLN uses 20 radial nodes per angular block. Reference cumulants 1--3 use Henze--Wagner's closed forms; all numerical cumulants and the reference fourth cumulant use the spectral truncations described in the text, not only the displayed eigenvalues.}\label{tab:bhep-d2}
  \begin{tabular*}{\linewidth}{@{\extracolsep{\fill}}rrccc@{}}
    \toprule
    Quantity & Multiplicity $m(\nu)$ & Reference
    & Rayleigh--Ritz & GLN \\
    \midrule
    $\nu_{1}$  & $2$ & $5.3010377\times 10^{-2}$
      & $5.3010377\times 10^{-2}$ & $5.3010377\times 10^{-2}$ \\
    $\nu_{2}$  & $1$ & $3.5858544\times 10^{-2}$
      & $3.5858544\times 10^{-2}$ & $3.5858544\times 10^{-2}$ \\
    $\nu_{3}$  & $2$ & $2.5124253\times 10^{-2}$
      & $2.5124253\times 10^{-2}$ & $2.5124252\times 10^{-2}$ \\
    $\nu_{4}$  & $2$ & $2.1286236\times 10^{-2}$
      & $2.1286236\times 10^{-2}$ & $2.1286236\times 10^{-2}$ \\
    $\nu_{5}$  & $2$ & $8.1306188\times 10^{-3}$
      & $8.1306188\times 10^{-3}$ & $8.1306188\times 10^{-3}$ \\
    $\nu_{6}$  & $2$ & $5.8500403\times 10^{-3}$
      & $5.8500403\times 10^{-3}$ & $5.8500404\times 10^{-3}$ \\
    $\nu_{7}$  & $1$ & $3.4833668\times 10^{-3}$
      & $3.4833668\times 10^{-3}$ & $3.4833661\times 10^{-3}$ \\
    $\nu_{8}$  & $4$ & $3.1056200\times 10^{-3}$
      & $3.1056200\times 10^{-3}$ & $3.1056200\times 10^{-3}$ \\
    $\nu_{9}$  & $2$ & $2.7050076\times 10^{-3}$
      & $2.7050076\times 10^{-3}$ & $2.7050072\times 10^{-3}$ \\
    $\nu_{10}$ & $4$ & $1.1862413\times 10^{-3}$
      & $1.1862413\times 10^{-3}$ & $1.1862413\times 10^{-3}$ \\
    \midrule
    $\kappa_{1}$ & -- & $2.9629630\times 10^{-1}$
      & $2.9629627\times 10^{-1}$ & $2.9629630\times 10^{-1}$ \\
    $\kappa_{2}$ & -- & $1.8698560\times 10^{-2}$
      & $1.8698560\times 10^{-2}$ & $1.8698560\times 10^{-2}$ \\
    $\kappa_{3}$ & -- & $3.1738431\times 10^{-3}$
      & $3.1738431\times 10^{-3}$ & $3.1738431\times 10^{-3}$ \\
    $\kappa_{4}$ & -- & $8.9596411\times 10^{-4}$
      & $8.9596411\times 10^{-4}$ & $8.9596411\times 10^{-4}$ \\
    \bottomrule
  \end{tabular*}
\end{table}

\begin{table}[htbp]
  \centering
  \small
  \caption{The 10 largest distinct BHEP eigenvalues, their multiplicities, and the first four cumulants of the limiting distribution for $d = 3$ and $\beta = 1$. Rayleigh--Ritz uses a total-degree cutoff of 46; GLN uses 24 radial nodes per angular block. Reference cumulants 1--3 use Henze--Wagner's closed forms; all numerical cumulants and the reference fourth cumulant use the spectral truncations described in the text, not only the displayed eigenvalues.}\label{tab:bhep-d3}
  \begin{tabular*}{\linewidth}{@{\extracolsep{\fill}}rrccc@{}}
    \toprule
    Quantity & Multiplicity $m(\nu)$ & Reference
    & Rayleigh--Ritz & GLN \\
    \midrule
    $\nu_{1}$  & $3$  & $3.6849251\times 10^{-2}$
      & $3.6849251\times 10^{-2}$ & $3.6849251\times 10^{-2}$ \\
    $\nu_{2}$  & $1$  & $2.7569203\times 10^{-2}$
      & $2.7569203\times 10^{-2}$ & $2.7569203\times 10^{-2}$ \\
    $\nu_{3}$  & $5$  & $1.6878810\times 10^{-2}$
      & $1.6878810\times 10^{-2}$ & $1.6878810\times 10^{-2}$ \\
    $\nu_{4}$  & $7$  & $1.3155617\times 10^{-2}$
      & $1.3155617\times 10^{-2}$ & $1.3155617\times 10^{-2}$ \\
    $\nu_{5}$  & $9$  & $5.0249987\times 10^{-3}$
      & $5.0249987\times 10^{-3}$ & $5.0249987\times 10^{-3}$ \\
    $\nu_{6}$  & $3$  & $4.0023224\times 10^{-3}$
      & $4.0023224\times 10^{-3}$ & $4.0023224\times 10^{-3}$ \\
    $\nu_{7}$  & $1$  & $2.5811913\times 10^{-3}$
      & $2.5811912\times 10^{-3}$ & $2.5811912\times 10^{-3}$ \\
    $\nu_{8}$  & $18$ & $1.9193787\times 10^{-3}$
      & $1.9193787\times 10^{-3}$ & $1.9193787\times 10^{-3}$ \\
    $\nu_{9}$  & $5$  & $1.8102684\times 10^{-3}$
      & $1.8102684\times 10^{-3}$ & $1.8102684\times 10^{-3}$ \\
    $\nu_{10}$ & $22$ & $7.3313744\times 10^{-4}$
      & $7.3313744\times 10^{-4}$ & $7.3313744\times 10^{-4}$ \\
    \midrule
    $\kappa_{1}$ & -- & $4.5472475\times 10^{-1}$
      & $4.5472472\times 10^{-1}$ & $4.5472475\times 10^{-1}$ \\
    $\kappa_{2}$ & -- & $1.5700689\times 10^{-2}$
      & $1.5700689\times 10^{-2}$ & $1.5700689\times 10^{-2}$ \\
    $\kappa_{3}$ & -- & $1.7005047\times 10^{-3}$
      & $1.7005047\times 10^{-3}$ & $1.7005047\times 10^{-3}$ \\
    $\kappa_{4}$ & -- & $3.2310998\times 10^{-4}$
      & $3.2310998\times 10^{-4}$ & $3.2310998\times 10^{-4}$ \\
    \bottomrule
  \end{tabular*}
\end{table}

\begin{table}[htbp]
  \centering
  \small
  \caption{The 10 largest distinct BHEP eigenvalues, their multiplicities, and the first four cumulants of the limiting distribution for $d = 5$ and $\beta = 1$. Rayleigh--Ritz uses a total-degree cutoff of 52; GLN uses 24 radial nodes per angular block. Reference cumulants 1--3 use Henze--Wagner's closed forms; all numerical cumulants and the reference fourth cumulant use the spectral truncations described in the text, not only the displayed eigenvalues.}\label{tab:bhep-d5}
  \begin{tabular*}{\linewidth}{@{\extracolsep{\fill}}rrccc@{}}
    \toprule
    Quantity & Multiplicity $m(\nu)$ & Reference
    & Rayleigh--Ritz & GLN \\
    \midrule
    $\nu_{1}$  & $5$   & $1.6878810\times 10^{-2}$
      & $1.6878810\times 10^{-2}$ & $1.6878810\times 10^{-2}$ \\
    $\nu_{2}$  & $1$   & $1.4967635\times 10^{-2}$
      & $1.4967635\times 10^{-2}$ & $1.4967635\times 10^{-2}$ \\
    $\nu_{3}$  & $14$  & $7.3735581\times 10^{-3}$
      & $7.3735581\times 10^{-3}$ & $7.3735581\times 10^{-3}$ \\
    $\nu_{4}$  & $30$  & $5.0249987\times 10^{-3}$
      & $5.0249987\times 10^{-3}$ & $5.0249987\times 10^{-3}$ \\
    $\nu_{5}$  & $55$  & $1.9193787\times 10^{-3}$
      & $1.9193787\times 10^{-3}$ & $1.9193787\times 10^{-3}$ \\
    $\nu_{6}$  & $5$   & $1.8102684\times 10^{-3}$
      & $1.8102684\times 10^{-3}$ & $1.8102684\times 10^{-3}$ \\
    $\nu_{7}$  & $1$   & $1.3410459\times 10^{-3}$
      & $1.3410459\times 10^{-3}$ & $1.3410459\times 10^{-3}$ \\
    $\nu_{8}$  & $14$  & $7.9220801\times 10^{-4}$
      & $7.9220801\times 10^{-4}$ & $7.9220801\times 10^{-4}$ \\
    $\nu_{9}$  & $121$ & $7.3313744\times 10^{-4}$
      & $7.3313744\times 10^{-4}$ & $7.3313744\times 10^{-4}$ \\
    $\nu_{10}$ & $195$ & $2.8003358\times 10^{-4}$
      & $2.8003358\times 10^{-4}$ & $2.8003358\times 10^{-4}$ \\
    \midrule
    $\kappa_{1}$ & -- & $7.0419708\times 10^{-1}$
      & $7.0419706\times 10^{-1}$ & $7.0419708\times 10^{-1}$ \\
    $\kappa_{2}$ & -- & $6.9644860\times 10^{-3}$
      & $6.9644860\times 10^{-3}$ & $6.9644860\times 10^{-3}$ \\
    $\kappa_{3}$ & -- & $2.9836842\times 10^{-4}$
      & $2.9836842\times 10^{-4}$ & $2.9836842\times 10^{-4}$ \\
    $\kappa_{4}$ & -- & $2.4833780\times 10^{-5}$
      & $2.4833780\times 10^{-5}$ & $2.4833780\times 10^{-5}$ \\
    \bottomrule
  \end{tabular*}
\end{table}

\begin{table}[htbp]
  \centering
  \small
  \caption{The 10 largest distinct BHEP eigenvalues, their multiplicities, and the first four cumulants of the limiting distribution for $d = 10$ and $\beta = 1$. Rayleigh--Ritz uses a total-degree cutoff of 64; GLN uses 24 radial nodes per angular block. Reference cumulants 1--3 use Henze--Wagner's closed forms; all numerical cumulants and the reference fourth cumulant use the spectral truncations described in the text, not only the displayed eigenvalues.}\label{tab:bhep-d10}
  \begin{tabular*}{\linewidth}{@{\extracolsep{\fill}}rrccc@{}}
    \toprule
    Quantity & Multiplicity $m(\nu)$ & Reference
    & Rayleigh--Ritz & GLN \\
    \midrule
    $\nu_{1}$  & $1$    & $2.4235628\times 10^{-3}$
      & $2.4235628\times 10^{-3}$ & $2.4235628\times 10^{-3}$ \\
    $\nu_{2}$  & $10$   & $2.0144217\times 10^{-3}$
      & $2.0144217\times 10^{-3}$ & $2.0144217\times 10^{-3}$ \\
    $\nu_{3}$  & $54$   & $8.2743867\times 10^{-4}$
      & $8.2743867\times 10^{-4}$ & $8.2743867\times 10^{-4}$ \\
    $\nu_{4}$  & $210$  & $4.5310385\times 10^{-4}$
      & $4.5310385\times 10^{-4}$ & $4.5310385\times 10^{-4}$ \\
    $\nu_{5}$  & $10$   & $2.1983766\times 10^{-4}$
      & $2.1983766\times 10^{-4}$ & $2.1983766\times 10^{-4}$ \\
    $\nu_{6}$  & $1$    & $2.1203357\times 10^{-4}$
      & $2.1203357\times 10^{-4}$ & $2.1203356\times 10^{-4}$ \\
    $\nu_{7}$  & $660$  & $1.7307027\times 10^{-4}$
      & $1.7307027\times 10^{-4}$ & $1.7307027\times 10^{-4}$ \\
    $\nu_{8}$  & $54$   & $9.1635737\times 10^{-5}$
      & $9.1635737\times 10^{-5}$ & $9.1635734\times 10^{-5}$ \\
    $\nu_{9}$  & $1992$ & $6.6106961\times 10^{-5}$
      & $6.6106961\times 10^{-5}$ & $6.6106961\times 10^{-5}$ \\
    $\nu_{10}$ & $10$   & $2.6293948\times 10^{-5}$
      & $2.6293948\times 10^{-5}$ & $2.6293954\times 10^{-5}$ \\
    \midrule
    $\kappa_{1}$ & -- & $9.5473251\times 10^{-1}$
      & $9.5473247\times 10^{-1}$ & $9.5473251\times 10^{-1}$ \\
    $\kappa_{2}$ & -- & $3.2137769\times 10^{-4}$
      & $3.2137769\times 10^{-4}$ & $3.2137769\times 10^{-4}$ \\
    $\kappa_{3}$ & -- & $1.2028050\times 10^{-6}$
      & $1.2028050\times 10^{-6}$ & $1.2028050\times 10^{-6}$ \\
    $\kappa_{4}$ & -- & $1.1231540\times 10^{-8}$
      & $1.1231540\times 10^{-8}$ & $1.1231540\times 10^{-8}$ \\
    \bottomrule
  \end{tabular*}
\end{table}

\FloatBarrier

The numerical results provide a strong consistency check for the spectral description in Theorem~\ref{thm:complete.spectrum}. The reference values agree closely with two independent finite-dimensional approximations, both for the leading eigenvalues and for the first four cumulants of the limiting distribution. Using the multiplicities supplied by the spherical-harmonic decomposition, the two numerical methods reproduce the predicted spectral levels and cumulants, including in dimensions where some unchanged levels have large multiplicities. This agreement persists as the dimension increases, even though the multiplicities of some unchanged spectral levels become large. Thus, the numerical calculations support the predicted spectral levels and, using the analytically established multiplicities, reproduce the reference cumulants.

\begin{figure}[t!]
\centering
\includegraphics[width=10cm]{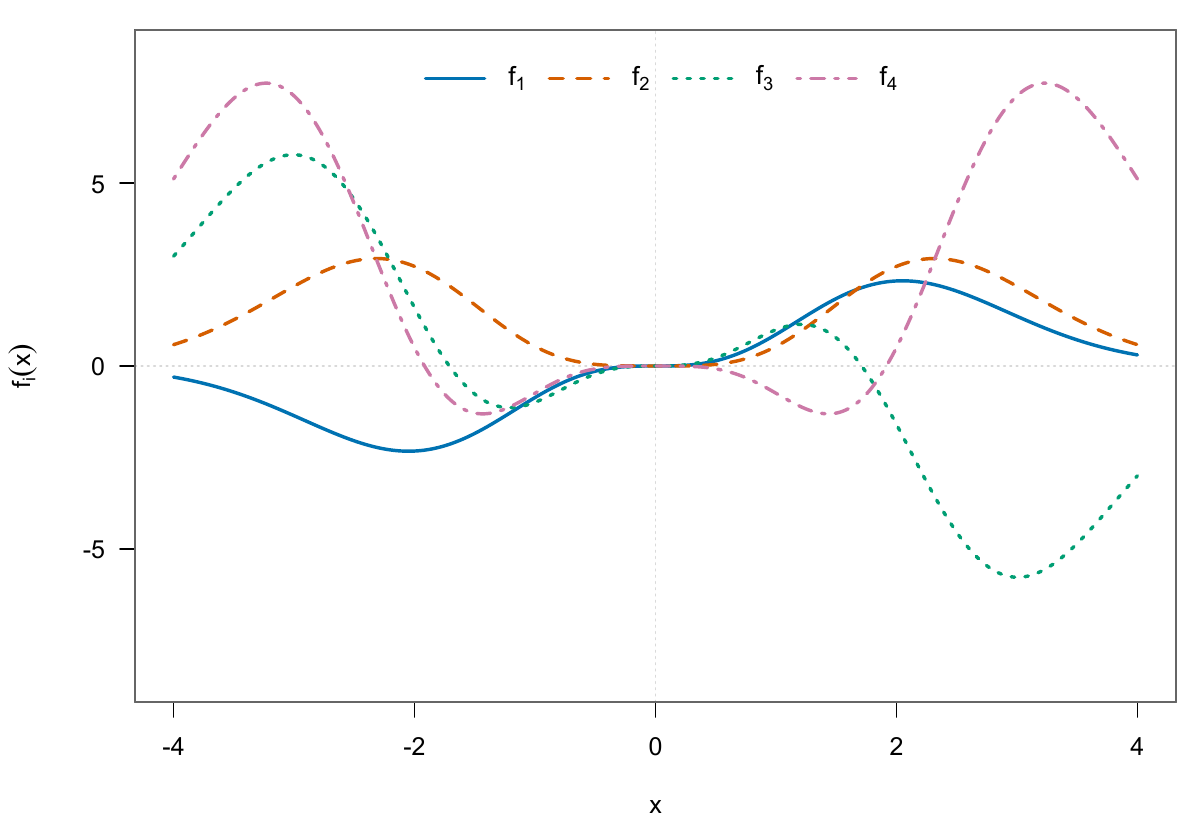}
\vspace{-3mm}
\caption{Eigenfunctions $f_1,\ldots,f_4$ of the BHEP covariance operator for $d=1$ and $\beta=1$, corresponding to the four largest eigenvalues in decreasing order.}\label{fig:eig_d1}
\end{figure}

\begin{figure}[H]
\centering
\includegraphics[width=15cm]{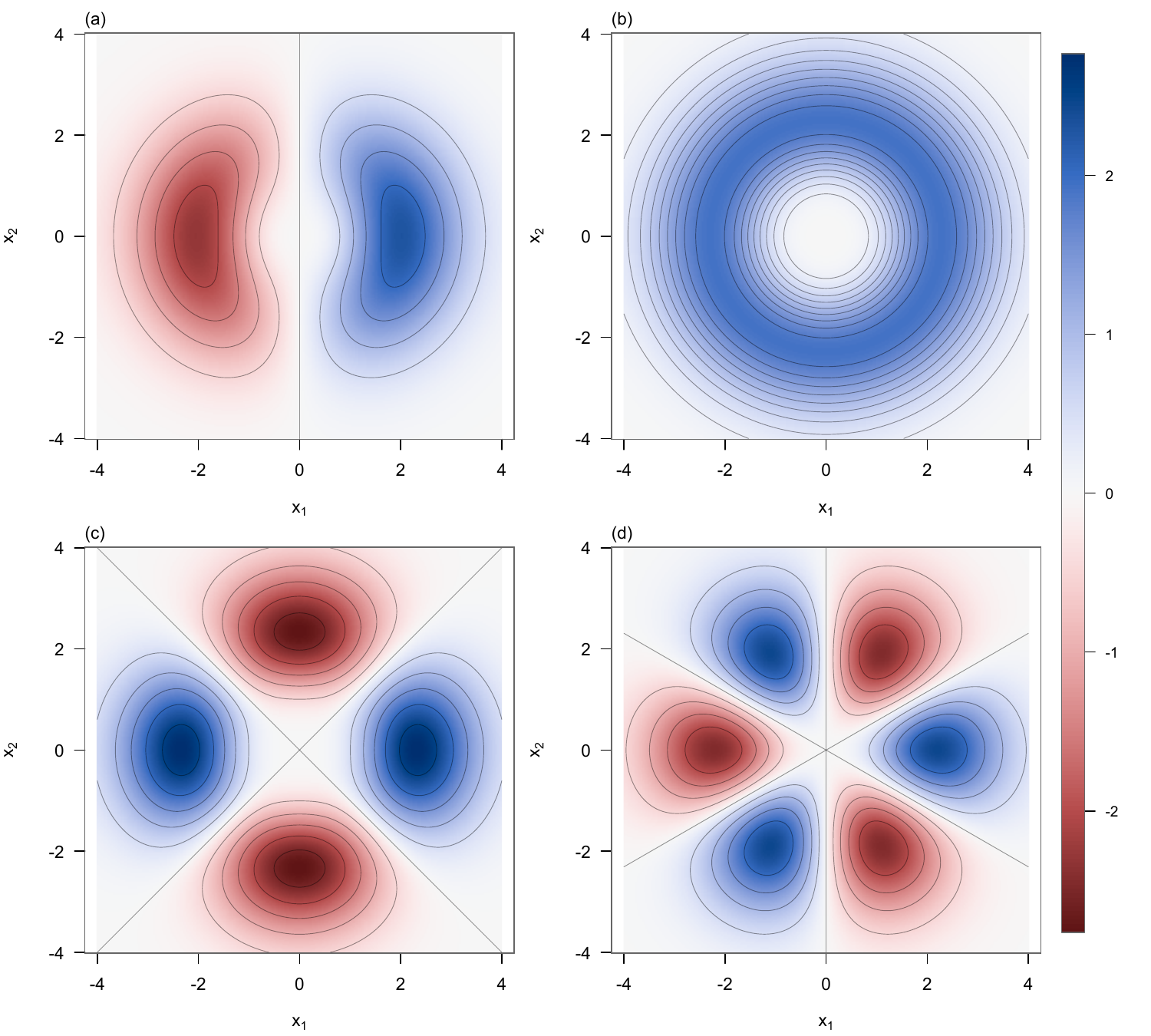}
\vspace{-3mm}
\caption{Eigenfunctions of the BHEP covariance operator for $d=2$ and $\beta=1$. Panels (a)--(d)  display one normalized real representative for each of the four largest distinct eigenvalues, in decreasing order. Their multiplicities are $2,1,2,2$, respectively, and the angular representative $j=1$ is selected.}\label{fig:eig_d2}
\end{figure}

Finally, we provide plots of eigenfunctions corresponding to the four largest distinct eigenvalues in decreasing order (the sign of each eigenfunction is arbitrary) for the dimensions $d=1$ (Figure \ref{fig:eig_d1}) and $d=2$ (Figure \ref{fig:eig_d2}). In each dimension, the eigenfunctions have unit norm with respect to the corresponding standard Gaussian measure.

\section*{Statement of AI use}
\addcontentsline{toc}{section}{Statement of AI use}

ChatGPT 5.6 Sol assisted in the initial discovery of the mathematical arguments and assisted with the numerical evaluations in Section~\ref{sec:numerical.evaluations}. The authors completely reworked and independently validated all proofs and assume full responsibility for their accuracy and rigor.

\section*{Funding}
\addcontentsline{toc}{section}{Funding}

B.\ Ebner is funded by the Deutsche Forschungsgemeinschaft (DFG, German Research Foundation) through grant 541565572. F.\ Ouimet is supported by the Natural Sciences and Engineering Research Council of Canada (NSERC) through Discovery Grant RGPIN-2026-04471 and Discovery Launch Supplement DGECR-2026-00449.



\addcontentsline{toc}{section}{References}
\setlength{\bibsep}{0pt plus 0ex}

\bibliographystyle{plainnat}
\bibliography{bib}

\end{document}